\documentclass{amsart}
\usepackage{graphicx}
\usepackage[all]{xy}
\usepackage{xy}

\usepackage{amsmath,amssymb}
\usepackage{amsfonts,amscd,latexsym}
\usepackage{euscript}

\newcommand{\labell}[1] {\label{#1}}

\newcommand{\sapprox}{\mbox{\tiny$\approx$}}
\newcommand{\Hsm}{\mbox{\it\tiny H}}

\newcommand{\bx}{{\bf x}}
\newcommand{\by}{{\bf y}}

\newcommand{\Mm}{{\mathcal M}}
\newcommand{\Nn}{{\mathcal N}}

\newcommand{\oo}{{\mathfrak o}}
\newcommand{\Qq}{{\mathcal Q}}
\newcommand{\Xx}{{\mathcal X}}

\newcommand{\Zz}{{\mathcal Z}}
\newcommand{\Vv}{{\mathcal V}}
\newcommand{\Uu}{{\mathcal U}}
\newcommand{\Mor}{{\rm Mor}}
\newcommand{\SSS}{{\smallskip}}

\newcommand{\Sect}{{\rm Sect}}
\newcommand{\supp}{{\rm supp}}

\newcommand{\HHat}{\widehat}
\newcommand{\GGL}{{\mathcal {GL}}}

\newcommand{\pr}{{\rm pr}}

  \newcommand{\Ss}{{{\mathcal S}}}

 \newcommand{\Bb}{{{\mathcal B}}}
 
\newcommand{\id}{{\rm id}}
\newcommand{\MS}{{\medskip}}

\newcommand{\N}{{\mathbb N}}

\newtheorem{theorem}{Theorem}[section]
\newtheorem{thm}[theorem]{Theorem}

\newtheorem{lemma}[theorem]{Lemma}

\newtheorem{prop}[theorem]{Proposition}
\newtheorem{defn}[theorem]{Definition}
\newtheorem{example}[theorem]{Example}
\newtheorem{rmk}[theorem]{Remark}

\numberwithin{figure}{section}
\numberwithin{equation}{section}
\numberwithin{table}{section}
\newcommand{\p}{{\partial}}
\newcommand{\al}{{\alpha}}

\newcommand{\be}{{\beta}}

\newcommand{\Om}{{\Omega}}

\newcommand{\eps}{{\varepsilon}}
\newcommand{\de}{{\delta}}

\newcommand{\ga}{{\gamma}}

\newcommand{\io}{{\iota}}
\newcommand{\ka}{{\kappa}}
\newcommand{\la}{{\lambda}}
\newcommand{\La}{{\Lambda}}
\newcommand{\si}{{\sigma}}
\newcommand{\Si}{{\Sigma}}

\newcommand{\Ll}{{\mathcal L}}
\newcommand{\Aa}{{\mathcal A}}
\newcommand{\Ii}{{\mathcal I}}
\newcommand{\Yy}{{\mathcal Y}}

\newcommand{\ov}{\overline}
\newcommand{\QED}{\hfill$\Box$\medskip}

\newcommand{\Cc}{{\mathcal C}}

\newcommand{\Tt}{{\mathcal T}}
\newcommand{\Ww}{{\mathcal W}}
\newcommand{\Ee}{{\mathcal E}}

\newcommand\im{\operatorname{Im}}
\newcommand\Int{\operatorname{Int}}

\newcommand\Diff{\operatorname{Diff}}

\newcommand{\Q}{{\mathbb Q}}

\newcommand{\R}{{\mathbb R}}
\newcommand{\Z}{{\mathbb Z}}

\newcommand{\un}{\underline}

\newcommand{\NI}{{\noindent}}
\begin{document}
 \title{Groupoids, branched manifolds and multisections}
 
 \author{Dusa McDuff}
\address{Department of Mathematics,
 Stony Brook University, Stony Brook, 
NY 11794-3651, USA}
\email{dusa@math.sunysb.edu}
\urladdr{http://www.math.sunysb.edu/~{}dusa}
\thanks{Partly supported by the NSF grant DMS 0305939.}
\date{18 September 2005, revised January 29 2006, July 16 2006 and December 23, 2006}
 
 \begin{abstract}
Cieliebak, Mundet i Riera and Salamon recently formulated a
 definition of branched submanifold of Euclidean space in connection 
with their discussion of multivalued sections and the Euler class.  
This note proposes an intrinsic definition of a weighted branched manifold $\un Z$
that  is obtained from  the usual definition 
of oriented  orbifold groupoid  by relaxing the properness condition
and adding a weighting.  We show that if $\un Z$ is compact, finite 
dimensional and oriented, then it carries a fundamental class $[Z]$.
Adapting a construction of Liu and Tian,
we also show that the fundamental class $[X]$ of
any oriented orbifold $\un X$ may be represented
 by a map $\un Z\to \un X$, where the branched manifold $\un Z$
  is unique up to a natural equivalence relation. This gives further 
  insight into the structure of the virtual moduli cycle 
  in the new polyfold theory recently constructed by
  Hofer, Wysocki and Zehnder.
\end{abstract}

 \maketitle

 %\begin{center} Working draft\end{center}
 
%%%%%%%%%%%%%%%%%%%%%%%%%%%%%%%%%%%%%%%%%%%%%%%%%%%%%%%%%%
\section{Introduction}
%%%%%%%%%%%%%%%%%%%%%%%%%%%%%%%%%%%%%%%%%%%%%%%%%%%%%%%%%%

Cieliebak--Mundet i Riera--Salamon~\cite{CRS} formulate the definition 
of branched submanifold of $\R^n$ in connection 
with their discussion of multivalued sections, and use it to represent
the Euler class of certain $G$-bundles (where $G$ is a compact Lie group).  
In Definition~\ref{def:brman} we propose
 an intrinsic definition of a weighted branched manifold 
 $\un Z$ generalizing that in Salamon~\cite{SalFl}.    It is obtained from 
 the usual definition 
of orbifold groupoid simply by relaxing the properness condition and adding
a weighting.  Proposition~\ref{prop:fclass} states that
 if $\un Z$ is compact, finite dimensional and oriented, then $\un Z$
 carries a fundamental class $[Z]$.
Our point of view allows us to deal with a few technical issues that 
arise when branched submanifolds are not 
embedded in a finite dimensional ambient space.  
 Our other main result can be stated informally as follows.
 (For a formal statement, see 
 Propositions~\ref{prop:resol} and \ref{prop:fclass}.)

\begin{thm}\labell{thm} Any compact oriented orbifold $\un Y$ has a \lq\lq 
     resolution" $\un\phi:\un Z\to \un Y$ by a branched manifold that is unique 
 up to a natural equivalence relation. Moreover ${\un\phi}\,\!_*([Z])$ is the
 fundamental class of $\un Y$.
 \end{thm} 

On the level of groupoids, the resolution 
 is constructed from an orbifold groupoid by refining the objects
 and also throwing away some of the morphisms; cf. 
 Example~\ref{ex:tear}.  One can think of
 resolutions as trading  orbifold
  singularities for branching. One application is to give a simple description of the (Poincar\'e dual of the) Euler class of a bundle over an orbifold
  as a homology class represented by a branched manifold; see Proposition~\ref{prop:euler} and~\S\ref{ss:euler}.    We work in finite dimensions but, as the discussion below indicates, 
  the result also applies in certain infinite dimensional situations.

 We  give two proofs of Theorem~\ref{thm}.  The first (in \S\ref{ss:resol}) 
 is an explicit functorial  construction that builds $\un Z$ from a set of
 local uniformizers 
 of $\un Y$. It is a groupoid version of  Liu--Tian's
construction of the virtual moduli cycle in ~\cite{LiuT}. (Also see Lu-Tian~\cite{LuT}.)  The second
 (in \S\ref{ss:multi}) constructs $\un Z$ as the graph of a multisection 
 of a suitable orbibundle  $\un E\to \un Y$.  It therefore relates to
 the construction of the Euler class in~\cite{CRS} and to
Hofer, Wysocki and Zehnder's new polyfold\footnote
{
For the purposes of the following discussion a polyfold can be understood 
as an orbifold in the category of Hilbert spaces. For more detail see~\cite{H,HWZ}.
}
approach to constructing the
 virtual moduli cycle of symplectic field theory.

 The situation here is the following. 
The generalized Cauchy--Riemann (or delbar) operator is a global Fredholm section $f$ of a polyfold bundle $\un E\to \un Y$, and if it were transverse 
to the zero section one would define the virtual moduli cycle to be its zero set.  However, in general $f$ is
 not transverse to the zero section.  Moreover, 
because $\un Y$ is an orbifold rather than a manifold, one can achieve transversality only by perturbing $f$ by a multivalued section $s$.
 Hence the virtual moduli cycle,
 which is defined to be the zero set of $f+s$, is 
a weighted branched submanifold of the infinite dimensional groupoid $\un E$: see~\cite[Ch~7]{HWZ}. Since $s$ is chosen so that
$f+s$ is Fredholm (in the language of~\cite{HWZ} it is an sc$^+$-section), this zero set is finite dimensional.

Although one can define the multisections $s$ fairly explicitly,
the construction in~\cite{HWZ} gives little insight into the
topological structure of the corresponding zero sets.  
This may be understood in terms of a 
resolution $\un\phi:\un Z\to \un Y$.  The pullback by $\un\phi:\un Z\to \un Y$ of an orbibundle $\un E\to \un Y$ is a  bundle ${\un\phi}^*(\un E)\to \un Z$ and we show in Lemma~\ref{le:multi1} that its
(single valued) sections $s$ push forward to give multivalued sections ${\un\phi}\,\!_*(s)$ of $\un E\to \un Y$ in the sense of~\cite{HWZ}.
%Since $s$ is single valued over $\un Z$ there are natural notions of 
%transversality etc.  
Proposition~\ref{prop:sect} shows that if $\un E\to \un Y$ has enough {\em local} sections  to achieve transversality
one can construct the resolution
to have enough {\em global} sections to achieve transversality.
Hence one can understand the virtual moduli cycle as the zero set of a
(single valued) section of $\un\phi^*(\un E)\to \un Z$. In particular, its branching is induced by that of $\un Z$.

This paper is organized as follows.  \S2 sets up the language in which to define
orbifolds in terms of groupoids.  It is mostly but not entirely review,
because we treat the properness requirements in a nonstandard way.
In \S3, we define weighted nonsingular branched 
groupoids and branched manifolds and establish their main properties.
\S4 gives the two constructions for the resolution.
 The relation to the work of
 Cieliebak--Mundet i Riera--Salamon~\cite{CRS}
is discussed in~\S\ref{ss:euler}.

I wish to thank Kai Cieliebak, Eduardo Gonzalez, Andr\'e Haefliger, Helmut Hofer  and Ieke Moerdijk  for some very pertinent  questions and comments  on  earlier versions of this paper.

 \tableofcontents

%%%%%%%%%%%%%%%%%%%%%%%%%%%%%%%%%%%%%%%%%%%%%%%%%%%%%%%%%%
\section{Orbifolds and groupoids}
%%%%%%%%%%%%%%%%%%%%%%%%%%%%%%%%%%%%%%%%%%%%%%%%%%%%%%%%%%

Orbifolds (or $V$-manifolds) were first introduced by Satake~\cite{Sat}.  The idea of describing them in terms of groupoids and categories is due to Haefliger~\cite{Hae1,Hae2,Hae3}.
Our presentation and notation is based on the survey by Moerdijk~\cite{Moe}.
Thus we shall denote the spaces of objects and morphisms of a 
small
topological\footnote
{
i.e. its objects and morphisms form topological spaces and all structure maps are continuous}
 category $\Xx$ by the letters $X_0$ and $X_1$ respectively.  The source and target maps are  $s,t:X_0\to X_1$ and the identity map $x\mapsto \id_x$ is $\id:X_0\to X_1$.
The composition map
$$
m: {X_1}\;_s\!\times_t X_1\;\to \;X_1,\quad (\de,\ga)\mapsto \de\ga
$$
has domain equal to the fiber product 
${X_1}\;_s\!\times_t X_1 = \{(\de,\ga): s(\de) = t(\ga)\}$.
We  denote the space of morphisms from $x$ to $y$ by $\Mor(x,y)$.

%%%%%%%%%%%%%%%%%%%%%%%%%%%%%%%%%%%%%%%%%%%%%%%%%%%%%%%%%%
 \subsection{Smooth, stable, \'etale (sse) Groupoids}
%%%%%%%%%%%%%%%%%%%%%%%%%%%%%%%%%%%%%%%%%%%%%%%%%%%%%%%%%% 
 
Throughout this paper we shall work in the smooth category, by which we mean the category of  finite dimensional second countable Hausdorff  manifolds.  If the manifolds 
have boundary, we assume 
that all local diffeomorphisms respect the boundary. Similarly, if they are oriented, all local diffeomorphisms respect the orientation.

% We begin with a 
%slightly nonstandard definition of  \'etale Lie groupoid.  Note that, although it is not mentioned explicitly in their name, by our definition \'etale groupoids always satisfy  the stability condition.
 
 \begin{defn}\labell{def:orb} An {\bf (oriented) sse  groupoid} $\Xx$ is a small topological 
 category such that the following conditions hold.\smallskip
  
 \NI
{\bf (Groupoid)}   All morphisms are invertible. More
  formally there is a structure map $\io:X_1\to X_1$ that takes each $\ga\in X_1$ to its inverse $\ga^{-1}$.
 \SSS
 
 \NI
{\bf (Smooth \'etale)}  The spaces $X_0$  of objects and $X_1$ of morphisms are (oriented) manifolds
 (without boundary), and
 all structure maps ($s,t,m,\io, \id$) are (oriented) local diffeomorphisms. 
 \SSS
 
 \NI
{\bf (Stable)}  For each $x\in X_0$ the set of self-morphisms $\Mor(x,x) =: G_x$ is finite.  \SSS
 
 \NI
 Groupoids that also satisfy the properness condition stated below will be called  {\bf ep groupoids}, where the 
initials ep stand for \lq\lq \'etale proper".\SSS
  
  \NI
  {\bf (Properness)}  The map $s\times t:X_1\to X_0\times X_0$ that takes a morphism to its source and target is proper.\SSS

\NI
The {\bf orbit space}
 $|X|$ of $\Xx$ is the quotient  of $X_0$ by the equivalence
  relation in which $x\sim y$ iff $\Mor(x,y)\ne \emptyset$.   
    \SSS
  
 \NI  Further, $\Xx$ is called \SSS

 \NI$\bullet$ 
 {\bf 
  nonsingular} if all stabilizers $G_x: = \Mor(x,x)$ are trivial. \SSS 
    
   \NI$\bullet$  {\bf effective} if for each $x\in X_0$ 
   and $\ga\in G_x$ 
each
    neighborhood $V\subset X_1$ of $\ga$ 
    contains  a morphism $\ga'$ such that $s(\ga')\ne t(\ga')$  (i.e. the action of $\ga$ is locally effective.)\SSS
  
  \NI$\bullet$ {\bf connected} if $|X|$ is path connected. %\SSS
  \end{defn} 
 
Unless there is specific mention to the contrary, all groupoids $\Xx$ considered in this paper are sse groupoids, understood in the sense defined above, and all functors are smooth, i.e. 
they are given by smooth maps on 
the spaces of objects and morphisms.\footnote
{
Moerdijk~\cite{Moe} formulates the smooth (or Lie) \'etale  condition in a slightly different but essentially equivalent way. 
In his context, equivalence is called Morita equivalence.  Note also that one can work with the above ideas in
categories other than that of finite dimensional Hausdorff manifolds and local diffeomorphisms. For example, as in~\cite{HWZ}, one can work with infinite dimensional $M$-polyfolds and sc-diffeomorphisms. Haefliger~\cite{Hae1,Hae2,Hae3} and Moerdijk--Mr\v cun~\cite{MM} consider Lie groupoids in which the space $X_1$ of morphisms is allowed to be a nonHausdorff manifold in order to accommodate examples such as the groupoid of germs of diffeomorphisms.  Foliation groupoids also need not be proper. Thus they also develop considerable parts of the theory of \'etale groupoids without assuming properness, though in the main they are interested in very different manifestations of nonproperness.}
Many authors call ep groupoids {\bf orbifold groupoids}.
Note that stability\footnote{This terminology, taken
 from~\cite{RobS}, is inspired by the finiteness condition  
satisfied by stable maps.}
 is a consequence of properness, 
 but we shall often assume the former and not the latter.
 
  Robbin--Salamon~\cite{RobS} use the stability condition to show that in any sse groupoid 
every point $x\in X_0$ has an open neighborhood $U_x$ such that
% $G_x\subset X_1$ has a neighborhood diffeomorphic to $G_x\times U_x$
% where
% $$
% s\times t: U_x\times G_x\to U_x\times U_x.
% $$
% In other words, one can choose $U_x$ so that
 each morphism $\ga\in G_x$ extends to a diffeomorphism of $U_x$ onto itself.\footnote
 {
 Note that each $\ga \in G_x$ extends to a local diffeomorphism  $\phi_\ga$
 of $X_0$ in the following way.
 Let $V$ be a neighborhood of $\ga$ in $X_1$ on which the source and target maps $s$ and $t$ are injective.  Then
 $\phi_\ga$ maps $s(V)$ to $t(V)$ by $s(\de)\mapsto t(\de), \de\in V$.
 Thus the point of Robbin--Salamon's argument is to show that we may assume that $s(V) = t(V)= U_x$.}
   However, in  general there could be many other morphisms 
 with source and target in $U_x$. Their Corollary 2.10 states that 
 $\Xx$ is proper iff $U_x$ can always be chosen so that this is not so, i.e. 
 so that $(s\times t)^{-1}( U_x\times U_x) \cong
 U_x\times G_x$.   Thus properness is equivalent to the existence of local uniformizers in the sense of Definition~\ref{def:atl} below.
 
Another well known consequence of properness is that the orbit space 
$|X|$ is Hausdorff. The next lemma shows that $|X|$ is Hausdorff
iff the
equivalence relation on $X_0$ is closed, i.e. the subset
$\{(x,y)\in X_0\times X_0: x\sim y\}$ is closed. 

\begin{lemma}\labell{le:Haus}  Let $\Xx$ be an sse groupoid. Then:
\SSS

\NI
{\rm (i)}  the projection $\pi:X_0\to |X|$ is open.  \SSS

\NI
{\rm (ii)} $|X|$ is Hausdorff iff $s\times t$ has closed image.
%\NI
%{\rm (iii)}  $\Xx$ is locally bounded iff the function $x\mapsto |G_x|$ is locally bounded on $X_0$ and the function $p\mapsto \#(\pi^{-1}(p))$ is locally bounded on $|X|$.\SSS

%\NI
%{\rm (iv)}  If $\Xx$ is locally bounded and effective, then $\Xx$ is proper iff 
%$s\times t$ has closed image.  FALSE
\end{lemma}
\begin{proof}  Let $U\subset X_0$ be open.  
We must show that $|U|: = \pi(U)$ is open in the quotient topology, i.e. that $\pi^{-1}(\pi(U)) 
= t(s^{-1})(U)$ is open in $X_0$.  But $s^{-1}(U)$ is  open in $X_1$ since $s$ is continuous
and $t$ is an open map because it is a local diffeomorphism. This proves (i).

Now suppose that
 $\im(s\times t)$ is closed and let $p,q$ be any two distinct points in $|X|$.  Choose $x\in \pi^{-1}(p)$ and $y\in \pi^{-1}(q)$.
 Then $(x,y)
\notin \im(s\times t)$ and so there is a neighborhood of $(x,y)$ of the form $U_x\times U_y$ that is disjoint from
$\im(s\times t)$.  Then  $|U_x|$ and $ |U_y|$ are open in $|X|$
by (i).  If
$|U_x|\cap |U_y|\ne \emptyset$  there is $x'\in \pi^{-1}(|U_x|), 
y'\in \pi^{-1}(|U_y|)$ such that $x'\sim y'$.  Then there is $x''\in U_x$ such that $x'\sim x''$ and  
$y''\in U_y$ such that $y'\sim y''$.  Therefore by transitivity
$x''\sim y''$, i.e.  $U_x\times U_y$ meets $\im(s\times t)$, a 
contradiction.  Hence $|U_x|$ and $ |U_y|$ are disjoint neighborhoods of $p,q$, and $|X|$ is Hausdorff.  
The proof of the converse is similar. Hence (ii) holds. \end{proof}
%
%(iii) is straightforward.  So consider (iv).
%Since the image of a proper map is always closed, it remains to show that if $\Xx$ is effective and locally bounded and if $\im(s\times t)$ is closed then $\Xx$ is proper.
%To see this, let $K\subset X_0\times X_0$ be compact.  Then $K': = K\cap \im(s\times t)$ is compact. For $i=1,2$ let $\pi_i: X_0\times X_0\to X_0$ denote the projection to the $i$th component. Since $X_0$ is locally compact, we may
% cover $\pi_1(K)\subset X_0$ by a finite number of
%  connected open sets $U_{j}$ 
% with compact closure
% such that  for each $j$ $\pi^{-1}(U_j)$ is a finite union of 
% connected sets $U_{jk}, 1\le k\le n_j,$ on which $s$ is injective.

% If some $U_{jk}$ does not have compact closure in $X_1$
% then $s(N_{jk})\ne N_j$, 
%and there is a sequence $\ga_n\in N_{jk}$ such that $s(\ga_n)$ converges to a point in $N_j\setminus N_{jk}$.  

%Note that 
%$s\bigl((s\times )^{-1}(K')\bigr) = \pr_1(K')\subset X_0$ 
%\end{proof}

 \NI
\begin{rmk}\labell{rmk:I}\rm (i) The properness condition
 is essential to distinguish 
orbifolds from branched manifolds.  For example, it is easy to define a (nonproper) sse groupoid $\Bb$  with orbit space $|B|$ equal to the quotient of the disjoint union $B_0: = (0,2)\sqcup (3,5)$  by the equivalence relation $x\sim y$ for $y=x+3, x\in (0,1)$;
cf. Fig.~\ref{fig0} and Example~\ref{ex:brorb} (i). 
%%D added
Note that $|B|$ is nonHausdorff, but can be made Hausdorff by identifying the points $1$ and $4$; cf. Lemma~\ref{le:maxH}. \SSS
The transverse holonomy groupoid of the Reeb foliation (cf. Haefliger~\cite{Hae2})  is also of this kind.

\begin{figure}[htbp] %  figure placement: here, top, bottom, or page
%xx   \centerline{\psfig{figure=orbfig0.jpg, width=2in}}
   \includegraphics[width=2in]{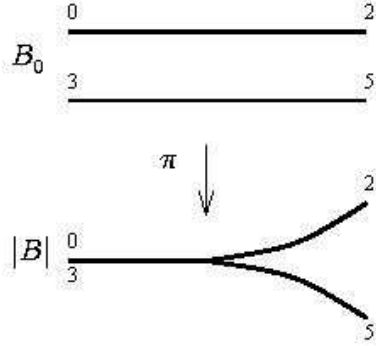} 
   \caption{The projection $\pi$ for a  nonproper groupoid}
   \label{fig0}
\end{figure}

\NI
(ii)  To say that $s\times t$ is a closed map is different from saying that its image is closed.  For example, one could define a noneffective groupoid $\Xx$ with objects $X_0 = S^1 = \R/\Z$ and with $s=t$
but such that $X_1$ has infinitely many components.  Such $\Xx$ 
can be stable, for example if $X_1 = X_0\sqcup(\sqcup_{k\ge 3}V_k)$
 where 
 %%D changed
 each element in $V_k$ has order $2$ and
$s(V_k) = t(V_k) =(\frac 1{3k+2},\frac 1{3k})$.
Then the image of $s\times t$ is the diagonal in $X_0\times X_0$ and so is closed. But the map $s\times t$ is not closed. 
(If $A = \{\ga_k\in V_k: s(\ga_k) = \frac1{3k+1}, 
k\ge 3\}$, $A$ is closed but its image under $s\times t$ is not.)
Note that, even when $\Xx$ is proper,
 the fibers of the projection $\pi:X_0\to |X|$ 
 need not be finite. \SSS

\NI(iii)
If an ep groupoid $\Xx$ is nonsingular, then the orbit space $|X|$ is a manifold.  
\end{rmk}

%Here is another example of an sse but nonproper groupoid.
%Let $Z_0=S^2$ with morphisms $Z_1: = Z_0\;\sqcup \;D\times \{\phi\}$,
%where the open disc $D$ is a subset of $S^2$ and the morphisms $(x,\phi),x\in D,$ have both source and target equal to $x$.  If we define $(x,\phi)\circ (x,\phi) = (x,\id_x)$ for all $x\in D$,
%these conditions define a groupoid $\Zz$.  Since the map 
%$s\times t: D\times \{\phi\}\to S^2\times S^2$ is not proper, $\Z$ is not proper.  It is also nonreduced.  In fact, it 
The groupoid in (ii) above
has the awkward property that, despite being connected, some of its points have 
effective stability groups and some do not.  
The following lemmas show that this cannot happen in the proper case.   Though they are well known, we include 
them for the sake of completeness.

\begin{lemma} \labell{le:red1}  Suppose that $\Xx$ is a connected  ep groupoid.
 Suppose further that for some $x\in X_0$ 
every element of the stabilizer group $G_x$ acts effectively.  Then $\Xx$ is effective.
\end{lemma}
\begin{proof} If every element of  $G_x$ acts effectively and $\pi(x) = \pi(y) \in |X|$, then the elements of $G_y$ act effectively.  Therefore we may partition $|X|$ into two disjoint sets  
$W_e$ and $W_n$, the first being the image of points where $G_x$ acts effectively and the second
 the image of points where the action of $G_x$ is  not effective. 
  The set $W_n$ is always open as it is the image 
 of the open subset 
$$
V: = \{\ga\in X_1: s(\ga') = t(\ga') \mbox{ for all } \ga' \mbox{ in some neighborhood of } \ga\}
$$
under the composite $\pi\circ s$, which is open by 
Lemma~\ref{le:Haus}(i).
If $\Xx$ is proper, then $W_n$ is also closed.
For if not, there is a convergent sequence $p_k\in W_n$ whose limit $p\notin W_n.$  Choose $x\in \pi^{-1}(p)$ and a sequence $x_k\in
\pi^{-1}(p_k)$ with limit $x$.   For each $k$ there is $\ga_k\in G_{x_k}\setminus \{\id\}$ 
that acts as the identity near $x_k$.
Since the sequence $(x_k,x_k)$ converges to $(x,x)$, properness of the map $s\times t$ implies that the $\ga_k$ have a convergent subsequence (also called $\ga_k$) with limit $\ga$.  Then $\ga\in G_x$ and hence has finite order.  On the other hand it extends to a local diffeomorphism near $x$ that equals the identity near the $x_k$.  
This is possible only if $\ga= \id$.  But this is impossible since $\ga_k\ne \id$ for all $k$ and the set of identity morphisms 
 is open in $X_1$.    Since $|X|$ is assumed connected and $W_e$ is nonempty by hypothesis, we must have $W_e=|X|$.
\end{proof}  

\begin{lemma} \labell{le:red2}  Suppose that $\Xx$ is a connected ep groupoid. Then: \SSS

\NI
{\rm (i)} the isomorphism class of the  subgroup $K_y$ of $G_y$ 
that acts trivially on $X_0$ is independent of $y\in X_0$; \SSS

\NI
{\rm (ii)} there is an associated effective groupoid $\Xx_{\it eff}$ 
with the same objects as $\Xx$ and where the morphisms from $x$ to $y$ may be identified with the quotient $\Mor_{\Xx}(x,y)/K_y$.
\end{lemma}
\begin{proof} 
 Given $\ga\in K_y$ denote by $V_\ga$ the component of $X_1$ containing $\ga$.  Then, by properness, the image of $s: V_\ga\to X_0$  is the component $U_y$ of $X_0$ containing $y$. Moreover, as in the previous lemma, $s=t$ on $V_\ga$.   
   Since this holds for all
$\ga\in K_y$, the groups $K_z, z\in U_y,$ all have the same number of elements.  Moreover, they are isomorphic because the operation of composition takes $V_\de\times V_\ga$ to $V_{\de\circ\ga}$. Statement (i) now follows because
$|X|$ is connected, and the groups $K_y$ are isomorphic 
as  $y$ varies in a fiber of $\pi:X_0\to |X|$.

To prove (ii),  define an equivalence relation on $X_1$ by 
 setting $\de\sim\de'$ iff the morphisms $\de,\de':x\to y$ are such that $\de\circ(\de')^{-1}\in K_y$.  It is easy to check that these equivalence classes form the morphisms of the category $\Xx_{\it eff}$.
 \end{proof}

The following definitions are standard.
 
  \begin{defn}\labell{def:refin}  Let $\Xx,\Xx'$ be sse groupoids.
A functor $F:\Xx'\to\Xx$ is said to be {\bf smooth} if the induced maps $X_i'\to X_i, i=0,1,$ are smooth.
The pair $(\Xx',F)$ is said to be a {\bf refinement} of $\Xx$
 if  \SSS

\NI
{\rm (i) } The induced map $F: X_0'\to X_0$ is a (possibly nonsurjective) local diffeomorphism 
 that induces a homeomorphism $|X'|\to |X|$;\SSS

\NI
{\rm (ii) }   For all $x'\in X_0'$, $F$ induces an isomorphism
$G_{x'}\to G_{F(x')}$.\SSS

\NI
Two sse groupoids are {\bf equivalent} if they have a common refinement.
\end{defn}

\begin{rmk}\labell{rmk:refin}\rm
If $(\Xx',F)$ refines $\Xx$ then the morphisms in 
$\Xx'$ are determined by the map $F:X_0'\to X_0$.  Indeed for any pair
$U,V$ of components of $X_0'$ the space $\Mor_{\Xx'}(U,V)$ of morphisms
 with source in $U$ and target 
in $V$ is
$$
\Mor_{\Xx'}(U,V) = \{(x,\ga,y)\in U\times X_1\times V|
s(\ga) = F(x), t(\ga) = F(y)\}.
$$
For short, we will say that the morphisms in $\Xx'$ are pulled back from those in $\Xx$.  
 Moreover $F:X_0'\to X_0$ can be any local diffeomorphism whose image surjects 
 onto $|X|$.  In particular, if $\Uu= \{U_i\}_{i\in A}$ is any collection of 
 open subsets of $X_0$ that projects to a covering of $|X|$, there is a unique refinement $\Xx'$ of $\Xx$ with objects
 $\sqcup_{i\in A} U_i$.  It follows that any category $\Xx''$ with the same objects as $\Xx$ but fewer morphisms is significantly 
 different from $\Xx$; cf. Example~\ref{ex:tear}.  Later we will see that under certain 
 conditions
the corresponding functor $\Xx''\to \Xx$ has the structure of a layered covering; cf.  Lemma~\ref{le:tame}.
\end{rmk}

The proof that the above notion of equivalence is
an equivalence relation on sse groupoids is based on the fact that if
 $F':\Xx'\to \Xx, F'': \Xx''\to \Xx$ are  two refinements of $\Xx$ their fiber product
$\Xx''\times_{\Xx}\Xx'$ refines both $\Xx'$ and $\Xx''$.  Here we use the so-called weak fiber product of Moerdijk--Mr\v cun~\cite{MM} (see also~\cite[\S2.3]{Moe})
with objects\footnote
{
 Note that we cannot take the objects to be  the usual (strict) fiber product 
 $$
X_0''\times_{X_0}X_0' =\{(x'', x')\in X_0''\times X_0': F''(x'') = F'(x')\}.
$$
This definition is useful only when one of the maps $X_0'\to X_0, X''\to X_0$ is surjective; for a general equivalence, the above space 
might be empty. We also cannot use  the fiber product
$$
X_0''\times_{|X|}X_0'' = \{(x'', x')\in X_0''\times X_0': \pi\circ F''(x'') = \pi\circ F'(x')\}
$$
since this is only guaranteed to  be a manifold if one of the projections $X_0'\to |X|, X_0''\to |X|$ is a local submersion.}
 given by the  homotopy pullback 
$$
\{(x'',\ga, x')\in X_0''\times X_1\times X_0': F''(x'')=\ga(F'(x')) \}.
$$
Morphisms 
$(x'',\ga, x')\to (y'',\de, y')$ are pairs $(\al'',\al')\in X_1''\times X_1'$, where $\al':x'\to y', \al'':x''\to y''$, such that the following diagram commutes:
$$
\begin{array}{ccc}
F'(x')&\stackrel{F'(\al')}\to &F'(y')\\
\ga\downarrow&&\de\downarrow\\
F''(x'')&\stackrel{F''(\al'')}\to &F''(y'').\end{array}
$$
Thus $(\al'',\al'): (x'',\ga, x')\to
 \bigl(\al''(x''), F''(\al'')\circ \ga\circ F'(\al')^{-1},\al'(x')\bigr)$.

 \begin{lemma}\label{le:weak}  Let $\Xx, \Xx',\Xx''$ be sse groupoids and $F':\Xx'\to\Xx, F'':\Xx''\to \Xx$ be smooth functors.\SSS
 
 \NI
{\rm (i)} If  $\Xx'$ is nonsingular  and $F''$ is injective 
on each group $G_{x''}$ then $\Xx''\times_{\Xx}\Xx'$ is nonsingular.\SSS

\NI
{\rm (ii)}  $F''$ is an equivalence iff  the projection $\Xx''\times_{\Xx}\Xx'\to \Xx'$ is an equivalence.
\end{lemma}
\begin{proof}  This is straightforward, and is left to the reader.
\end{proof}

%%%%%%%%%%%%%%%%%%%%%%%%%%%%%%%%%%%%%%%%%%%%%%%%%%%%%%%%%%%%
\subsection{Orbifolds and atlases}\labell{ss:orb}
%%%%%%%%%%%%%%%%%%%%%%%%%%%%%%%%%%%%%%%%%%%%%%%%%%%%%%%%%%%%

\begin{defn}\labell{def:orbstr}
An {\bf orbifold structure} on a 
 paracompact Hausdorff 
 space $Y$ is a pair $(\Xx, f)$ consisting of an ep groupoid $\Xx$ together with a 
homeomorphism $f:|X|\to Y$.  Two orbifold structures $(\Xx, f)$ and
$(\Xx', f')$ are equivalent if they have a common refinement, i.e. if there is a third structure $(\Xx'', f'')$ and refinements $F:\Xx''\to\Xx, F':\Xx''\to \Xx'$ such that 
$f'' = f\circ |F| = f'\circ |F'|$.
  
  An {\bf orbifold} $\un Y$ is a second countable
paracompact
Hausdorff space $Y$ equipped 
with an equivalence class of orbifold structures. 
An {\bf orbifold  map} $\un \phi: \un X\to \un Y$ is 
an equivalence class of functors  $\Phi:\Xx\to 
\Yy$, where $(\Xx, f)$ and
$(\Yy, g)$  are orbifold structures on 
$\un X$ and $ \un Y$, respectively. The equivalence relation is generated by the obvious notion of refinement of functors.  Thus if $F:\Xx'\to \Xx$ and $F':\Yy'\to \Yy$ are refinements, $\Phi':\Xx'\to \Yy'$  is said to refine $\Phi:\Xx\to \Yy$ if there is a natural transformation between the two functors $\Phi\circ F, F'\circ\Phi': \Xx'\to \Yy$.
 \end{defn}
 
 Each orbifold  map $\un \phi: \un X\to \un Y$ induces a well defined
continuous map $\phi:X\to Y$ on the spaces $X,Y$ underlying $\un X, \un Y$. 
  Note that it is possible for different
  equivalence classes of functors to induce the same 
map $X\to Y$. For more details, see
 Moerdijk~\cite[\S2,3]{Moe}.  As pointed out by Lerman~\cite{Ler}, 
 one really should take a more sophisticated approach to 
 defining orbifolds and orbifold maps in order for 
 them to form some kind of category.
 Since the focus here is on defining a new class of objects, 
 we shall not deal with these subtleties. 
 
\begin{example}\labell{ex:tear}\rm  {\bf The teardrop orbifold and its resolution.}  The teardrop orbifold  has underlying space $Y=S^2$ and one singular point $p$ at the north pole of order $k$. 
Cover $Y$ by two open discs $D_+, D_-$ of radius $1+\eps$ that intersect in the annulus $A = (1-\eps,1+\eps)\times S^1$ and are invariant by rotation about the north/south poles. Denote by $\phi:A\to A$ the $k$-fold covering map given in polar coordinates by $(r,\theta)\mapsto (2-r,k\theta)$ and by $R_{t}$  rotation through
the angle $2\pi t$. 
An orbifold structure on $Y$ is provided by the proper groupoid $\Xx$ whose space of objects $X_0$ is the disjoint union $D_+\sqcup D_-$, with morphisms 
$$
X_1: = \bigl(X_0\times\{1\} \bigr)\;\sqcup\; \bigl(D_+\times 
\{\ga_1,\dots,\ga_{k-1}\}\bigr)\; \sqcup \; 
\bigl( A\times\{\phi_+, \phi_-\}\bigr).
$$
Here $1$ acts as the identity and $(x,\ga_j), x\in D_+,$ denotes the morphism with target $x\in D_+$ and source $R_{j/k} (x)$ in $D_+$.  Further  the pair 
$(x,\phi_+)$ is the unique morphism
with source  $x\in A_+: =A\cap D_+$ and target 
 $\phi (x)\in A_-: = A\cap D_-$,  By definition $(x,\phi_-): = (x,\phi_+)^{-1}.$ Thus $|X| = D_+\cup D_-/\sim$ where $x\in D_+$ is identified with $R_{j/k}(x)$ and $D_+$ is attached to $D_-$ over $A$ by a $k$ to $1$ map. 

\begin{figure}[htbp] %  figure placement: here, top, bottom, or page
% xx  \centerline{\psfig{figure=orbfig1.jpg,width=3in}}
\centering   \includegraphics[width=3.5in]{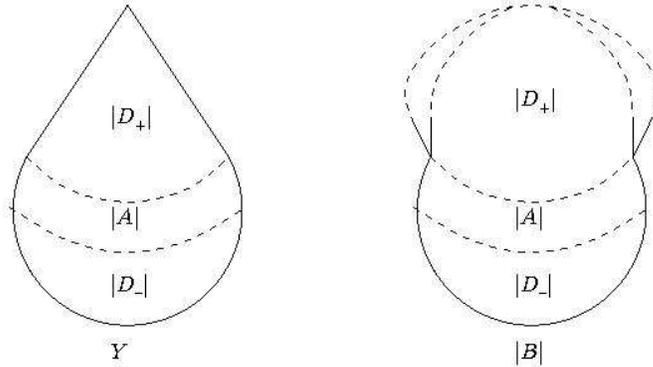} 
   \caption{The teardrop orbifold and its resolution $\Bb$ when $k=2$}
   \label{fig1}
\end{figure}

By way of contrast, consider the groupoid $\Bb$ formed from $\Xx$ by omitting the morphisms $(D_+\setminus A_+)\times 
\{\ga_1,\dots,\ga_{k-1}\}$.  This groupoid is nonsingular, but is no longer proper since the restriction of $(s,t): B_1\to B_0\times B_0$ to the component
of morphisms $ A_+\times\{\ga_1\}$  is not proper.  Note that $|B|$ is a branched manifold:   $k$ (local) leaves come together 
along the boundary $\p |D_-|\subset |B|$.
Further, if one weights the leaves over $|D_+|$ by $1/k$ the induced map $|B|\to |X|$ represents the fundamental class of $|X|=Y$. We call $\Bb$ the resolution of $Y$.
\end{example}

 Often it is convenient to describe an orbifold structure on $\un Y$ in terms of local charts.
 Here are the relevant definitions. 
 
\begin{defn}\labell{def:atl} A {\bf local uniformizer} $(U_i, G_i, \pi_i)$ for an (orientable) orbifold $\un Y$ is a triple consisting of
 a connected  open subset $U_i\subset \R^d$, a finite group 
 $G_i$ that acts by (orientation preserving) diffeomorphisms of $U_i$
and a map $\pi_i: U_i\to Y$ that factors through a homeomorphism from the quotient $U_i/G_i$ onto an open subset $|U_i|$ of $Y$.
Moreover, this uniformizer determines the smooth structure of $\un Y$
over $|U_i|$ in the sense that for one (and hence any) orbifold structure $(\Xx,f)$ on $\un Y$ the projection $f\circ \pi:X_0\to Y$ lifts  to a local diffeomorphism $(f\circ \pi)^{-1}(|U_i|) \to U_i$.

A {\bf good atlas} for $\un Y$ is a collection $\Aa=\{(U_i, G_i, \pi_i), i\in A\}$
of local uniformizers whose images
 $ \{|U_i|, i\in A\}$ form a locally finite covering of 
 $Y$. \end{defn}

It is shown in Moerdijk--Pronk~\cite[Corollary~1.2.5]{MPr}
that every effective orbifold has such an atlas.
The argument works equally well in the general  case
since it is based on choosing an adapted triangulation of $|X|$.
In fact, one can also assume that the $U_i$ and their images $|U_i|$ are contractible and that the $|U_i|$ are closed under intersections.   
Here we assume as always that $|X|$ is finite dimensional.  The arguments of~\cite{MPr} do not apply in the infinite dimensional case.
However 
Robbin--Salamon~\cite[Lemma~2.10]{RobS} show that in any category an orbifold has a good atlas. 
They start from a groupoid structure $(\Xx,f)$ on $\un Y$ and 
construct for each $x\in X_0$ a local uniformizer $(U,G,\pi)$ that embeds in $\Xx$ in the sense that $U\subset X_0$ and
the full subcategory of $\Xx$ with objects $U$ is isomorphic to $U\times G$.   

Let $(U_i, G_i, \pi_i), i\in A$, be a locally finite cover of $Y_0$ by such uniformizers, and denote by $\Xx'$ the full
 subcategory of $\Xx$ with objects $\sqcup_{i\in A} U_i$.
 Since this is a refinement of $\Xx$, it is another orbifold structure on $\un Y$.  Observe that the morphisms in $\Xx'$ with source and target $U_i$ can be identified with $U_i\times G_i$.  
 In this situation we say that $\Xx'$ is 
 {\bf constructed from the good atlas}
$(U_i, G_i, \pi_i), i\in A$.  
Many interesting orbifolds, for example the ep polyfold groupoids considered in~\cite{HWZ}, are constructed in this way.

The following definition will be useful.

\begin{defn}\labell{def:sm} Let $\Xx$ be an ep groupoid.  
A point $p\in |X|$ is said to be smooth if for one (and hence every)
point $x\in \pi^{-1}(p)\subset X_0$ every element $\ga\in G_x$ acts
trivially near $x$.
\end{defn}

Thus when $\Xx$ is effective $p$ is smooth iff $G_x = \{\id\}$
for all $x\in \pi^{-1}(p)$.  Note that 
for any  local uniformizer  $(U, G,\pi)$ 
 the fixed point set of each $\ga\in G$ is either the whole of $U$ or is nowhere dense in $U$.  It follows
 that the smooth points $|X|^{sm}$
form an open and dense subset of $|X|$.

 We now show that if $\un Y$ is effective its orbifold structure is uniquely determined by the charts in any good atlas; there is no need for further explication of how these charts fit together.  This result is well known: cf.~\cite[Prop.~5.29]{MM}. 
We include a proof to clarify ideas.  It is particularly relevant 
in view of the work of Chen--Ruan~\cite{CR} and Chen~\cite{Ch} on orbifold Gromov--Witten invariants, that discusses orbifolds from the point of view of charts.
Note that the argument applies to orbifolds in any category, and in particular to polyfold groupoids. 

\begin{lemma}\labell{le:red3} An effective
 orbifold $\un Y$  is uniquely determined by
the charts  $(U_i, G_i, \pi_i), i\in A,$ of 
a good atlas.
\end{lemma} 

\begin{proof}  Here is an outline of the proof.  We shall construct an
orbifold structure $(\Xx,f)$ on $\un Y$, such that $\Xx$ is a groupoid whose objects $X_0$ are 
the disjoint union of the sets $U_i, i\in A,$ and whose morphisms
with source and target in $U_i$ are given by $G_i$. 
The construction uses the fact that $\un Y$ has an effective orbifold structure
$(\Zz,g)$, but the equivalence class of 
$(\Xx,f)$ is independent of the choice of $(\Zz,g)$.  It follows that $(\Xx,f)$ and $(\Zz,g)$ are equivalent, and that the orbifold structure of $\un Y$ is determined by the local charts.

We described above the objects in $\Xx$ and some of the morphisms. To complete the construction, we must add to $X_1$ some
 components $C^{ij}$ of morphisms from 
 $U_j$ to $U_i$ for all $i,j\in A$
  such that $|U_i|\cap |U_j|\ne \emptyset$.  
%Because the $G_i$ act effectively,  
%each component $U_i$ of $X_0$ is the closure of the subset $U_i^{sm}$ %consisting of smooth points, i.e. those with trivial stabilizers.  
%Moreover, because the fixed point set of any nontrivial orientation 
%preserving diffeomorphism has codimension at least two, each set 
%$U_i^{sm}$ is connected.
For smooth points 
$x\in U_i^{sm}, y\in U_j^{sm}$ the set $\Mor(x,y)$ has at most one element
and is non empty iff $\pi_i(x) = \pi_j(y)\in |X|=Y$.  Hence the given data determine the set $X_1^{sm}$ of all morphisms whose source or target is smooth.  We shall see that there is a unique way to complete the morphism space 
$X_1$. 

%Since the $|U_i|$ are closed under finite intersections, 
To begin, fix $i\ne j$ such that $|U_{ij}|: = |U_j|\cap |U_i|\ne \emptyset$.  Given $y\in U_j$ denote by $V_{ji}^y\subset U_j$ the connected component of $U_j\cap \pi^{-1}(|U_{ij}|)$ that contains
$y$.
 The key point  is that for any
 point $x\in U_i$ with $\pi_i(x) = \pi_j(y)$ there is 
 a local diffeomorphism $\phi_{xy}$ 
 from a neighborhood $N_j(y)$ of $y$ in $U_j$ to a neighborhood $N_i(x)$ of $x$ in $U_i$ such that $\pi_j = \pi_i\circ \phi$. This holds because there is an orbifold structure $(\Zz,g)$ on $Y$. Namely, choose $z\in Z_0$ so that $g(z) = \pi_i(x) = \pi_j(y)$.  Shrinking $N_j(y)$ and $N_i(x)$ if necessary, we may suppose that there is a neighborhood $N_z$ of $z$ in $Z_0$ such that 
 $g(N_z) = \pi_i(N_i(x)) = \pi_j(N_j(y))$.  Since $(U_i, G_i,\pi_i)$ is a local chart for $\pi_i(U_i)\subset Y$, there is a local diffeomorphism $\psi_i: N_z\to U_i$ such that $g = \pi_i\circ\psi_i$.
 Moreover, since $\pi_i$ is given by quotienting by the action of $G_i$ there is $\ga\in G_i$ such that 
 $\psi_{xz}(z): = \ga\circ\psi_i(z) = x$. Denoting by  $\psi_{yz}$
 the similar map for $y\in U_j$, we
 may take $\phi_{xy}: = \psi_{xz}\circ (\psi_{yz})^{-1}$.
 
Since $\phi_{xy}$ is determined by its restriction to the dense open set of smooth points, it follows that
 there are precisely $|G_y|$ such local diffeomorphisms, namely the composites $\phi_{xy}\circ \ga, \ga\in G_y$. These local diffeomorphisms form the local sections of a sheaf over $V_{ji}^y$,  whose stalk at $y$ consists of the $|G_y|$ elements $\phi_{xy}\circ \ga$.  
 %Since $U_j$ is connected (but perhaps not simply connected), 
 The existence and uniqueness properties of the $\phi_{xy}$ imply that
 each local section is the restriction of a unique, but possibly multivalued, global section $\si$ of this sheaf. We may identify the graph of $\si$ (which is a manifold) with 
% Each such global section $\si$ defines 
 a component $C_\si$ of the space of morphisms from
 $U_j$ to $U_i$.   Since $|G_y|$ is finite, the source map
  $s:C_\si\to V_{ji}^y$ 
  is a surjective and finite-to-one covering map.    
 The set $\{\si_1,\dots,\si_\ell\}$ of all such global sections 
 (for different choices of $x\in U_i, y\in U_j$) is
  invariant under the action  of $G_j$ by precomposition and of $G_i$ by postcomposition.
  We define the space of morphisms in $\Xx$  from $U_j$ to $U_i$
 to consist of the $\ell$ components $C_{\si_k}, 1\le k\le \ell.$ 
 
 This defines the morphisms in $\Xx$ from  $U_j$ to $U_i$.
 We then complete $X_1$ 
by adding the inverses to these elements and all composites.
The resulting composition operation is  associative because its restriction to the smooth elements is associative. 
Thus $\Xx$ is an sse Lie groupoid.  Since the atlas is locally finite, the projection $X_0\to |X|$ is finite to one.  Moreover, the induced map $|X|\to Y$ is a homeomorphism.  Hence $\Xx$ is proper
because $Y$   is Hausdorff. (See also~\cite[Cor~2.11]{RobS}.) 
% To see it is proper, note first that because the covering $|U_i|$ is locally finite, it suffices to check the properness of $s\times t$ on each component of $X_1$.  We saw above that was true for components $C$ such that $|s(C)|
%\subset |t(C)|$.  The fact that it holds for all components follows because $Y$ is Hausdorff. (See also~\cite[Cor~2.11]{RobS}.) 
Hence $\Xx$ is an ep groupoid.  

Moerdijk--Mr\v cun  prove that $\Zz$ and $\Xx$ are equivalent by 
looking at the corresponding groupoids of germs of diffeomorphisms.
Alternatively, define $\Zz'$ to be the refinement of $\Zz$ with 
objects $\sqcup_{z\in A} N_z$, where $A\subset Z_0$ is large enough that the sets $g(N_z), z\in A$, cover $Y$.  For  each $z\in A$
choose
one of the corresponding local diffeomorphisms
  $\psi_{xz}$ and call it $f_z$.  Then define $F:\Zz'\to \Xx$ by setting 
  $$
  F|_{N_z} = f_z,\qquad F(\ga) = f_w\circ \ga\circ (f_z)^{-1},
  $$
  where $s(\ga) = z$ and $t(\ga) = w$. Then  $F$ is an equivalence, as required.
 \end{proof}

\begin{rmk}\labell{rmk:nonred0}
\rm   This lemma is false when $\un Y$ is not effective.
To see this, let $K$ be the cyclic group $\Z/3\Z$ and consider
two groupoids $\Zz,\Xx$ both with objects $S^1$. We assume  in both cases that at all objects the stabilizer group is $K$ and that there 
are no other morphisms. This implies that  $s=t$. We define $\Zz$ to be topologically trivial, with 
$Z_1 = S^1\times K$ (i.e. three copies of $S^1$) and $s=t$ equal to the projection $S^1\times K\to S^1$.  On the other hand, we define $X_1$ to be the disjoint union of two copies of $S^1$.
In this case, $s=t:X_1\to X_0$ is the identity on one circle (the one
 corresponding to the identity morphisms) and is a double cover on the second.  It is easy to check that these groupoids are not equivalent.
(One way to see this is to notice that their classifying 
spaces\footnote{
The classifying space $B\Xx$ is the realization of the nerve of
the category, and is denoted $|\Xx_{\bullet}|$ in~\cite{Moe}.}
$B\Zz, B\Xx$ are not homotopy equivalent.)
However they have the same local uniformizers over any proper open subset of $S^1$.  

These groupoids $\Zz,\Xx$ are totally noneffective. According to Henriques--Metzler~\cite{HM}, the best way to understand their structure is to think of them as a special kind of gerbe.
\end{rmk}

%  Then this groupoid   
% satisfies conditions (i) and (ii) above.  Although 
% it may not satisfy (iii), the components of $X_1'$ do project via $s$ to a locally finite covering of $X'_0$.  Hence, as in the remarks after Definition~\ref{def:refin}, 
% there is a refinement $\Xx$ of this subcategory, with $X_0 = X'_0$, that satisfies all three conditions.  (In forming $\Xx$ one should only subdivide those components of  $X_1'$ with source and target in different $U_i$.)
% Since $\Xx$ is
% equivalent to $\Xx'$, it provides an orbifold structure on $\un Y$.
% 
%Now consider a general orbifold $\un Y$ with 
%orbifold structure $(\Yy, g)$.  Choose a good atlas 
% $(U_i, G_i, \pi_i), i\in A',$ for $\un Y$.
 % such that the covering
% of $Y$ by the sets
% $|U_i|, i\in A',$ is subordinate to its
%  covering by the components of $Y_0$. 
%Construct an orbifold structure $\Xx_{\it eff}$ for 
%$(\un Y)_{\it eff}$ as above
% from the uniformizers corresponding to any subset $A$ of $A'$ that gives a locally finite covering of $\un Y$.

% 
% This clearly satisfies (i) and satisfies (ii) by construction.   However there may be components of
% $(X'')_1$ for which (iii) does not hold, since the covering map $\rho:Y_1\to (Y_{\it eff})_1$ need not be trivial over 
% the full image of $\Phi'$.  But again the components of $(X'')_1$ project via $s$ to a locally finite cover of $(X'')_0$,
% and so  there is a refinement $\Xx$ of $\Xx_{\it eff}$ 
% that also satisfies (iii).  
%  This completes the proof.
%\end{proof}

%%%%%%%%%%%%%%%%%%%%%%%%%%%%%%%%%%%%%%%%%%%%%%%%%%%%%%%%%%%%
\subsection{Fundamental cycles and cobordism}\labell{ss:fnorb}
%%%%%%%%%%%%%%%%%%%%%%%%%%%%%%%%%%%%%%%%%%%%%%%%%%%%%%%%%%%%

We now summarize known facts about homology and cobordism in the context of orbifolds, since we will later generalize   them to the branched case.
Note first that if a compact space $Y$ admits the structure of an 
oriented $d$-dimensional orbifold (without boundary) there is a distinguished $d$-dimensional 
cycle $[Y]$ in the singular chain complex of $Y$
that we will call the {\bf fundamental cycle} of $Y$.  One way to see this is to give $Y$ an adapted triangulation as in Moerdijk--Pronk~\cite{MPr}.
In such a triangulation the open simplices of dimensions $d$ and $d-1$ lie in the smooth\footnote
{
If $\un Y$ is oriented its nonsmooth points have codimension at least $2$ since the fixed 
point submanifold of an orientation preserving linear transformation of $\R^n$ 
has codimension at least $2$.}
part $Y^{sm}$ of $Y$, and $[Y]$ is represented by the sum of the $d$-simplices.  Another way to define $[Y]$ 
is to interpret the inclusion $Y^{sm}\to Y$ as a pseudocycle: see Schwarz~\cite{Sch} or Zinger~\cite{Z}.
These definitions imply that the fundamental cycle of the orbit space $|X|$ of  $\Xx$  is the same as that of its reduction $\Xx_{\it eff}$.\footnote
{Note that here we are talking about the ordinary (singular)
cohomology of $Y$, not the stringy cohomology defined by Chen--Ruan~\cite{CR} (cf also Chen--Hu~\cite{CHu}).}

A smooth map of an orbifold $\un Y$ to a manifold $M$ is a map $\phi: Y\to M$ such that for one (and hence every) 
groupoid structure $\Xx$ on $\un Y$ the induced map $X_0\to M$ is smooth.  Given two such maps $\phi:Y\to M, \phi':Y'\to M$, one can  perturb $\phi$ to be transverse to $\phi'$. This means that there are adapted triangulations on
$Y,Y'$ that meet transversally.  In particular, if 
 $Y$ and $Y'$ are oriented and have complementary dimensions, their images in $M$ intersect in smooth points and  
  the intersection number 
$\phi_*([Y])\cdot \phi'_*([Y'])$ is defined.   
Further, if $\al$ is a $d$-form on $M$ that is Poincar\'e dual to the class 
 $\phi'_*([Y'])$ then 
 $$
 \phi_*([Y])\cdot \phi'_*([Y']) = \int_Y\phi^*(\al) : =  \int_{Y^{sm}}
 \phi^*(\al).
 $$

Suppose now that $\Ww$ is an ep groupoid in the  category of oriented $(d+1)$-dimensional manifolds with boundary.  Then $\Ww$ has a well defined boundary groupoid $\p\Ww$.  If in addition $|W|$
is compact and connected,
the fundamental cycle of $|W|$ generates $H_{d+1}(|W|,\p|W|)$. 
We may therefore make the usual definition of cobordism for orbifolds.
We say that two $d$-dimensional oriented orbifolds $Y_1,Y_2$ are 
{\bf cobordant} if there is a $(d+1)$-dimensional 
groupoid $\Ww$ whose boundary 
$\p \Ww$ is the union of two disjoint oriented
groupoids $-\Xx_1, \Xx_2$, that represent $-Y_1, Y_2$ respectively.
(Here $-\Xx$ denotes the groupoid formed from $\Xx$ by reversing the orientation of $X_0$.)
As usual, if there is a class $\al\in H^{d}(|W|)$ whose 
restriction to $|\p W| = -|X_1|\sqcup |X_2|$ is $-\al_1+\al_2 $, where $\al_i\in H^d(Y_i) \equiv H^{d}(|X_i|)\subset H^d(|\p W|)$, then
$
\al_1([Y_1]) = \al_2([Y_2]).
$

%%%%%%%%%%%%%%%%%%%%%%%%%%%%%%%%%%%%%%%%%%%%%%%%%%%%%%%%%%
\section{Weighted nonsingular branched groupoids}
%%%%%%%%%%%%%%%%%%%%%%%%%%%%%%%%%%%%%%%%%%%%%%%%%%%%%%%%%%

%%%%%%%%%%%%%%%%%%%%%%%%%%%%%%%%%%%%%%%%%%%%%%%%%%%%%%%%%%
\subsection{Basic definitions}
%%%%%%%%%%%%%%%%%%%%%%%%%%%%%%%%%%%%%%%%%%%%%%%%%%%%%%%%%%

We saw in Remark~\ref{rmk:I}(ii) 
above that the category of nonproper sse groupoids 
 contains some rather strange objects, not just branched 
objects such as the teardrop but also  groupoids that 
are not effective 
in which the trivially acting part $K_x$ of the 
stabilizer subgroups is not locally constant.  In order to get a good class of branched groupoids that carries a fundamental cycle we need to impose a weighting condition and to ensure that the groupoids are built by assembling well behaved pieces called local branches.
Definition~\ref{def:brorb} adapts
 Salamon~\cite[Def~5.6]{SalFl} to the present context; it is also 
 very close to the definition in
 Cieliebak--Mundet i Riera--Salamon~\cite{CRS} of a branched submanifold of Euclidean space.  For the sake of simplicity, we shall 
 restrict to the nonsingular and oriented case.  
  
  We shall use the concept of the {\bf maximal Hausdorff quotient} $Y_{\Hsm}$ of a nonHausdorff topological space $Y$.  This is a pair
  $(f, Y_{\Hsm})$ that satisfies the following conditions:
  \SSS
  
  \NI
  (i)  $f:Y\to Y_{\Hsm}$  is a  surjection
  and  $Y_{\Hsm}$ has the corresponding quotient topology,\SSS
  
  \NI
  (ii) $Y_{\Hsm}$ is Hausdorff,\SSS
  
  \NI
  (iii)  any continuous surjection $f_\al:Y\to Y_\al$ with
  Hausdorff  image  factors through $f$.
  
\begin{lemma}\labell{le:maxH}  Every topological space $Y$ has a maximal 
Hausdorff quotient $(f, Y_{\Hsm})$.  
%\NI
%{\rm (ii)}  If $V$ is an open subset of $Y$ that is homeomorphic to an open subset of $\R^d$, then  $f$ induces  a homeomorphism $V\to f(V)$.
 \end{lemma}
  \begin{proof}
  To construct $Y_{\Hsm}$, denote by $A$ the set of all equivalence classes of pairs $(g_\al, Y_\al)$ where $Y_\al$ is a
Hausdorff topological space and $g_\al:Y\to Y_\al$ is a continuous surjection. 
(Two such pairs are equivalent if there is a homeomorphism $\phi:Y_\al\to Y_\be$ such that $\phi\circ g_\al = g_\be$.)  The product $Y_\pi: = \prod_{\al\in A} Y_\al$ is Hausdorff in the product topology and there is a continuous map $f_\pi:Y\to Y_\pi$.
Define $Y_{\Hsm}: = {\im f_\pi}$ with the subspace topology $\tau_s$.  Then 
conditions (ii) and (iii) above hold, and there is a continuous surjection $f:Y\to Y_{\Hsm}$. 
We must show that (i) holds, i.e. that $\tau_s$ coincides with the quotient topology
$\tau_q$.  Since the identity map $(Y_{\Hsm},\tau_q)\to (Y_{\Hsm},\tau_s)$ is continuous, $(Y_{\Hsm},\tau_q)$ is Hausdorff.
Therefore by (iii), the map $f:Y\to (Y_{\Hsm},\tau_q)$ factors through $(Y_{\Hsm},\tau_s)$.  Thus the identity map 
$(Y_{\Hsm},\tau_s)\to (Y_{\Hsm},\tau_q)$ is also continuous; in other words, the topologies $\tau_s$ and $\tau_q$ coincide.
\end{proof}
%To prove (ii) note first that for any point $x\in U$ there is a continuous function $g:Y\to [0,1]$ 
%such that $f(x) = 0$ and $f(z) = 1$ for all $z\notin U$. It follows that 
%$f^{-1}(f(x)) = \{x\}$ for all $x\in U$, and $f^{-1}(f(V))=V$ for all open $V$ in $U$.  Therefore $f$ is an injective open mapping. The result follows.

If $\Bb$ is a  nonproper sse groupoid, then
$|B|$ may not be Hausdorff; cf. Lemma~\ref{le:Haus}. It will be useful to consider the Hausdorff quotient $|B|_{\Hsm}$.\footnote
{Note that $|B|= |B|_{\Hsm}$ if $\Bb$ is proper.}
We shall denote by $|\pi|_{\Hsm}$ the projection $|B|\to |B|_{\Hsm}$
and by $\pi_{\Hsm}$ the projection $B_0\to |B|_{\Hsm}$.  Further 
$|U|_{\Hsm}: = \pi_{\Hsm}(U) $ denotes the image 
of $U\subset B_0$ in $|B|_{\Hsm}$.

\begin{defn}\labell{def:brorb}
A  {\bf weighted nonsingular branched groupoid} (or {\bf wnb groupoid} for short)
is a pair $(\Bb,\La)$
%triple $(\Bb,\la,\La)$, 
consisting of an 
oriented,  nonsingular sse Lie groupoid together with a weighting function
%s $\la:\pi_0(B_0)\to (0,\infty)$  
$\La:|B|_{\Hsm}\to (0,\infty)$ that satisfies the following
 compatibility conditions.
For each $p\in |B|_{\Hsm}$ there is
 an open neighborhood $N$ of $p$ in $|B|_{\Hsm}$, 
 a collection  $U_1,\dots,U_\ell$ of 
 %(possibly disconnected) 
 disjoint open subsets of $
 \pi_{\Hsm}^{-1}(N)\subset B_0$ (called {\bf local branches})
 and a set of positive weights $m_1,\dots,m_\ell$
such that: \SSS

\NI
{\bf (Covering) } $|\pi|_{\Hsm}^{-1}(N) = |U_1|\cup\dots \cup 
|U_\ell| \subset |B|$;\SSS

\NI
{\bf (Local Regularity)}  for each $i=1,\dots,\ell$ the projection $\pi_{\Hsm}: U_i\to 
|U_i|_{\Hsm}$  is %%Dinjective and its image is a relatively
a homeomorphism onto a relatively
 closed subset of $N$;\SSS
 % that contains $p$

\NI
{\bf (Weighting)}  
%each local branch $U_i$ is contained 
%in a single component of $B_0$, and, 
for all $q\in N$, $\La(q)$ is the sum of the weights of the local
branches whose image contains $q$:
$$
\La(q) = 
\underset{i:q\in |U_i|_{\Hsm}}\sum m_i.
$$

\NI
The tuple 
$(N,U_i, m_i)$ is said to form a {\bf local branching structure} at $p$.  
Sometimes we denote it by $(N^p, U_i^p, m_i^p)$  
to emphasize the dependence on $p$. $\Bb$ is called {\bf compact}
if its Hausdorff orbit space $|B|_{\Hsm}$ is compact.
The points $p\in |B|_{\Hsm}$  that have more than one inverse image in 
$|B|$ will be called {\bf branch points}.
\end{defn}

\begin{figure}[htbp] %  figure placement: here, top, bottom, or page
 %for:xxfig  \centerline{\psfig{figure=orbfig2.jpg,width=3in}}
      \centering
   \includegraphics[width=3in]{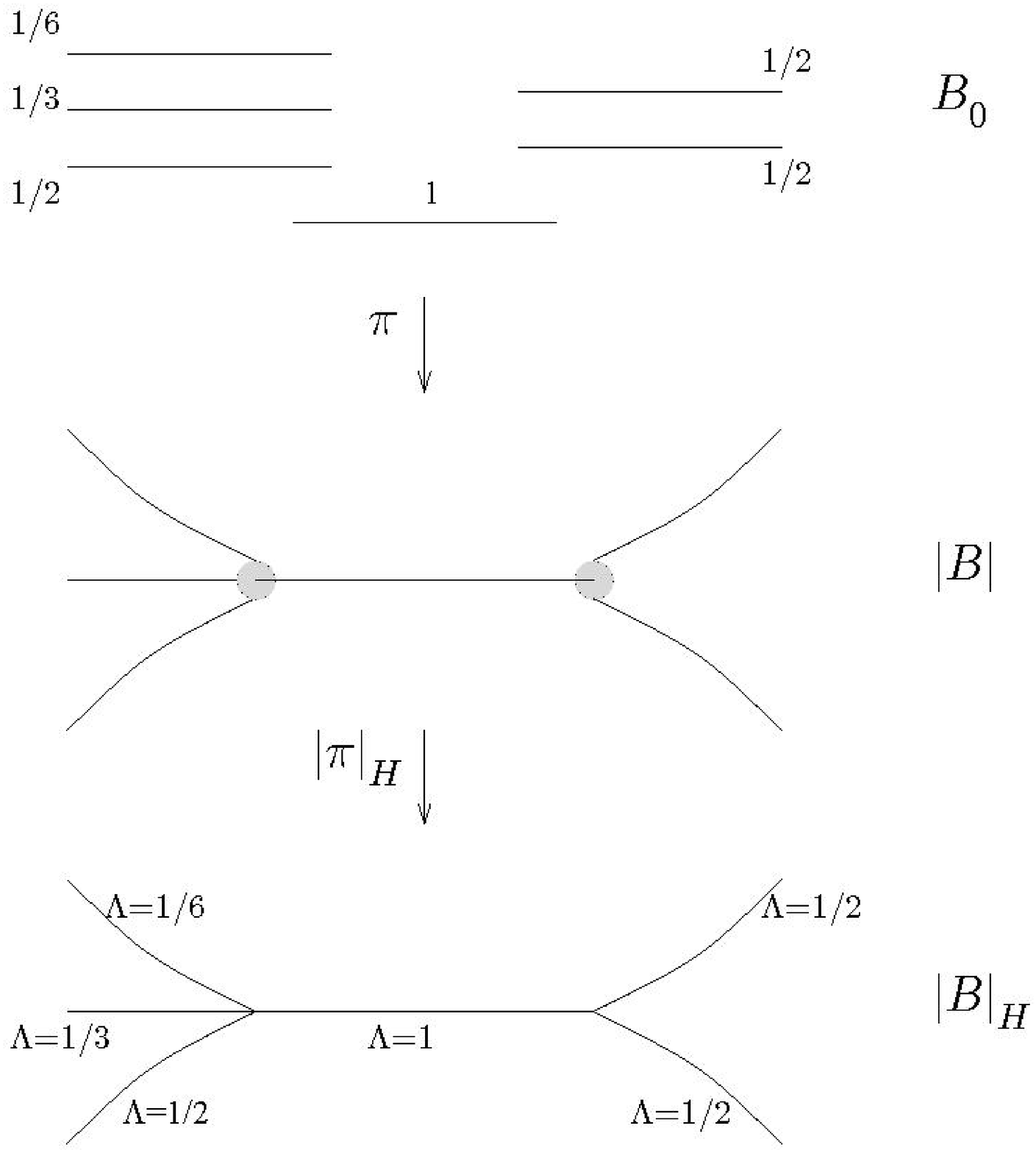} 
   \caption{Local branching structures of a weighted branched groupoid}
   \label{fig2}
\end{figure}

 \begin{example}\labell{ex:brorb}\rm (i) 
% We shall consider
% an  oriented ep groupoid $\Xx$  as a wnb groupoid with  $\la, \La$ identically equal to $1$. 
%(In this case, we take $k=1$ at each point $p\in |X|$.)\SSS
%Observe also that by Lemma~\ref{le:cts} below  
%every weighting $\La$ on an ep groupoid must be constant.
The groupoid $\Bb$ of Remark~\ref{rmk:I}(i)
has a weighting  in which 
%$\la= 1/2$ on each component of
%$B_0$ and 
$\La = 1$ on the image of $(0,1]$ in $|B|_{\Hsm}$ and $=1/2$ on the rest of $|B|_{\Hsm}$.  There is one branch point at the image $p_1$ of $1$, and the weighting condition 
implies that the weights satisfy Kirchoff's law at this branch point.
Moreover we may take the local branches at $p_1$ to be
the two components of $B_0$ each weighted by $m_1=m_2=1/2$.
Similarly,
the groupoid $\Zz$ defined in Example~\ref{ex:tear} has a weighting function in which  $\La = 1$ on the closure $cl(|D_-|_{\Hsm})$ of $|D_-|_{\Hsm}$
  and $\La=1/k$ elsewhere on $|Z|_{\Hsm}$.
 The branch locus is formed by the boundary of $|D_-|_{\Hsm}$. 
\SSS

\NI
(ii) If $(\Bb,\La)$ is compact and zero dimensional, then it is 
necessarily proper.  Hence $|B|_{\Hsm} = |B|$
consists of a finite number of equivalence classes $p$ of the relation $\sim$, each with 
a label $\La_p$.  Since the morphisms in $\Bb$ preserve orientation,
all the points $x\in B_0$ in the equivalence class $p$
have the same orientation.  Therefore, $p$ has a second label 
$\oo_p$ consisting of a sign $\oo_p=\pm$ that describes its orientation.\footnote
{
For further remarks about orientations, see
Remark~\ref{rmk:HWZ} and \S\ref{ss:euler}.}

  Thus the {\bf number of points} in $\Bb$ may be defined as
$\sum_{p\in |B|_{\Hsm}} \oo_p\La(p)$.  \SSS

\NI(iii)
 In both of the above cases
$\La$ is induced by a weighting 
function $\la$ defined on the components of $B_0$.  However, there are 
wnb groupoids for which such a function $\la$ is not uniquely defined by $\La$: see Remark~\ref{rmk:brorb0} (ii).
\end{example}

The following remark explains some of the technicalities of the above
 definition.

\begin{rmk}\labell{rmk:brorb}\rm 
(i) The local branches usually do not comprise
 the whole of the inverse image of $\pi_{\Hsm}^{-1}(N)$, but it is important that they map onto the full inverse image of $N$ in $|B|$.
 This condition, together with the fact that each branch injects into $|B|_{\Hsm}$, implies that the projection $|B|\to |B|_{\Hsm}$ is finite-to-one.
 \SSS
 
 \NI
(ii)  The local  regularity condition
 rules out the trivial example of 
an arbitrary  nonsingular sse groupoid $\Bb$ with 
connected object space $B_0$ (weighted by $1$) 
and with $\La \equiv 1$.  
Note that requiring the $U_i$ to inject into $|B|_{\Hsm}$ 
is considerably stronger than requiring that they inject into $|B|$.
 For example, consider the groupoid $\Bb_\phi$ in which
$B_0$ is the disjoint union of two copies 
 $R_1, R_2$ of $\R^2$ and the nonidentity 
morphisms  consist of two components each diffeomorphic to the open unit disc $D_1$ in $R_1$, one corresponding to a smooth embedding
 $\phi$ of $D_1$ onto a precompact subset of $R_2$
and the other corresponding to $\phi^{-1}$. This groupoid  can be weighted by setting $\La(p) = 1$ if $\pi_{\Hsm}^{-1}(p) $ intersects both $R_1$ and $R_2$ and $=1/2$ otherwise.
  Taking $N_p = |B|_{\Hsm}$ with local branches $R_1, R_2$ one finds that all the other 
conditions of Definition~\ref{def:brorb} 
are satisfied.  However the map $U_i\to |U_i|_{\Hsm}$ is injective
 iff $\phi$ extends to a homeomorphism of the closure
$cl(D_1)$ onto its image in $R_2$.  In this example, because $\phi(D_1)$ is precompact, it is enough to assume only  the  injectivity of
$U_i\to |U_i|_{\Hsm}$.  However, in general, to avoid pathologies 
we must assume that this map is a homeomorphism onto its image.
\SSS

\NI (iii) The other 
 part of regularity is that the branches have 
closed image in $N$.  This is an essential 
ingredient of the proof of Lemma~\ref{le:brstr}
and of the proof in Proposition~\ref{prop:approx} that
$\La$ is continuous on a dense open subset of $|B|_{\Hsm}$.
To see it at work,
consider the groupoid $\Bb$ whose objects are the disjoint union of the plane $\R^2$ and the unit disc (each weighted by $1/2$) with two components of nonidentity morphisms given by
the inclusion $D\subset \R^2$ and its inverse.  
Then $|B|_{\Hsm} = \R^2$ and all conditions except for properness are satisfied if
 we set $\La = 1$ on $D$ and $=1/2$ elsewhere.  Since $\Bb$ does not
 accord with our intuitive idea of a branched manifold (for example, if compactified to $S^2$ it does not carry a fundamental class), it is important to rule it out. \SSS
 
 \NI
 (iv)  Because we do not require the local branches to be
 connected, local branching structures $(N, U_i, m_i)$ behave well under restriction of $N$: if $(N, U_i,m_i)$ is a local branching structure at $p$
 then, for any other neighborhood $N'\subset N$ of $p$, the
 sets $U_i\cap \pi_{\Hsm}^{-1}(N')$ are the local branches 
 of a branching structure over $N'$.  Thus each point in $|B|_{\Hsm}$ has a neighborhood basis consisting of sets that support a local branching structure.  Note also that, because $|U_i|_{\Hsm}$ is assumed closed in $N$, $U_i$ is closed in $\pi_{\Hsm}^{-1}(N)$ and so is a union of components of 
 $\pi_{\Hsm}^{-1}(N)$.\SSS
 
 %%D added
 \NI (v)  We have chosen to impose rather few regularity conditions in order to make the definition of a wnb groupoid as simple and general as possible. However, in order for the quotient space $|B|_{\Hsm}$ to have a 
 reasonable smooth structure (so that, for instance, 
 one can integrate over it) one needs 
 more control over the morphisms; cf. the tameness conditions of Definition~\ref{def:tame}.  
  \end{rmk}

We now explain  the structure of the Hausdorff quotient $|B|_{\Hsm}$ for wnb groupoids. For  general spaces $Y$ it seems very hard to give a constructive definition of its maximal Hausdorff
quotient $Y_{\Hsm}$.  However, the covering and local 
regularity conditions are so strong that the quotient map $B_0\to |B|_{\Hsm}$ has the following explicit description.

 We define $\approx$ to be
the equivalence relation on $B_0$ generated by setting
$x\,\approx\, y$ if $|x|,|y|$ 
do not have disjoint 
open neighborhoods in $|B|$. In particular, if $x\sim y$ 
(i.e. $\pi(x) = \pi(y)\in |B|$) then $x\,\approx\,y$.

\begin{lemma}\labell{le:approx} Let $(\Bb, \La)$ be a 
 wnb groupoid. Then the fibers of $\pi_{\Hsm}: B_0\to |B|_{\Hsm}$ 
are the equivalence classes of $\;\approx$.
\end{lemma}
\begin{proof}   Since $|B|_{\Hsm}$ is Hausdorff, any 
two points in $|B|$ 
that do not have disjoint neighborhoods  must have the 
same image in $|B|_{\Hsm}$. Hence the equivalence classes of $\approx$ are contained in the fibers of $\pi_{\Hsm}$. 
Therefore there is a continuous surjection $B_{\sapprox}\to |B|_{\Hsm}$,
where $B_{\sapprox}$ denotes the  quotient of $B_0$ by 
the relation $\approx$, and it suffices to show that $B_{\sapprox}$ is Hausdorff. 
 Because of the covering property of the local branches,
 it suffices to work locally in subsets $W$ of $B_0$  of the form 
 $W: = \pi_{\Hsm}^{-1}(N) = \cup_{i=1}^\ell U_i$, 
where
$(N, U_i, m_i)$ is a
local branching structure.  We show below 
that for each such $W$, the quotient $W_{\sapprox}$ is Hausdorff.
This will complete the proof.

Consider any pair of distinct local branches
$U_i, U_j$ and the 
%%D corresponding 
set $M_{ji}: = \Mor(U_i,U_j)$ of morphisms 
from $U_i$ to $U_j$.  Since $U_i, U_j$ inject into $|B|$ 
both maps $s: M_{ji}\to U_i, t: M_{ji}\to U_j$ are injective and hence are diffeomorphisms onto their images.  Call these $V_{ji}\subset U_i$ and $V_{ij}\subset U_j$. 
%%D As in Remark~\ref{rmk:brorb} (ii),
%because $U_i, U_j$ inject into $|B|_{\Hsm}$ the  diffeomorphism  
Denote by $cl(V_{ji})$  the closure of $V_{ji}$ in $U_i$.
Since $|U_j|_{\Hsm}$ is relatively closed in $N\supset |U_i|_{\Hsm}\cup
|U_j|_{\Hsm}$, the set $\pi_{\Hsm}(V_{ji})$  is contained in $|U_j|_{\Hsm}$.
Hence, because $\pi_{\Hsm}$ is a homeomorphism on each $U_i$,
 the  diffeomorphism
$t\circ s^{-1}:V_{ji}\to V_{ij}$ extends to a homeomorphism $cl(V_{ji})\to cl(V_{ij})$.  
Define the $ji$ branch set as
$$
Br_{ji}: = \p(V_{ji}): = cl(V_{ji})\setminus V_{ji}\subset U_i.
$$
 Then $Br_{ji}$  is closed in $U_i$.
Further there is a homeomorphism $\phi_{ji}:Br_{ji}\to Br_{ij}$.

Now observe that if $x\in Br_{ji}$ then there is a unique 
$y\in Br_{ij}$ such
that $x\approx\, y$, namely $\phi_{ji}(x)$.  Conversely, 
if $x\in U_i, y\in U_j$ and $x\approx\, y$ then either $|x|=|y|$ or there are convergent sequences $x_n, y_n\in B_0$ with limits $x_\infty, y_\infty$
such that $x_n\sim y_n, x_\infty\sim x, y_\infty\sim y$.  
The morphism $\ga$ from $x_\infty$ to $x$ extends to a neighborhood of
$x_\infty$ and so transports (the tail of) the sequence $x_n$  to a sequence $x_n'\in U_i$ that converges  to $x$. Similarly, we may suppose that $y_n\in U_j$.  Hence 
%%D$x\in B_{ij}, y\in B_{ji}$.
$x\in Br_{ji}, y\in Br_{ij}$.
Let us write in this case that $x\approx_{ij} y$.
%%DWe also write $x\approx_{ij} y$ if $|x|=|y|$ in $|U_i|\cap |U_j|$.

The equivalence relation $\approx$ therefore has the following description on $W$.  Given
 $x\in U_i, z\in U_k$,  $x\approx z$ iff 
 either $|x|=|z|$ or
there is a finite sequence $i_1: = i, i_2,\dots, i_n = k$
of indices (with $n>1$) and elements 
%%D$x_j \in B_{i_j\,i_{j+1}}$ for $j=1,\dots,n-1$ such that
$x_j \in Br_{i_{j+1}i_j}$ for $j=1,\dots,n-1$ such that
 $$
 x_j\,\approx_{i_j\,i_{j+1}} \,x_{j+1},\quad\mbox{ for all } j.
 $$
 Note that we may assume that all the indices in this chain are different.  For, because the maps $\pi_{\Hsm}: U_i\to |B|_{\Hsm}$ are injective
 and constant on equivalence classes,
if $i_j = i_{j'}$ for some $j'>j$ then $x_j = x_{j'}$
so that the intermediate portion of the chain can be omitted.
It follows that the number of nonempty chains of this form is bounded 
by a number depending only on $\ell$, the number  of local branches.  
Moreover, for each such chain
$I: = i_1, \dots, i_n$ the map that takes its initial point to its endpoint is a homeomorphism $\phi_I$ from
a subset $X_I\subset Br_{i_n i_1}$ to a subset $Z_I\subset Br_{i_{1}\,i_{n}}$, where $\phi_I: = \phi_{i_{n}\,i_{n-1}} \circ\dots\circ \phi_{i_2\, i_1}$
is the restriction of $\phi_{i_n i_1}$. Note that $X_I$ and $Z_I$ are closed in $U_{i_1}$ and $U_{i_n}$ respectively.   
For each $i\ne k$, let  $C_{ki}$ be the set of chains from $i$ to $k$, and consider the set
$$
X_{ki}': = \Bigl(\bigl(\pi^{-1}(|U_{k}|)\bigr)\cap U_i\Bigr)\;\cup\; \bigcup_{I\in C_{ki}} X_I 
$$ 
of all points in $U_i$ that are equivalent to a point in $U_k$. 
The above remarks imply that this is closed in $U_i$. 

We now return to the quotient $W_{\sapprox}$. To see this is Hausdorff, note first that each equivalence class $\bx$
 contains at most one element from each $U_i$.  Hence  we may write $\bx: = ({x_i})_{i\in J_\bx}$, where 
$x_i\in U_i$ and $J_\bx\subset \{1,\dots,\ell\}$. 
  Suppose that $\bx\ne \by$.
Then $x_i\ne y_i$ for all $i$.  We 
construct disjoint $\approx$-saturated\footnote
{
A set is said to be $\approx$-saturated if it is a 
union of equivalence classes.}
 neighborhoods $V_\bx, V_\by$ as follows.

Suppose first that  $J_\bx\cap J_\by\ne \emptyset.$ By renumbering 
we may suppose that
 $1\in J_\bx\cap J_\by$.  
Choose disjoint open neighborhoods $V_{x1}, V_{y1}$
 of  $x_1, y_1$ in $U_1$
and define $V_{x1}^S, V_{y1}^S$ 
to be their saturations under $\approx$.
Then $V_{x1}^S\cap V_{y1}^S=\emptyset$. To see this note that any point
 $z\in V_{x1}^S\cap V_{y1}^S$ is equivalent to some point in $V_{x1}$ 
 and some point in $V_{y1}$. Since $z$ is equivalent to at most one
  point in $U_1$ this implies that $V_{x1}\cap V_{y1}\ne \emptyset$, a contradiction. 
  
  We now enlarge $V_{x1}^S, V_{y1}^S$ to make them open.
   To this end, consider the smallest integer $i$ such that
  either $V_{x1}^S\cap U_i$ or  $V_{y1}^S\cap U_i$ is not open.
 Then $i>1$ by construction,  and  for every $j<i$
  the disjoint sets $V_{x1}^S\cap X_{ji}'$,
  $V_{y1}^S\cap X_{ji}'$ are relatively open in the closed set $X_{ji}'\subset U_i$.  Therefore
  we may construct open disjoint neighborhoods $U_{xi}\subset U_i$ of $V_{x1}^S\cap U_i$ 
  and $U_{yi}\subset U_i$ of
  $V_{y1}^S\cap U_i$  by adding points
  in $U_i\setminus \bigl(\cup_{j<i} X_{ji}'\bigr)$.  Now 
  consider the sets $V_{xi}: = V_{x1}^S\cup
  V_{xi}^S$ and $V_{yi}: = V_{y1}^S\cup
  V_{yi}^S$, where $A^S$ denotes the saturation of $A$.
   These sets are disjoint as before.  Moreover, their intersections with the $U_j, j\le i,$ are open.   Hence after a finite number of similar steps
  we find suitable disjoint open $V_\bx: = V_{xn}$ and 
  $V_\by: = V_{yn}$.

Next suppose  that $J_\bx$ and $J_\by$ are disjoint.  If $i\in J_\bx$ then $x_i\notin X_{ji}'$ for any $j\notin J_\bx$.
 Therefore the set  
$$
V_\bx: = \bigcup_{i\in J_\bx} \Bigl(U_i\setminus \bigl(\cup_{j\notin J_\bx} X_{ji}'\bigr)\Bigr)
$$
is an open neighborhood of $\bx$ in $W = \cup_iU_i$.
Since $V_\bx$ is
 saturated under $\approx$ by construction, it projects to
an open neighborhood of $\bx$ in $W_{\sapprox}$.  Finally note that
$V_\bx$ is disjoint from the similarly defined open set
$$
V_\by: = \bigcup_{i\in J_\by} \Bigl(U_i\setminus \bigl(\cup_{j\notin J_\by} X_{ji}'\bigr)\Bigr).
$$
This completes the proof.
\end{proof}

\begin{prop}\labell{prop:approx}
%%D{\rm (i)} (later on changed numbers without note)
%Each local branch $U_i$ is homeomorphic to its image $|U_i|_{\Hsm}$ in $|B|_{\Hsm}$.\SSS

\NI{\rm (i)} For all $p\in |B|_{\Hsm}$, any open neighborhood in $|B|$ of the fiber over $p$ contains a saturated open neighborhood $|W|$.\SSS

\NI{\rm (ii)}
 $|B|_{\Hsm}$ is second countable and locally compact.
%It is path connected if $|B|_{\Hsm}$ may be covered by 
%open sets $N_p, p\in P,$ that each support local 
% branch structures whose branches $|U_i^p|_{\Hsm}$ all contain $p$.
 \SSS

\NI
 {\rm (iii)} 
 The branch points form a closed and nowhere dense subset of $|B|_{\Hsm}$.\SSS

\NI
{\rm (iv)}  The weighting function $\La$ is locally constant except possibly at branch points.
\end{prop}

\begin{proof}
%  We saw in the proof of the  lemma how to extend an open subset $V$ of $U_1$ to an open and saturated subset $V^S$ of $B_0$ by adding points from the sets 
%$U_i\setminus \bigl(\cup_{j<i}X_{ij}'\bigr)$ in turn. Then $V^S\cap U_1 = V$.  Hence the restriction of $\pi_{\Hsm}$ to $U_1$ is an open map.  (i) readily follows.
To prove (i) let $(N_p, U_i, m_i)$ be a local branching structure at $p$, and $|V|
\subset |B|$ be any open neighborhood of the fiber at $p$.  Then 
$|V_i|: = |V|\cap |U_i|$ is open in $|U_i|$.  If $V_i$ denotes the corresponding open subset of $U_i$, the set
$|W|: = \cap_i |V_i^S|$ satisfies the requirements.

Each open set $N$ that supports a local branching structure is a finite union of locally compact closed sets $|U_i|_{\Hsm}$, and hence is locally compact.  Moreover,  the $|U_i|_{\Hsm}$ are second 
countable in the induced topology.  Hence so is $N$.  Since we assumed $B_0$ has a countable dense subset, the same is true for $|B|_{\Hsm}$.  Therefore, $|B|_{\Hsm}$ is the union of countable many open sets $N_i$ and so is itself second countable.  This proves (ii).

Denote by $|Br|_{\Hsm}$ the set of branch points in $|B|_{\Hsm}$. To prove (iii) it suffices to show that for any $N$ as above
the intersection 
$|Br|_{\Hsm}\cap N$ is closed and nowhere dense in $N$.
It follows from the proof of Lemma~\ref{le:approx} that
$$
|Br|_{\Hsm}\cap |U_i|_{\Hsm} = \pi_{\Hsm}(\cup_{j\ne i}Br_{ij}).
$$
We saw earlier that for each $i$ the set of branch points 
$\cup_{j\ne i}Br_{ij}$ in $U_i$ is relatively closed.  It is nowhere dense by construction. 
Since $|U_i|_{\Hsm}$ is closed in $N$, and $U_i$ is homeomorphic to $|U_i|_{\Hsm}$, $ \pi_{\Hsm}(\cup_{j\ne i}Br_{ij})
$
is closed in $N$ for all $i$.  Since there are a finite number of local branches, (iii) holds. 

Consider $|Br|: = |\pi|_{\Hsm}^{-1}(|Br|_{\Hsm})$, the set of points in $|B|$
on which $|\pi|_{\Hsm}$ is not injective, and denote by $|N|$ the open set $|\pi|_{\Hsm}^{-1}(N)\subset |B|$.
(iv) will follow if we show that for each connected component 
$|V|$ of $|N|\setminus |Br|$ and each local branch $U_i$ over $N$ either
$|U_i|\cap |V|=\emptyset$ or $|V|\subset |U_i|$.  But $|U_i|$ is open in $|N|$ by Lemma~\ref{le:Haus}. Its intersection with $|N|\setminus |Br|$ is also closed since it is the inverse image of the relatively closed subset 
$|U_i|_{\Hsm}\setminus |Br|_{\Hsm}$ of $|N|\setminus |Br|$.  Hence $|U_i|\setminus |Br|$ is a union of components of $|N|\setminus |Br|$, as required.
\end{proof}

%%%%%%%%%%%%%%%%%%%%%%%%%%%%%%%%%%%%%%%%%%%%%%%%%%%%%%%%
\subsection{Layered coverings}
%%%%%%%%%%%%%%%%%%%%%%%%%%%%%%%%%%%%%%%%%%%%%%%%%%%%%%%%

There are two useful kinds of functors for wnb groupoids, 
those that induce 
homeomorphisms on the orbit space $|B|_{\Hsm}$ and those 
that simply induce surjections on the orbit space.  In the first case we require $\La$ to be preserved while in the second we expect the induced map $|F|_{\Hsm}:|B'|_{\Hsm}\to |B|_{\Hsm}$ to push $\La'$ forward to $\La$. In other words we expect the identity  
\begin{equation}\labell{eq:Fbr}
(|F|_{\Hsm})_*(\La')(p): = \sum_{q: |F|_{\Hsm}(q) = p}\,\La'(q) = \La(p) 
\end{equation}
to hold at all points $p\in |B|_{\Hsm}$.

\begin{defn}\labell{def:brref} Let  $(\Bb,\La)$ and $(\Bb',\La')$ be wnb groupoids.  A refinement $F: \Bb'\to \Bb$ is said to be {\bf weighted}  if   $\La\circ |F|_{\Hsm} = \La'$. \SSS

\NI
A smooth functor $F: \Bb'\to \Bb$ is said to be a {\bf layered covering}
 if 
 \SSS
 
\NI {\bf (Covering)}  $F$ is a local diffeomorphism on objects and induces a surjection $|F|: |B'|\to |B|$,\SSS
 
  \NI {\bf (Properness)} the induced map
 $|F|_{\Hsm}: |B'|_{\Hsm}\to |B|_{\Hsm}$ is proper.\SSS

 \NI
 {\bf (Weighting)}  $(|F|_{\Hsm})_*(\La') = \La$.\SSS
 
 \NI
 Two  wnb groupoids are {\bf equivalent} if they have a common weighted refinement.  They are {\bf commensurate} if they have a common layered covering. Finally, two compact wnb groupoids 
 $(\Bb, \La)$ and $(\Bb',\La')$ (without boundary)
  are {\bf cobordant} if there is a compact $(d+1)$-dimensional wnb groupoid with boundary
 $-(\Bb,\La)\sqcup (\Bb',\La')$.
  \end{defn}

\begin{figure}[htbp] %  figure placement: here, top, bottom, or page
 %xx  \centerline{\psfig{figure=orbfig6.jpg,width=4in}}
  \centering
  \includegraphics[width=4in]{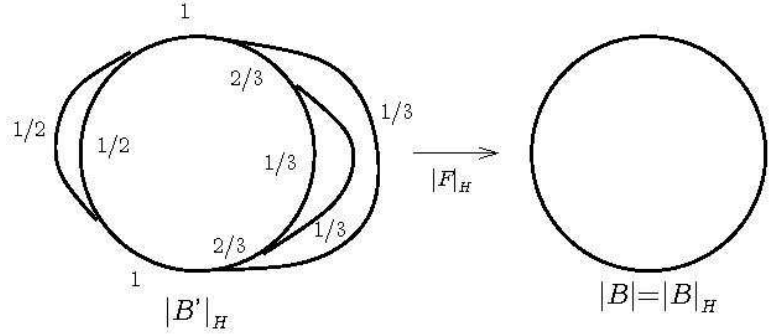} 
   \caption{A layered covering of the circle; $|B|$ has $\La\equiv 1$, the values of $\La'$ are marked.}
   \label{fig6}
\end{figure}

\begin{rmk}\labell{rmk:brorb0} \rm (i)
 Any wnb groupoid $(\Bb,\La)$ is equivalent to a
wnb groupoid $(\Bb',\La')$  in which all 
local branches needed to describe its branching structure
are unions  of components of $B_0'$.
To see this, choose a locally finite cover $N_p,  i\in A_p,$
of $|B|_{\Hsm}$ by sets that support local branching structures 
and let $U_i^p, i\in A_p$, be the corresponding set of local branches.
Then define $\Bb'$ to be the  refinement of $\Bb$
with objects $\sqcup_{i,p} U_i^p$ as in Remark~\ref{rmk:refin}.
Then $|B'|=|B|$ and $|B'|_{\Hsm} = |B|_{\Hsm}$, and so we may define $\La': = \La$.\SSS

\NI
(ii)  The condition that $|F|:|B'|\to |B|$ is surjective 
does not follow from the other conditions for a layered covering.  
Consider for example the wnb groupoid $\Cc$ with objects two copies of $S^1$ 
and morphisms $C_1= C_0\cup A^\pm$, where $A^+: = S^1\setminus \{0\}$ identifies one circle with the other except over $0$
 and 
$A^-: = (A^+)^{-1}$.  Then $|C|_{\Hsm} = S^1$ and we give it weight $\La\equiv 1$.  
The functor $F:\Cc\to \Cc$ that identifies both copies of $S^1$
 to 
the same component of $C_0$ satisfies all  conditions for a layered covering except that the induced map $|C|\to |C|$ is not surjective.
This wnb groupoid also illustrates the fact that 
the local weights $m_i$ 
on the local branches $U_i$ need not be uniquely 
determined by $\La$, and hence that $\La$ may not lift to a well defined function on $|C|$.\SSS

%\NI
%(iii)  The fact that in the last example $\La$ is not well defined on $|C|$  occurs because 
%in this case 
%any local branching structure must contain different local branches whose images in $|C|_{\Hsm}$ coincide. At the opposite extreme, suppose 
%$(\Bb',\La')$ is a wnb groupoid  as described in (i)  in which all
%local branches needed to describe the branching structure
%are unions  of components of $B_0'$.  Then we may define an auxiliary weighting function $\la:\pi_0(B_0')\to (0,\infty)$ by setting
%$\la([C]): = m_i^p$ where the component $C$ is contained in $U_i^p$.

%
%\SSS

\NI
(iii)
If $F$ is a layered covering then although the induced map $|F|: |B'|\to |B|$ is surjective,
the induced map on objects need not be surjective.
Also, the properness requirement is important.  Otherwise,
given any wnb groupoid $(\Bb, \La)$ define 
$(\Bb',\La')$ as in (i) above and then consider
$(\Bb'',\La'')$, where $\Bb''$ 
has the same objects as $\Bb'$ but only identity morphisms.
Then $|B''|_{\Hsm} = |B''| = \sqcup_{i,p} U_i^p$ and we may define $\La'': = m_{i,p}$ on $U_i^p$.
The inclusion $\Bb''\to \Bb$ satisfies all
the conditions for a layered covering except for properness.
Since we want commensurate wnb groupoids to have the same fundamental class, we cannot allow this behavior.
\end{rmk}

Our next aim is to  show that
layered coverings have  the expected functorial properties. In particular, commensurability is an equivalence relation.
For this we need some preparatory lemmas. 

 We shall say that a local branching structure 
%%D $(N, U_i, m_i),_{i\in I},$ at $p$
 $(N, U_i, m_i)_{i\in I}$ at $p$
 is {\bf minimal} if the fiber 
 $ (|\pi|_{\Hsm})^{-1}(p)\subset |B|$ over $p$ is a collection of distinct points
$|x_i|, i\in I$, where $x_i\in U_i$.  Thus in this case there is a bijective correspondence between the local branches and the points of the fiber.

 \begin{lemma}  Let $(\Bb,\La)$ be a wnb groupoid. Then every point $p\in |B|_{\Hsm}$ has a minimal local branching structure.
 \end{lemma}
 \begin{proof}  Choose any local branching structure
 $(N, U_i, m_i)$ at $p$.  For each point $w_\al\in (|\pi|_{\Hsm})^{-1}(p)\subset |B|$,
 let $I_\al$ be the set of indices $i$ such that $w_\al\in |U_i|$.  
 For each such $\al$, choose one element $i_\al\in I_\al$, and define
 $$
 m_\al: = \sum_{i\in I_\al} m_i,\qquad |V_\al|: = 
 \bigcap_{i\in I_\al} |U_i|.
 $$
Each $|V_\al|$ is open.  Hence $\cup_\al |V_\al|$ is an open neighborhood of the fiber
 $(|\pi|_{\Hsm})^{-1}(p)$ and so  by Proposition~\ref{prop:approx}
  it contains a saturated open neighborhood
 $|W|: = (|\pi|_{\Hsm})^{-1}(N')$.  
Define $U_\al': = U_{i_\al}\cap \pi^{-1}(|W|)$.  Then
 $(N', U_\al', m_\al)$ is a minimal local branching structure
 at $p$.
 \end{proof}

    We shall say that two (possibly disconnected) open subsets 
$U_0, U_1$ of the space of objects $X_0$ of an sse groupoid $\Xx$ are 
$\Xx$-{\bf diffeomorphic} if there is a subset $C\subset X_1$ 
of $s^{-1}(U_0)$ such that the maps $s:C\to U_0$ and 
$t:C\to U_1$ are both diffeomorphisms.  In this situation
we also say that there is a diffeomorphism $\phi: U_0\to U_1$ in $\Xx$. 
If the $U_i$ both inject into $|X|$, there is such $\phi$ iff $|U_0| =
|U_1|\in |X|$.
 
 \begin{lemma}\labell{le:brstr}  Suppose that $F:\Aa\to \Bb$ is a layered covering.
 Let $p\in |B|_{\Hsm}$ and denote the points in $(|F|_{\Hsm})^{-1}(p)$ by
 $q_\al, 1\le \al\le k$.  Then there are  minimal
  local branching structures 
  $$
  (N_p, U_i^p, m_i^p)_{i\in I}\; \mbox{ at } p,\qquad
(N_\al, U_j^\al, m_j^\al)_{j\in J_\al}\; \mbox{ at } q_\al,
$$
 %at $p$ and the points $q_\al\in |A|_{\Hsm}$ such that $|F|_{\Hsm}(q_\al) = p$,
 that are compatible  in the sense that each 
  $F(U_j^\al)$
 is $\Bb$-diffeomorphic to some local branch $U_{i_j}^p$. Moreover $U_{i_j}^p$
 is unique.
% Moreover, we may suppose that $(N_p, U_i^p, m_i)$ is minimal.
 \end{lemma}
 \begin{proof}  
 Let $\{q_1,\dots,q_m\}= (|F|_{\Hsm})^{-1}(p)\subset |A|_{\Hsm}$.
 Choose minimal local branching structures $(N_p, U_i^p, m_i)$ at $p$
 and $(N_\al, U_j^\al,m_{j\al})$ at $q_\al$ for  $1\le\al\le m$.  
 Since $|A|_{\Hsm}$ is
 Hausdorff, we may suppose that the $N_\al$ are pairwise disjoint. 
 Moreover, because $|F|_{\Hsm}$ is proper, the union $\cup_\al
 |F|_{\Hsm}(N_\al)$ is a neighborhood of $p$ in $|B|_{\Hsm}$. 
 (Otherwise, there would be a sequence of points $q_n\in |A|_{\Hsm}$
 lying outside $\cup_\al N_\al$ but such that $|F|_{\Hsm}(q_n)$
 converges to $p$, which contradicts properness.)
%By shrinking $N_p$ we may arrange that each 
%intersection $|F|_{\Hsm}(N_\al)\cap N_p$
%is a relatively closed subset of $N_p$.

Fix $\al$, and for each $j\in J_\al$ denote by $x_j^\al$ the point
in  $U_j^\al$ that  projects to $q_\al$; this exists by minimality.
Since $(|\pi|_{\Hsm})^{-1}(N_p) = \cup_i |U_i^p|$,
for each branch $U_j^\al$ there is a morphism $\ga_\al\in B_1$ with
source $F(x_j^\al)$ and target in some local branch, say
$U_{i_j}^{p,\al}$,
% where $i_j\in I$, 
 at $p$.  By minimality at $p$, the index $i_j$ is unique.
  This morphism extends to a diffeomorphism $\phi_\al$ from an
 open subset $F(V_j^\al)$ of $F(U_j^\al)$ onto an open subset
 $V_{i_j}^p$ of $U_{i_j}^p$.  Since $\cup_{j\in J_\al} |U_j^\al|$ is a
 neighborhood in $|A|$ of the fiber over $q_\al$, so is the open set
 $\cup_{j\in J_\al} |V_j^\al|$.  Hence, by
 Proposition~\ref{prop:approx}(i), there is an open neighborhood $N_\al'$
 of $q_\al$ such that $|\pi|_{\Hsm}^{-1}(N_\al')\subset \cup_{j\in
 J_\al} |V_j^\al|$.  By shrinking the sets $V_j^\al$ we may therefore
 suppose that $|\pi|_{\Hsm}^{-1}(N_\al')=\cup_j |V_j^\al|$, i.e. that
 $(N_\al', V_j^\al, m_{j\al})$ is a local branching structure at
 $q_\al$.  We also shrink the $V_{i_j}^{p,\al}$ so that they remain
 $\Bb$-diffeomorphic to the sets $F(V_j^\al)$, $j\in J_{\al}$.

Now observe that because $|F|:|A|\to |B|$ is surjective, 
$$
(|\pi|_{\Hsm})^{-1}(p)\;\;\subset\;\; |W|: = 
\bigcap_{\al,\;j\in J_\al}
|V_{i_j}^{p,\al}| \;\;\subset\;\; |B|.
$$
%$$
%(|\pi|_{\Hsm})^{-1}(p)\;\;\subset\;\; |W|: = \bigcup_{w\in
%(|\pi|_{\Hsm})^{-1}(p)}\left(\bigcap_{\al:|F|(x_j^{\al})=w,\;j\in J_\al}
%|V_{i_j}^{p,\al}|\right) \;\;\subset\;\; |B|.
%$$
%contains
%the fiber $(|\pi|_{\Hsm})^{-1}(p)$ over $p$.  
The set $|W|$ is open since both $F$ and
  the map $B_0\to |B|$ are open.  Hence it contains a set of the form 
  $|\pi_{\Hsm}|^{-1}(N_p')$ for some neighborhood $N_p'$ of $p$.
  Shrinking the sets $N_\al'$ and $V_j^\al$
   further as necessary, we may suppose that
 $|F|_{\Hsm}^{-1}(N_p') = \cup_\al N_\al'$.   
  Then the local branching structures $(N_\al', V_j^\al, m_{j\al})$ and
  $(N_p', V_i^p\cap (\pi_{\Hsm})^{-1}N_p', m_i)$ 
  satisfy all requirements.
\end{proof}

\begin{prop}\labell{prop:commen}  Let $\Aa, \Bb,\Cc$ be wnb groupoids.
Then:\SSS

\NI {\rm (i)}  If $F:\Aa\to \Bb$ and $G:\Bb\to \Cc$ are layered coverings so is the composite $G\circ F$.\SSS

\NI
{\rm (ii)}
If $F:\Aa\to \Bb$ and $G:\Cc\to \Bb$ are layered coverings,
the (weak) fiber product $\Zz: = \Aa\times_\Bb\Cc$ of $\Aa$ and $\Cc$ over
$\Bb$ is a wnb groupoid and the induced functors
$\Zz\to \Bb, \Zz\to \Aa$ are layered coverings.\SSS 

\NI
{\rm (iii)}  Any two compact and commensurate wnb groupoids are cobordant through a 
wnb groupoid.
\end{prop}

\begin{proof}  The proof of (i) is straightforward, and is left to the reader.  As for (ii), observe first that $\Zz$ is a nonsingular
 sse Lie groupoid with objects $Z_0$ contained in the product
$A_0\times B_1\times C_0$.\footnote
{ 
Since $\Bb$ is nonsingular, we may identify $Z_0$ with the 
strict fiber product $A_0\times_{|B|} C_0$.  However,  
later on we will apply this construction to certain nonsingular $\Bb$,
and so it is convenient to retain more general language here.}
%
%Hence we define 
%$\la_Z$ on $\pi_0(Z_0)$ as the composite 
%$$
%\la_Z: \pi_0(Z_0)\to  \pi_0(A_0)\times \pi_0(C_0) \stackrel{\la_A\times \la_C} \to [0,\infty),
%$$
%where $\la_A\!\times\!\la_C ([U]\!\times \![V]): = \la_A([U])\cdot\la_C([V]).$
The orbit space $|Z|$ is the (strict) fiber product  $|A|\times_{|B|}|C|$.  It follows easily that
$|Z|_{\Hsm}$ is the fiber product $|A|_{\Hsm}\times_{|B|_{\Hsm}} |C|_{\Hsm}$, since the latter space has the requisite universal property.
Hence $|Z|_{\Hsm} = \{(a,c): |F|_{\Hsm}(a) = 
|G|_{\Hsm}(c)\in |B|_{\Hsm}\}$ and we set
$$
\La_Z (a,c) : = \frac{\La_A(a)\, \La_C(c)}{\La_B\bigl(|F|_{\Hsm}(a)\bigr)}. 
%: |Z|_{\Hsm}\;\to\;  |A|_{\Hsm}\times |C|_{\Hsm} \;\stackrel{\La_A\times \La_C} \to\;[0,\infty).
$$
Observe that the projections of $|Z|_{\Hsm}$ to  $|A|_{\Hsm}$ and $|C|_{\Hsm}$ push 
$\La_Z$ forward to $\La_A$ and $ \La_C$ respectively.  
 Hence (ii) will follow once we show that $\Zz$ has the requisite  local branching structures.

Given $z=(a,c)\in |Z|_{\Hsm}$, let $p = |F|_{\Hsm}(a) = |G|_{\Hsm}(c)\in |B|_{\Hsm}$. 
Applying Lemma~\ref{le:brstr}
 first to $\Aa\to \Bb$ and then to $\Cc\to\Bb$ 
and then  restricting as in Remark~\ref{rmk:brorb}(iv), we
can find minimal local branching structures 
 $(N_a, U_j^a, m_{ja})_{j\in J}$, $(N_c, U_k^c,m_{kc})_{k\in K}$ and $(N_p, U_i^p, m_i)$ satisfying the following compatibility conditions: for each $j, k$
 there are unique $i_j, i_k$ and local diffeomorphisms $\ga_j, \ga_k$ in $\Bb$ such that
 $$
F(U_j^a) =\ga_j(U_{i_j}^p),\qquad \ga_k(G(U_k^c)) =U_{i_k}^p.
 $$

We now show that there is a local branching structure over 
$N_z: = N_a\times_{|B|_{\Hsm}} N_c$.    
By  minimality, for each point $(w_a, w_c) \in |A|\times_{|B|}|C|$,
there is a unique pair of local branches $(U_j^a, U_k^c)$ such that
$w_a\in |U_j^a|, w_c\in |U_k^c|$.  Since $|F|(w_a) = |G|(w_c)\in |B|$,
the corresponding indices $i_j, i_k$ coincide.  Hence there is an open set 
$$
U_{jk}^{ac}: = \{(x,\ga,y): x\in U_j^a, \ga = \ga_j\circ\ga_k|_y, y\in U_k^c\}\subset Z_0.
$$
It remains to observe that the collection of all such sets, with weights 
$$
m_{jk}^{ac}: = \frac{m_{ja}\, m_{kc}}{m_{i_j}},
$$
forms a local branching structure over 
 $N_z$. This proves (ii).

To prove (iii), it suffices by (ii) to consider the case when
there is a layered covering $F:\Bb'\to \Bb$.  Define  $\Ww$ by setting $W_0: = B_0'\times [0,3/4)\; \sqcup \; B_0\times (1/4,1]$ and 
$$
W_1: = \Bigl(B_1'\times [0,3/4)\Bigr) \;\sqcup \; \Bigl(B_1\times (1/4,1]\Bigr)\;\sqcup \;\Bigl((B_0')_\pm\times (1/4,3/4)\Bigr),
$$
where $(x,t)_+\in (B_0')_+\times (1/4,3/4)$ denotes a morphism from $(x,t)\in B_0'\times (1/4,3/4)$ to $(F(x),t)\in B_0\times (1/4,3/4)$
and $(x,t)_-$ is its inverse. The nonsingularity of $\Bb',\Bb$ implies that $\Ww$ is a nonsingular groupoid in the category of oriented manifolds with boundary.  In fact, its boundary $\p\Ww$ projects to $\p \Ii$ under the obvious projection  $\Ww\to \Ii$, where $\Ii$ is the category with objects $[0,1]$ and only identity morphisms.
Further $\p\Ww$ may obviously be identified with $-\Bb'\sqcup \Bb$. 
Therefore $\Ww$ is a cobordism from  $\Bb'$ to $\Bb$.

It remains to check that $\Ww$ is a wnb groupoid.  Since
$$
|W|_{\Hsm} \cong \Bigl(|B'|_{\Hsm}\times [0,1/4)\Bigr)\cup 
\Bigl(|B|_{\Hsm}\times [1/4,1]\Bigr),
$$
we define $\La_W$ to be the pullback of 
$\La_B'$ on $|B'|_{\Hsm}\times [0,1/4)$ and the pullback of $\La_B$ on 
$|B|_{\Hsm}\times [1/4,1]$.  It is obvious that $(\Ww,\La)$ has local
branching structures at all points of $|W|_{\Hsm}$ except possibly for
$p\in |B|_{\Hsm}\times \{1/4\}$.  But here one constructs suitable
local branches as in the proof of (ii), using
Lemma~\ref{le:brstr} as before.
%
% if we set $\la_W: = \la_{B'}$ on 
%$B_0'\times [0,3/4)$ and $: = \la_{B}$ on 
%$B_0\times [0,3/4)$, then this pushes forward to a weighting function $\La_W$ on $|W|_{\Hsm}$ that restricts to $\La_{B'}, \La_B$ on $\p(|W|_{\Hsm})$. 
\end{proof}

%%%%%%%%%%%%%%%%%%%%%%%%%%%%%%%%%%%%%%%%%%%%%%%%%%%%%%%%%%%%%%%%
\subsection{Branched manifolds and resolutions}
%%%%%%%%%%%%%%%%%%%%%%%%%%%%%%%%%%%%%%%%%%%%%%%%%%%%%%%%%%%%%%%%

\begin{defn}\labell{def:brman}
A {\bf branched manifold structure} on a topological 
space $Z$ consists of a wnb groupoid
$(\Bb,\La)$ together with a 
homeomorphism $f:|B|_{\Hsm}\to Z$. 
Two such structures $(\Bb,\La, f)$ and
$(\Bb',\La', f')$ are {\bf equivalent} if they have a 
common weighted equivalence, i.e. if there is a third structure 
$(\Bb'',\La'', f'')$ and weighted equivalences 
$F:\Bb''\to\Bb, F':\Bb''\to \Bb'$ such that
$f'' = f\circ |F|_{\Hsm} = f'\circ |F'|_{\Hsm}$.
%Note that equivalent structures induce the same function $\La_Y:Y\to 
%(0,\infty)$.  Therefore we define 

A {\bf branched manifold} $(\un Z, \La_Z)$ is a
pair, consisting of 
%%Da second countable paracompact Hausdorff 
topological space $Z$ together with a function $\La_Z: Z\to (0,\infty)$, in which $Z$ is equipped with an equivalence class of 
branched manifold structures that induce the function $\La_Z$.

Two $d$-dimensional branched manifolds $(\un Z, \La_Z)$ and 
$(\un Z', \La_Z')$ are {\bf commensurate} if they have commensurate branched manifold structures. They are {\bf cobordant} if there is a 
$(d+1)$-dimensional branched manifold $(\un W, \La_W)$ in the 
category of smooth manifolds with boundary whose
 (oriented) boundary decomposes into the disjoint union $(-\un Z, \La_Z)\sqcup (\un Z', \La_Z')$.  
\end{defn}

Proposition~\ref{prop:commen} implies that
if $(\un Z, \La_Z)$ and $(\un Z', \La_Z')$ are commensurate any pair of 
branched manifold structures on them have a common layered covering.
Further any pair of commensurate branched manifolds are 
cobordant via a branched manifold.

We now consider maps from branched manifolds to orbifolds.
These are induced by smooth functors 
$F:(\Bb, \La)\to \Xx$ where $\Xx$ is an ep groupoid.
It is convenient to consider $\Xx$ as a weighted but {\it singular} branched groupoid with weighting function that is identically equal to $1$.
In this case $|X|_{\Hsm} = |X|$ and every point $p$ has a local branch 
structure with one branch, namely any open set 
containing a point $x\in \pi^{-1}(p)$.  This groupoid satisfies
all conditions of Definition~\ref{def:brorb}, except that it is singular and
  the local branches do not inject into $|X|$.

The following definition is motivated by the example of the inclusion 
$F:\Zz\to \Xx$ of Example~\ref{ex:tear}, where $\Xx$ is the teardrop groupoid
and $(\Zz,\La_Z)$ is as in Example~\ref{ex:brorb}(i). We want 
to consider this as some kind of equivalence, but note that $F$ 
does not push $\La_Z$ forward to $\La_X\equiv 1$ at all points.  Rather
$(|F|_{\Hsm})_*(\La')(p)=1$ except at the singular point of $|X|$
where it equals $1/k$.  The fact that this is the unique nonsmooth
point in the sense of Definition~\ref{def:sm} motivates the following definition.

 \begin{defn}\labell{def251}  Let $(\Bb, \La)$ be a wnb groupoid 
 and $\Xx$ an ep groupoid.  A  functor $F:(\Bb, \La)\to \Xx$ is 
 said to be a {\bf resolution}
of $\Xx$ if the following conditions hold:\SSS

\NI {\bf (Covering)}  $F$ is a local diffeomorphism on objects,
\SSS
 
  \NI {\bf (Properness)} the induced map
 $|F|_{\Hsm}: |B|_{\Hsm}\to |X|$ is proper.\SSS

 \NI
 {\bf (Weighting)}  $(|F|_{\Hsm})_*(\La') = 1$ at all smooth points of $|X|$.\SSS

\NI Similarly, a map 
 $\un\phi:(\un Z,\La_Z)\to (\un Y,\La) $ from a  branched manifold 
 to an orbifold is  called 
 a {\bf resolution of $(\un Y,\La)$} if it is induced by a resolution 
 $F:(\Bb,\La)\to \Xx$ where $(\Bb, \La, f)$ is a branched 
 manifold structure on $(\un Z,\La_Z)$ and $(\Xx, f')$ is 
 an orbifold structure on  $(\un Y,\La)$.  Here we require that
 the diagram
 $$
 \begin{array}{ccc} |B|_{\Hsm} &\stackrel{|F|_{\Hsm}}\longrightarrow&
     |X|\vspace{.08in}\\
 f\downarrow\;\;\;\; &&  f'\downarrow\;\;\;\;\\
 Z &\stackrel{\phi}\longrightarrow & Y
 \end{array}
 $$
 commute.
 \end{defn}

Note that for any resolution  $F:\Bb\to \Xx$ 
the induced map $|F|_{\Hsm}: |B|_{\Hsm}\to |X|\equiv |X|_{\Hsm}$ is surjective; its
image is closed by the properness assumption, and is dense by the
weighting property and the fact that the smooth points are dense in
$|X|$.  Moreover, the following analog of Lemma~\ref{le:brstr} holds.

\begin{lemma}\labell{le:brstr1}
For any resolution  $F:\Bb\to \Xx$ each
$x\in X_0$ has a neighborhood $U$ that is evenly covered by $F$ in the
sense that for all $q_\al\in |F|_{\Hsm}^{-1}(|x|)$ there are local
branching structures $(N_\al, U_j^\al, m_j^\al)_{j\in J_\al}$ at
$q_\al$ such that $|F|_{\Hsm}(|U_j^\al|) = |U|$ for all $\al, j\in
J_\al$.  Moreover each set $F(U_j^\al)$ is $\Xx$-diffeomorphic to $U$.
 \end{lemma}
\begin{proof}  Since $|F|_{\Hsm}:|B|_{\Hsm}\to |X|$ is proper, it is finite to one, 
    and also open (because $|X|$ is locally compact and normal).
In particular, there are a finite number of points $q_\al$.  
For each $q_\al$ choose $N_\al'$ so that it supports a local branching 
structure at $q_\al$.  Choose 
a connected neighborhood $U$ of $x$ so that $|U|\equiv |U|_{\Hsm}\subset \cap_\al \bigl(|F|_{\Hsm}(N_\al')\bigr)$,
and then define
$N_\al: = N_\al'\cap (|F|_{\Hsm})^{-1}(|U|)$. Then $|F|_{\Hsm}(|U_j^\al|_{\Hsm})$ is closed
in $|U|$ since $|U_j^\al|_{\Hsm}$ is closed in $N_\al$ and $|F|_{\Hsm}$ is proper.
But $|F|_{\Hsm}(|U_j^\al|_{\Hsm})$ is also open in $|U|$,
since it is the image of $U_j^\al$ under the composite of the two open maps $F:B_0\to X_0$ and $\pi:X_0\to |X|$.
%
%for it is the image of 
%the open set $|U_j^\al|\subset |B|$, and $|F|:|B|\to |X|$ is open
%since it is  finite to one and induced by the local diffeomorphism
%$F:B_0\to X_0$.
Hence $|F|_{\Hsm}(|U_j^\al|_{\Hsm}) = |U|$ since $|U|$ is connected. This proves 
the first statement.  Since there are only finitely many  pairs $(\al,j)$ 
and $\Xx$ is sse, one can now restrict $U$ further so that the second 
statement holds.
\end{proof}

\begin{lemma}\labell{le:comm} {\rm (i)}
If $F:\Aa\to \Xx$ and $G:\Cc\to \Xx$ are resolutions of the ep groupoid
$\Xx$ then
the (weak) fiber product $\Zz: = \Aa\times_\Xx\Cc$ of $\Aa$ and $\Cc$ over
$\Xx$ may be given the structure of a wnb groupoid, and the induced functors
$\Zz\to \Bb, \Zz\to \Aa$ are layered coverings.
\SSS

\NI{\rm (ii)}  If $F:\Bb\to \Xx$ is a resolution and $G:\Xx'\to \Xx$
is an equivalence then $ \Bb': = \Bb\times _\Xx\Xx'$ is a wnb groupoid. Moreover the induced functor $F':\Bb'\to\Xx'$
is a resolution, while $G':\Bb'\to \Bb$ is an equivalence.
\end{lemma}
\begin{proof} The proof of (i)  is very similar to that of 
Proposition~\ref{prop:commen} (ii), using
 Lemma~\ref{le:brstr1} instead of Lemma~\ref{le:brstr}.
Its details will be left to the reader. 

To see that $\Bb'$ in (ii) is a wnb groupoid one need only choose 
the set $U$ of Lemma~\ref{le:brstr1} so small that it is $\Xx$-diffeomorphic to a subset of $G(X'_0)$.  Then it can be lifted into $\Xx'$, and the rest of the proof is clear.
  \end{proof}

Note that the diagram
$$
\begin{array}{ccc}\Bb': =  \Bb\times _\Xx\Xx'&\stackrel{F'}\longrightarrow& \Xx'\vspace{.05in}\\
G'\downarrow\;\;\;\;&& G\downarrow\;\;\;\;\\
\Bb&\stackrel{F}\longrightarrow& \Xx,
\end{array}
$$
considered in Lemma~\ref{le:comm}(ii) does {\it not} in general commute strictly, but only \lq\lq up to homotopy".  In particular,
the map $F': B_0'\to X_0'$ is always 
surjective, while $F:B_0\to X_0$ may not be; cf. 
Remark~\ref{rmk:multi}(i).

The
 following result restates the main assertion of Theorem~\ref{thm}.   
 
\begin{prop}\labell{prop:resol}  Every ep groupoid $\Xx$ has a resolution that is
is unique up to commensurability and hence up to cobordism.
 \end{prop}

The  uniqueness statement follows immediately from Lemma~\ref{le:comm}
and Proposition~\ref{prop:commen}.
We give two proofs of the existence statement in \S\ref{sec:resol}.  The first 
 gives considerable control over the  branching structure of 
the resolution, while the second, which constructs the resolution as a multisection of a bundle $\Ee\to \Xx$, is perhaps more direct. However, it gives a resolution as defined above 
only when 
%%D$\Ee$ is nonsingular.
$\Xx$ acts effectively on $\Ee$.

\begin{rmk}\labell{rmk:nonred}\rm  One might argue that the above definition of resolution
 is not the most appropriate for  groupoids that are not effective; 
 cf. the example in Remark~\ref{rmk:resol} (ii).  We are mostly interested in resolutions because they give easily understood representatives for the fundamental class of an orbifold.  But we have defined the fundamental class of $\Xx$ so that it is
 the same as that of $\Xx_{\it eff}$; see \S\ref{ss:fnorb}.  Therefore we could define a resolution of $\Xx$ simply to be a resolution  $F:(\Bb,\la,\La)\to \Xx_{\it eff}$ of $\Xx_{\it eff}$.
 The information about the trivially acting part $K_x$ of the stabilizer groups $G_x$ of $\Xx$ would then be recorded  in the (strict) pullback\footnote
 {
 One could also use the weak pullback $\Bb\times_{\Xx_{\it eff}}\Xx$,
 but this is somewhat larger.
 }
  $(\Bb',\la',\La')$ of the fibration $\Xx\to \Xx_{\it eff}$ by $F$.  Here $\Bb'$ has the same objects, orbit space and weighting function as $\Bb$. But the morphisms in $\Bb'$ from $x$ to $y$ equal the set $\Mor_{\Xx}(F(x), F(y))$ if there is a morphism in $\Bb$ from $x$ to $y$, and equal the empty set otherwise. If  composition is defined by pull back from $\Xx$, one readily checks that $\Bb'$ is a groupoid.
 Since the morphisms in $\Bb'$ act trivially on the objects, one can easily extend the definition of wnb groupoid  to include this case.  
Hence
  $(\Bb',\la',\La')$ may also be considered as a kind of resolution of $\Xx$.  
  \end{rmk}

%%%%%%%%%%%%%%%%%%%%%%%%%%%%%%%%%%%%%%%%%%%%%%%%%%%%%%%%%%
\subsection{The fundamental class}
%%%%%%%%%%%%%%%%%%%%%%%%%%%%%%%%%%%%%%%%%%%%%%%%%%%%%%%%%%

We now show  that each compact  branched
manifold $(\un Z,\La_Z)$ of dimension $d$ carries a fundamental class $[Z]\in
H_d(Z,\R)$ which is compatible with resolutions; that is, if $\un\phi: (\un
Z,\La_Z)\to \un Y$ is any resolution, then $\phi_*([Z]) = [Y]$, where
$[Y]$ is the fundamental class of the orbifold $\un Y$ discussed in
\S\ref{ss:fnorb}.  We shall define $[Z]$ as a singular cycle using triangulation, but also show in Proposition~\ref{prop:fclass}
that in \lq\lq nice" (i.e. tame) 
cases, one can integrate over $[Z]$.   In order to
construct a suitable integration theory we shall need to consider
smooth partitions of unity.

\begin{defn}\labell{def:smooth}
Let $(\un Z,\La_Z)$ be a branched manifold and $M$ a smooth
manifold.  A map $g:Z\to M$ is  smooth iff for one
(and hence any) branched 
%%Dgroupoid 
manifold structure $(\Bb,\La,f)$ on $\un Z$
the composite
$$
g_0: B_0\stackrel{\pi_{\Hsm}}\to |B|_{\Hsm}\stackrel f\to Z\stackrel g\to M
$$
 is  smooth.
This is equivalent to saying that $g$ is induced
by a smooth functor $\Bb\to \Mm$, where $\Mm$ is the category with
objects $M$ and only identity morphisms. 

A {\bf smooth partition of unity  subordinate to the 
covering} $\Nn = \{N_k\}_{k\in A},$ of $Z$ is a family $\be_k, k\in
A,$ of smooth functions $Z\to \R$ such that \SSS

\NI
{\rm (i)} $\supp\, \be_k\subset N_k$ for all $k$,\SSS

\NI
{\rm (ii)} for each $z\in Z$ only finitely many of the $\be_k(z)$ are
nonzero, and \SSS

\NI {\rm (iii)} $\sum_k \be_k(z) = 1$.
\end{defn}

The above smoothness condition is quite strong.  For example, if
$\phi:D_1\to R_2$ in the example of 
Remark~\ref{rmk:brorb}(ii) does not
extend to a smooth function near $\p D_1$, there is no smooth function
$g:|B_\phi|_{\Hsm}\to \R^{2}$ that is injective over the image of
$R_1$.  We shall deal with these problems by introducing the notion of
tameness.

%%D added
\begin{defn} Let $\Om$ be a precompact open subset of $\R^d$. 
 It is said to have {\bf piecewise smooth boundary} if for every 
 $x\in \p \Om$ there is a neighborhood $U_x$ of $x$ and smooth functions $f_1,\dots,f_k:U_x\to \R$ such that\SSS
 
 \NI
 {\rm (i)}  $U_x\cap \ov\Om = \{y\in U_x : 
 f_i(y)\le 0,$ for all $ i=1,\dots,k\}$;\SSS
 
 \NI
 {\rm (ii)} $x\in \p \Om \Longleftrightarrow f_i(x) = 0$ for 
 at least one $i$;\SSS
 
 \NI
 {\rm (iii)} Let $x\in \p\Om$ and write $I_x = \{i: f_i(x) = 0\}$.
 Then the set of vectors $df_i(x),i\in I_x,$ is linearly independent.
 \SSS
 
 \NI
 Moreover, $\Om$ is said to have 
 {\bf piecewise smooth boundary} over an open set $N\subset\R^d$ if the above conditions hold for every $x\in N\cap \p\Om$.
\end{defn}

Note the following facts about domains $\Om$
 with piecewise smooth boundary.
\MS

\NI$\bullet$  Any smooth function defined on the closure $\ov\Om$ 
has a smooth extension to $\R^d$;\SSS

\NI$\bullet$  The boundary $\p \Om$ has zero measure;\SSS

\NI$\bullet$ The intersection $\cap_{j\in J}\Om_j$ of a generic finite collection $\Om_j, j\in J,$ of such domains also has 
piecewise smooth boundary. More precisely, suppose that each $\ov\Om_j$ is contained in the interior of the ball $B_R$ of radius $R$.  Then the set of diffeomorphisms $\phi_j, j\in J$ of $B_R$ such that $\cap_{j\in J}\phi_j(\Om_j)$ has piecewise smooth boundary
has second category in the group $\prod_{j\in J}\Diff(B_R)$.
Sets $\Om_j, j\in J,$ with this property will be
 said to be in {\bf general position} or to {\bf intersect transversally}.
\MS

%%changed

%\begin{defn}\labell{def:tame}
%  A $d$-dimensional wnb groupoid $(\Bb, 
%  \La)$ is said to be {\bf tame}
%iff the following conditions are satisfied.\SSS

%\NI
%{\rm (i)}  One can choose the local branch structures so that each local branch is a union of components of $B_0$.\SSS

\begin{defn}\labell{def:tame}
  A $d$-dimensional wnb groupoid $(\Bb, 
  \La)$ is said to be {\bf tame}
iff if it has local branching structures 
$(N_{p_k}, U_i^k,m_i^k)_{k\in K}$ such that 
the following conditions are satisfied.\SSS

\NI
{\rm (i)}  $|B|_{\Hsm} = \cup_{k\in K}N_{p_k} $.\SSS

\NI
{\rm (ii)} For each local branch $U_i^k$
there
is an injective smooth map  $\phi_i^k:
U_i^k\to \R^d$ onto the interior of a 
%%Dclosed 
compact domain
$\ov\Om\,\!^k_i$ in $\R^d$ with piecewise smooth boundary 
such that the composite $\pi_{\Hsm}\circ (\phi_i^k)^{-1}$ extends to an injection
$\rho_i^k:  \ov\Om\,\!^k_i\to 
|B|_{\Hsm}$.  
\SSS

\NI
{\rm (iii)}  For each pair $U_i^k, U_j^\ell$
the set $\ov\Om\,\!^k_i\cap (\rho_i^k)^{-1}(\rho_j^\ell(\ov\Om\,\!^\ell_j))$ has
piecewise smooth boundary and
 the transition map 
$
(\rho_j^\ell)^{-1}\circ \rho_i^k$ extends smoothly to a local 
diffeomorphism defined on a neighborhood of this set.
\SSS

If $\Bb$ is a groupoid in the category of manifolds 
with boundary we replace $\R^d$ in the above by the 
half space ${\mathbb H}^d = \{x\in \R^d: x_1\ge 0\},$ and require 
that all sets meet $\p {\mathbb H}^d$ transversally.
\end{defn}

The main point of this definition is to ensure that
that the branch locus is piecewise smooth
and hence has zero measure.

%Here condition (i) is not essential, and was included to make it 
%easier to state the subsequent conditions.  It can always be 
%achieved by refining $\Bb$: cf. Remark~\ref{rmk:brorb0}(ii). 
%The other conditions are more meaningful since they imply 
%that the branch locus is piecewise smooth
%and hence has zero measure.
%
% As we saw in Lemma~\ref{le:branch} the latter property is a 
% consequence of condition (b) in Definition~\ref{def:brorb}. 
% However the proof of Lemma~\ref{le:tame} below shows that even
% if we do not assume (b) at the outset we can achieve it for a 
% commensurate groupoid by suitable pruning of the morphisms.

 \begin{lemma}\labell{le:tame0} Every weighted refinement of a tame wnb
 groupoid may be further refined to be tame.
 \end{lemma}
  \begin{proof}  Suppose that $(\Bb,\La)$ is tame.
 Since $|B|_{\Hsm}$ is unchanged under refinement, 
 the branching locus of any refinement  $F:\Bb'\to\Bb$  is piecewise smooth. 
 If we use the pullback local branching structures on $\Bb'$, then the
 taming condition (ii) 
 may not be preserved because it 
 concerns the objects rather than the morphisms of $\Bb'$. 
% While taming conditions (i) and (ii) 
% may not be preserved by an arbitrary refinement, because 
% these conditions concern the objects rather than morphisms of $\Bb$ any refinement
Nevertheless it is easy to see that there is
 a further refinement $\Bb''\to\Bb'$ for which it holds.
 Then (iii) also holds.
 \end{proof} 
 
In view of the above lemma, we shall say that the {\bf branched manifold $(\un
Z, \un \La_Z)$ is tame} if its structure
may be represented by a tame wnb
groupoid.

In the following lemma we write
$V\sqsubset U$ to mean that the closure $cl(V)$ of $V$ is contained in
$U$.

%%D revised
\begin{lemma}\labell{le:tame}  Every 
 wnb groupoid $(\Bb,\La_B)$ has a layered covering that is tame.
\end{lemma}

\begin{proof} Choose a locally finite cover of $|B|_{\Hsm}$ by sets
$N_{k}, k\in \N$, as in Definition~\ref{def:brorb}.  For each $k$ 
let $U_i^k, i\in A_k,$ be the corresponding local branches. 
Without loss of generality we may suppose that each
local branch $U_i^k$ is identified with a subset of $\R^d$.
Since $|B|_{\Hsm}$ is normal we may choose open subsets $W_k', N_k', W_k$ of
$N_{k}$ such that the $N_k'$ cover $|B|_{\Hsm}$ and
$$
W_k' \,\sqsubset\, N_k'\,\sqsubset\, W_k \,\sqsubset\, N_{k},\qquad\mbox{ for all } k.
$$
Since $\pi_{\Hsm}: U_i^k\to |U_i^k|_{\Hsm}$ is a homeomorphism for each $i$, it
follows that
$$
U_i^k\cap\pi_{\Hsm}^{-1}(W_k') \,\sqsubset\, U_i^k\cap
%%D(\pi_{\Hsm})^{-1}(N_k')\,\sqsubset\,
\pi_{\Hsm}^{-1}(N_k')\,\sqsubset\,
U_i^k\cap \pi_{\Hsm}^{-1}(W_k)
\,\sqsubset\,U_i^k \sqsubset \R^d.
$$

%In particular, we may also choose the $N_k'$ so that the intersection of
%$\pi_{\Hsm}^{-1}(N_k')$ with each local branch $U_i^k$ has piecewise
%smooth boundary.
 
We now
%tame the groupoid $\Bb$ over the sets 
%$N_k'$ for $k=1,2,\dots$ in turn,
construct a sequence of
commensurate wnb groupoids $\Bb^0\supset\Bb^1\supset\Bb^2\supset\dots
$ and sets $Z_1\subset Z_2\subset \dots \subset |B|_{\Hsm}$ 
such that, for all $k\ge 1$,  $\Bb^k$ is tame over 
 $Z_k \supset \cup_{j=1}^k W_j'$ and also $\Bb^{k+1}
= \Bb^{k}$ over $Z_{k}$.  Then  $\Bb^\infty: = \lim \Bb^k$ will be the
desired tame groupoid.  (We weight all these groupoids by pulling back
$ \La_B$.  Further to say that $\Bb^k$ is tame over $Z_k$ means 
that
its full subcategory with objects in $\pi_{\Hsm}^{-1}(Z_k)$ is tame.)

We define $\Bb^0$ to be the full subcategory of $\Bb$ with objects 
$$
B^0_0: = \underset{k}\sqcup \Bigl(\underset{i\in A_k}\sqcup\;U_i^k\Bigr).
%\;U_i^k\cap\pi_{\Hsm}^{-1}(N_k')\Bigr).
$$
Because the sets $N_p$ cover $|B|_{\Hsm}$,  the inclusion 
$\Bb^0\to \Bb$ is 
a weighted equivalence. 
%Note that $\Bb^0$ satisfies the taming condition (i).% and (ii).
 
The next step is to
correct the morphisms over $N_1'$. 
%To do this,
%consider the full subcategory $\Ss$ of $\Bb^0$ with objects the sets $U_i^1, i\in A_1 = \{1,\dots,N_1\}$.
Let 
$$
\Ii_k: = \{I\subset A_k: |U_I^k|_{\Hsm}: = \bigcap _{i\in I} |U_i^k|_{\Hsm}\ne \emptyset\}.
$$
Order the sets in $\Ii_1$ by inclusion.  Starting with 
a maximal $I\in \Ii_1$ and then working down, 
choose sets $|V_I^1|_{\Hsm}\subseteq |U_I^1|_{\Hsm}$ satisfying 
the following conditions for all $I,J$:
\MS

\NI$\bullet$ $I\supset J  \Longrightarrow |V_I^1|_{\Hsm} \subseteq |V_J^1|_{\Hsm}$\,;\SSS

\NI$\bullet$ $ |V_I^1|_{\Hsm}\cap (N_1\setminus W_1) = |U_I^1|_{\Hsm}$\,;\SSS

\NI$\bullet$
for one (and hence every) $i\in I$ the subset $\pi_{\Hsm}^{-1}(|V_I^1|_{\Hsm})$ of $U_i^1\subset \R^d$ has piecewise smooth boundary over $\pi_{\Hsm}^{-1}(N_1')$.\SSS

\NI$\bullet$  $|V_i^1|_{\Hsm} = |U_i^1|_{\Hsm}$ for all $i$.
\MS

For all $k>1$ and $I\subset A_k$ set $|V_I^k|_{\Hsm}: = |U_I^k|_{\Hsm}$. Then
define $\Bb^1$ to be the
 subcategory of $\Bb^0$ with the same objects, labelled for convenience as $V_i^k$ instead of $U_i^k$,  and with morphisms 
 determined by the identities
  $$
 \bigcap_{i\in I} |V_i^k|_{\Hsm} = |V_I^k|_{\Hsm},\quad\mbox{for all } I\subset A_k, k\ge 1.
 $$
Then $\Bb^1$ is a wnb groupoid commensurate to $\Bb^0$. (Note that the inclusion $|B^1|_{\Hsm}\to |B^0|_{\Hsm}=|B|_{\Hsm}$ is proper because we 
did not change the morphisms near the boundary $\p N_1$.)  Now choose $Z_1$ so that $W_1'\sqsubset Z_1\sqsubset N_1'$ and so that for each $i\in A_1$ its pullback $\pi_{\Hsm}^{-1}(Z_1)\cap U_i^1$ has piecewise smooth boundary 
and is transverse to the sets $\pi_{\Hsm}^{-1}(V_I^1)\cap U_i^1$ for all $I\in \Ii_1$.  
Then,   $\Bb^1$ is tame over $Z_1$.  

 We next repeat this cleaning up process over $N_{2}$, making no changes
 to the morphisms lying over a neighborhood of $Z_1\cup (N_{2}\setminus
 W_2)$ and taming the morphisms over $N_2'$. We then choose a suitable set $Z_2\supset
 Z_1\cup W_2'$ to obtain
a groupoid
 $\Bb^2$ that is tame over $Z_2$.  Continuing this way
 we construct the $\Bb^k$ and hence $\Bb^\infty$.\end{proof}

\begin{lemma}\labell{le:part}  Let  $(\Bb,\La)$ be a tame
    wnb groupoid and let $\Nn$ be any open cover of $|B|_{\Hsm}$.  Then
    $|B|_{\Hsm}$ has a smooth partition of unity subordinate to $\Nn$.
%\SSS

%\NI
%{\rm (ii)}  If the orbit space $|B|_{\Hsm}$ of a wnb groupoid $(\Bb,\la,\La)$ admits smooth partitions of unity 
% then there is a branched partition of $\La$.
\end{lemma}
\begin{proof}  By Remark~\ref{rmk:brorb}(iv), we may suppose that 
    $\Nn$ consists of sets 
 of the form $N_p$, where $(N_p, U_i^p, m_i)$ are local branching
 structures.  Pick out a countable subset $A$ such that the sets $N_p,
 p\in A,$ form a locally finite covering of $|B|_{\Hsm}$.  
 Since
 $|B|_{\Hsm}$ is normal by Proposition~\ref{prop:approx}, there are open
 sets $N'_p\sqsubset N''_p\sqsubset N_p$ such that $\{N'_p\}_{p\in A}$
 is an open cover of $|B|_{\Hsm}$.  For each $p\in A$ we shall
 construct a smooth function $\la_p:N_{p}\to [0,1]$ that equals $1$ on
 $N_p'$ and has support in $N_p''$.  Then $\la: = \sum_{p\in A} \la_p:
 |B|_{\Hsm}\to \R$ is everywhere positive and smooth.  Hence the
 functions $\be_p: = \la_p/\la$ form the required  partition of unity.
 
For each $p$, we construct $\la_p: N_p\to \R$ inductively over its subsets
$Q_m^p: = \cup_{i=1}^m |U_i^p|_{\Hsm}$, where $U_1^p,\dots,U_k^p$ are
the local branches at $p$.  To begin, choose a smooth function
$f_1^p: U_1^p\to [0,1]$ that equals $1$ on $\pi_{\Hsm}^{-1}(cl
(N'_p))$ and has support in $\pi_{\Hsm}^{-1}(N_p'').$ This exists
because the map $\pi_{\Hsm}: U_1^p\to N_p$ is proper.  Since $U_1^p$
injects into $N_p$, we may define $\la_p(q)$ for $q\in Q_1^p$ by
 $$
 \la_p(q): = f_1^p(x),\quad\mbox{where }\; \pi_{\Hsm}(x) = q.
 $$ 
 By the tameness hypothesis, the pullback by $\pi_{\Hsm}$ of $\la_p$ to
 $U_2^p$ (which is defined over $U_2^p\cap
 \pi_{\Hsm}^{-1}(|U_1^p|_{\Hsm})$) may be extended over $U_2^p$ to a
 smooth function $f_2^p$ that equals $1$ on $\pi_{\Hsm}^{-1}(cl
 (N'_p))$ and has support in $\pi_{\Hsm}^{-1}(N_p'').$ Now extend
 $\la_p$ over $Q_2^p$ by setting it equal to the pushdown of $f_2^p$ over
 $ \pi_{\Hsm}(U_2^p)$.  Continuing in this way, one extends $\la_p$ to
 a function on the whole of $N_p$ that equals $1$ on
 $\pi_{\Hsm}^{-1}(cl(N'_p))$ and has support in
 $\pi_{\Hsm}^{-1}(N_p'').$
 
 It remains to check that $\la_p:|B|_{\Hsm}\to \R$ is smooth.  Its pull
 back to any local branch $U_i^p$ is smooth by construction.  Consider
 any other point $x\in B_0$ such that $\pi_{\Hsm}(x) \in N_p$.  Then
 there is some point $y\in U_i^p$ such that $|x|=|y|$ by the covering
 property of the local branches.  Hence there is a local diffeomorphism
 $\phi_{yx}$ of a neighborhood of $x$ to a neighborhood of $y$ and so
 $\la_p\circ \pi_{\Hsm}$ is smooth near $x$ because it is smooth near
 $y$.
\end{proof}

%Consider a compact branched orbifold $(\un Y, \La_Y)$.

%
%Although one could adapt the simplicial techniques 
%of Moerdijk--Pronk~\cite{MPr} to construct a simplicial fundamental cycle on 
% $(\un Y, \La_Y)$, we shall rather work in the smooth category.  

Let $(\Bb,\La,f)$ be a compact tame $d$-dimensional branched manifold
structure on $(\un Z, \La_Z)$.  Choose a smooth partition of unity
$\{\be_p\}$ on $|B|_{\Hsm}$ that is subordinate to a covering by sets $N_p$ that
support local branching structures $(N_p, U_i^p, m_i)$.  If $g: Z\to
M$ is any smooth map into a manifold, and $\mu$ is a closed $d$-form
on $M$ (where $d = \dim Z$), we define
\begin{equation}\label{eq:int}
\int_{Z} g^*\mu: = \sum_{p,i}m_i\int_{U_i^p}
(\pi_{\Hsm})^* \be_p\, (g\circ f\circ
\pi_{\Hsm})^*\mu.
\end{equation}
As we explain in more detail below, the reason why this is well
defined and independent of choices is that, because $\Bb$ is tame, its
branching locus is a finite union of piecewise smooth manifolds of
dimension $d-1$ and so has zero measure.\footnote
{  
Achieving an analog of this is also a crucial step
in the work of K. Cieliebak, I. Mundet i Riera and D. Salamon, though
they use rather different methods to justify it: cf.~\cite[Lemma~9.10]{CRS}.}

\begin{lemma} {\rm (i)} The number $\int_{Z} g^*\mu$ defined above is
independent of the choice of partition of unity. \SSS
    
    \NI
    {\rm (ii)} If $g$ is bordant to $g': Z'\to M$ by a
bordism through a tame wnb groupoid, then $\int_{Z}
g^*\mu=\int_{Z'} (g')^*\mu$.  In particular, it is
independent of the choice of $(\Bb,\La)$.
\end{lemma}
  \begin{proof}
      First suppose that the partitions of unity $\be_p, \be_p'$ are
  subordinate to the same covering and consider the product groupoid
  $\Bb\times \Ii$ where $\Ii$ has objects $I_0:=[0,1]$ and only identity
  morphisms and $\La_I \equiv 1$.  There is a partition of unity
  $\{\be''_p\}$ on $\Bb\times \Ii$ that restricts on the boundary to the
  two given partitions of unity and is subordinate to the cover
  $N_p\times |I|$.  Hence by Stokes' theorem, it suffices to show that
 $$
 \sum_{p,i}\int_{U_i^p} m_i (\pi_{\Hsm})^*(\be_p - \be_p')\, (\pi_{\Hsm}\circ pr)^*\mu_B
 =
  \sum_{p,i}m_i\int_{U_i^p\times I_0} (\pi_{\Hsm})^*(d\be_p'')\, 
  (\pi_{\Hsm})^*\mu_B = 0,
  $$
where $\mu_B: =  (g\circ f)^*\mu$ and $\pr:B_0\times I_0\to B_0$ is the projection.  Let $V\subset (|B|_{\Hsm}\setminus |Br|_{\Hsm})$ be a component of the complement of
the branching locus, and let $\La(V)$ be the constant value of $\La$ on
$V$ (cf.  Proposition~\ref{prop:approx}).  Because the branching locus in each
$U_i^p$ has zero measure it suffices to show that the sum of the
integrals over $\bigl(U_i^p\cap \pi_{\Hsm}^{-1}(V)\bigr)\times I_0$
vanishes for each $V$.  But $V\cap N_p$ is diffeomorphic to $U_i^p\cap
\pi_{\Hsm}^{-1}(V)$ for every $i$ for which the intersection is
nonempty.  Hence, because $\supp(\be_{p})\subset N_p$, it makes sense
to integrate over $V$ and $V\times I_0$, and we find that
\begin{eqnarray*}
\sum_{p,i}m_i\int_{(U_i^p\cap \pi_{\Hsm}^{-1}(V))\times I_0}
(\pi_{\Hsm})^*(d\be_p'')\; (\pi_{\Hsm}\circ pr)^*\mu_B & = &
\sum_{p}\La(V)\int_{V\times I_0}
d\be_p''\; pr^*\mu_B
%\\
%& = & \La(V)\int_{V\times I_0}
%\sum_p d\be_p''\; pr^*\mu_B 
 \end{eqnarray*}
 which vanishes because $\sum_p \be_p'' = 1$.
 
 The proof of Proposition~\ref{prop:commen}
 shows that any two covers that support local branching
 structures have a common refinement that supports a local branching
 structure.  Hence to prove (i) it suffices to show that if $\{\be_p\}$
 is subordinate to $\{N_p\}$ and the cover $\{N_q'\}$ refines $\{N_p\}$
 then there is some partition of unity subordinate to $\{N_q'\}$ for
 which the two integrals are the same.  The previous paragraph shows
 that the first integral may be written as
  $$
  \sum_V \sum_p \int_V \La(V) \be_p\; \mu_V,
  $$
  and so this statement holds by the standard arguments valid for manifolds.
  
  To prove (ii), note that if $(\Ww,\La_W)$ is a tame cobordism from
  $(\Bb,\La)$ to $(\Bb',\La')$ then every partition of unity on its
  boundary extends to a partition of unity on the whole groupoid. 
  Moreover we can construct the extension to be subordinate to any
  covering that extends those on the boundary.  The rest of the details
  are straightforward, and are left to the reader.
\end{proof}

\begin{prop}\labell{prop:fclass} Let $(\un Z, \La_Z)$ be 
    a compact $d$-dimensional branched manifold with boundary.  
    Then the singular homology group $H_d(Z, \p Z;\R)$ contains
    an element $[Z]$ called the fundamental class with the 
    following properties: \SSS

\NI
{\rm (i)}  If the weights of all the branches of $Z$ are rational then
$[Z]\in H_d(Z;\Q)$.\SSS

\NI
{\rm (ii)}  If $\un \phi: (\un Z', \La_{Z'})\to (\un Z, \La_Z)$ is any 
layered covering, $\phi_*([Z']) = [Z].$\SSS

\NI
{\rm (iii)}  The image of $[Z]$ under the boundary map
$H_d(Z, \p Z;\R)\to H_{d-1}(\p Z;\R)$ is $[\p Z]$.\SSS

\NI
{\rm (iv)} If $\un\phi:(\un Z,\La_Z)\to \un Y$ is any 
resolution of the orbifold $\un Y$, then
$\phi_*([Z]) = [Y]\in H_d(Y;\R)$.\SSS

\NI
{\rm (v)} Suppose further that $\un Z$ is tame.  Then for any smooth map
$g:Z\to M$ of $Z$ into a smooth manifold $M$ and any closed $d$-form
$\mu$ on $M$,
$$
[\mu](g_*([Z])) = \int_{Z} g^*(\mu).
$$
\end{prop}
\begin{proof} By Lemma~\ref{le:tame} we may assume that $(\un Z, \La_Z)$
is commensurate with a branched manifold  $(\un B, \La_B)$ 
with structure given
 by the tame nonsingular wnb  groupoid $(\Bb,\La_B)$.  The
tameness condition implies that $|B|_{\Hsm}$ may be triangulated
in such a way that both the branching locus $|Br|_{\Hsm}$ and the boundary of
$|B|_{\Hsm}$ are contained in the $(d-1)$-skeleton.  More precisely,
we can arrange that any (open) $(d-1)$-simplex that intersects this
$|Br|_{\Hsm}\cup \p(|B|_{\Hsm})$ is entirely contained in this set, and
that no open 
$d$-simplex meets it.  By first triangulating the boundary we may 
assume that similar statements (with $d$ replaced by $d-1$) hold for
its branching locus. 
Then, $\La_B$ is constant on each open $d$-simplex $\si$ in the
triangulation $\Tt$. To simplify the proof below we will refine $\Tt$ 
until each of its
$(d-1)$-simplices $\rho$ lies in the support $N_p$ of a local branching structure such that $N_p$ contains all the $d$-simplices that meet $\rho$.

Suppose that $\Bb$ has no boundary.  Then we claim that the singular
chain defined on $|B|_{\Hsm}$ by
$$
c(|B|_{\Hsm}): = \sum_{\si\in \Tt} \La_B(\si)\, [\si]
$$ 
is a cycle. 
To see this, consider an open $(d-1)$-simplex $\rho$.  If
$\rho$ lies in $|B|_{\Hsm}\setminus |Br|_{\Hsm}$ then it is in the boundary of
precisely two oppositely oriented $d$-simplices with the same weights. 
Hence it has zero coefficient in $\p c(|B|_{\Hsm})$.  
Suppose now that $\rho$
lies in the branching locus.   Choose  a 
local branching structure
$(N_p,U_i^p, m_i)$ such that $\rho$ and all the $d$-simplices that meet it are contained in $N_p$.  Each simplex $\si$ whose boundary contains $\rho$
lies in a component of $|B|_{\Hsm}\setminus |Br|_{\Hsm}$ and so, for each $i$,
$\si\cap
|U_i^p|_{\Hsm}$ is either empty or is the whole of $\si$.  Moreover,
$$
\La_B(\si) =\sum_{i:\si\cap |U_i^p|_{\Hsm}\ne \emptyset} m_i.
$$
Hence  $\La_B(\si) [\si]$ is the pushforward of the 
chain 
$$
\sum_{i:\si\cap |U_i^p|_{\Hsm}\ne \emptyset} m_i \; [(\pi_{\Hsm})^{-1}
(\si)\cap U_i^p]\quad\mbox{on } B_0.
$$
Because the simplices in $U_i^{p}$ cancel each other out in pairs in the
usual way, $\rho$ again makes no contribution
 to $\p c(|B|_{\Hsm})$.
We now define $[|B|_{\Hsm}]$ to be the singular  homology class
represented by $c(|B|_{\Hsm})$, and $[Z]$ to be 
its pushforward by $|B|_{\Hsm}\to Z$.

Note that if $\un Z$ and hence $\Bb$ has boundary, then (iii) holds 
for the cycles $[|B|_{\Hsm}]$ and
$[Z]$ by our choice of triangulation.  Any two triangulations of 
$|B|_{\Hsm}$ can be considered as a triangulation of 
$|B|_{\Hsm}\times \{0,1\}$ and then extended over  
$|B|_{\Hsm}\times [0,1]$.  Applying (iii) to $[|W|_{\Hsm}]$
where $\Ww = \Bb\times \Ii$, one easily sees that 
$[Z]$ is independent of the choice of triangulation.
A similar argument shows that 
 $[Z]$ is independent of the choice of representing groupoid $(\Bb,\La_B)$, since any two such groupoids are cobordant
 by Proposition~\ref{prop:commen}(iii).  
 Therefore the singular homology
 class $[Z]\in H_d(Z)$ is independent of all choices. It satisfies (i) by definition.

The other statements follow by standard arguments.  In
 particular (ii) holds by a cobordism argument and
 (v) holds because we can assume that both the triangulation
 and the partition of unity are subordinate to the same covering $\{N_p\}$.
 Hence each $d$-dimensional simplex in $\Tt$ lies in some component $V$ of $|B|_{\Hsm}\setminus |Br|_{\Hsm}$
 and so we can reduce this to the usual statement for manifolds.
\end{proof}

\begin{rmk}\rm   Above we explained the class $g_*([Z])$ 
    for a smooth map $g:Z\to M$ in terms of integration.  
    However, it can also be understood in terms of intersection theory.  
In fact, because $[Z]$ is represented by the singular cycle $f:
c(|B|_{\Hsm})\to Z$, the number $[\mu](g_*([Z]))$ may be calculated by
counting the signed and weighted intersection points of the cycle
$g\circ f: c(|B|_{\Hsm})\to M$ with any singular cycle in $M$
representing the Poincar\'e dual to $[\mu]$; cf. 
Example~\ref{ex:brorb}(ii).  Detailed proofs of very similar
statements may be found in~\cite{CRS}.
\end{rmk}

%%%%%%%%%%%%%%%%%%%%%%%%%%%%%%%%%%%%%%%%%%%%%%%%%%%%%%%%%%
\section{Resolutions}\labell{sec:resol}
%%%%%%%%%%%%%%%%%%%%%%%%%%%%%%%%%%%%%%%%%%%%%%%%%%%%%%%%%%

Our first aim in this section is to show that every  orbifold 
$\un X$
has a resolution.  We then discuss the relation between resolutions and the (multi)sections of orbibundles. Although we shall assume that $\un X$ is finite dimensional, many arguments apply to orbifolds in any category.

%%%%%%%%%%%%%%%%%%%%%%%%%%%%%%%%%%%%%%%%%%%%%%%%%%%%%%%%%%
\subsection{Construction of the resolution}\labell{ss:resol}
%%%%%%%%%%%%%%%%%%%%%%%%%%%%%%%%%%%%%%%%%%%%%%%%%%%%%%%%%%

Let $\un Y$ be  a (possibly not effective) orbifold.
Choose a good atlas  $(U_i, G_i,\pi_i), i\in A,$  for $\un Y$ 
and use it to construct 
 an orbifold structure $\Xx$ on $\un Y$ with objects $\sqcup U_i$;
cf. the discussion after Definition~\ref{def:atl}.  For each finite subset $I\subset A$, we denote
$$
|U_I|: = \bigcap_{i\in I} |U_i|\subset |X|,\qquad G_I: = \prod_{i\in I} G_i.
$$
We shall not assume that the sets $|U_I|$ are connected, 
although it would slightly simplify the subsequent argument to do so.
We shall identify the countable set $A$ with a subset of $\N$, and shall 
write $I$ as
$\{i_1,\dots,i_k\}$, where $i_1<i_2<\dots < i_k$. The 
length of $I$ is
$|I|: = k$.

We now define some sets $\HHat U_I$.
If $I= \{i\}$ we set $\HHat U_i: = U_i$.  When $|I|>1$
we define $\HHat U_I$ to be the set of composable tuples 
$(\de_{k-1},\dots,\de_1)$ of morphisms, where 
$$
s(\de_{j}) \in U_{i_{j}}, \mbox{ for } 1\le j\le k-1,\;\;\,\mbox{ and }  t(\de_{k-1})\in U_{i_{k}}.
$$
Since $s$ and $t$ are local diffeomorphisms, $\HHat U_I$ is a manifold.\footnote
{
It is important that $\HHat U_I$ consists of morphisms rather than being the corresponding fiber product of the $U_i$ over $|U_I|$, since that is not a manifold in general.  For example, 
suppose that $U_1=U_2 = \R$ and $G_1=G_2 = \Z/2\Z$ acting by $x\mapsto -x$ and that each $U_i$ maps to $|X| =  [0,\infty)$ by the obvious projection.
Then $\HHat U_{12}$ can be identified with the two disjoint lines $\{(x,x): x\in \R\}$  and
$\{(x,-x): x\in \R\}$.
 On the other hand, 
the fiber product $U_1\times _\pi U_2$ consists of two lines that intersect at $(0,0)$.  Thus this step reformulates 
Liu--Tian's concept of desingularization in the language of orbifolds: see \cite[\S4.2]{Mcv}.}
Further,  $\HHat U_I$ supports an  action of the group $G_I$
%: = \prod_{i\in I} G_i$ 
via:
$$
(\de_{k-1},\dots,\de_1)\cdot (\ga_k,\dots,\ga_1) = 
(\ga_k^{-1} \de_{k-1} \ga_{k-1},\dots, \ga_2^{-1} \de_1\ga_1).
$$
The
 action of $G_I$   
 is not in general free.  Indeed if $x=s(\de_1)$ the
stabilizer of $(\de_{k-1},\dots,\de_1)$ is a subgroup of $G_I$ isomorphic to $G_x$. 
%%D added
(For example, given $\de\in \HHat U_{12}$ and $g_1\in G_{s(\de)}$ there is a unique $g_2\in G_{t(\de)}$ such that $\de = g_2^{-1}\de g_1$.)
 On the other hand,
for each $\ell\in I$
the group 
\begin{equation}\labell{eq:G'}
G'_{\ell I}: = \prod_{i\in I\setminus \ell}G_i
\end{equation}
 does act freely on $\HHat U_I$.
The obvious projection $\pi_I: \HHat U_I\to |U_I|\subset |X|$ identifies
 the  quotient $\HHat U_I/G_I$ with $|U_I|$.  
 Note that each component of  $\HHat U_I$ surjects onto a component of $|U_I|$.
% We shall denote the components of $|U_I|$ by $|U_I|^\al$ 
% and shall write $\HHat U_I^\al$ for the inverse image of $|U_I|^\al$. 
%Further, 
%because we assumed that  $|U_I|$ is connected, each component of $\HHat U_I$ surjects onto $|U_I|$. 

The resolution is a groupoid with objects contained in the sets 
$\HHat U_I$ and morphisms given by certain projections $\HHat\pi_{JI}$
 that we now explain. 
 For each $i_j\in I$ there is a
projection $\HHat\pi_{i_jI}: \HHat U_I\to \HHat U_{i_j} = U_{i_j}$ defined as
\begin{eqnarray*}
\HHat\pi_{i_jI}(\de_{k-1},\dots,\de_1) &= &s(\de_j), \quad \mbox { if } j<k\\
&= & t(\de_{k-1}), \quad \mbox { if } j=k.
\end{eqnarray*}
If $k=|I|\ge 3$ and $J = I\setminus\{i_\ell\}$ we define a projection $\HHat\pi_{JI}:\HHat U_I\to \HHat U_J$ as follows:
\begin{eqnarray*}
\HHat\pi_{JI} (\de_{k-1},\dots,\de_1) & = & (\de_{k-1},\dots,\de_2)\quad \mbox{ if } \ell=1,\\
& = & (\de_{k-1},\dots, \de_\ell\de_{\ell-1},\dots,\de_1)\quad 
\mbox{ if } 1<\ell<k,\\
 & = & (\de_{k-2},\dots,\de_1)\quad \mbox{ if } \ell=k.
\end{eqnarray*}
This map is equivariant with respect to the actions of $G_I$ and $G_J$, and identifies the image $\HHat\pi_{JI}(\HHat U_I)$ as the quotient of $\HHat U_I$ by a free action of $G_{i_\ell}$.
If $J$ is an arbitrary subset of $I$ with $|J|>1$ we define $\HHat\pi_{JI}$ as a composite of these basic projections. 
If $J = \{j\}$ we define $\HHat\pi_{JI}: = \HHat\pi_{jI}$.
Clearly $\HHat\pi_{jJ}\circ\HHat\pi_{JI} =  \HHat\pi_{jI}$,
whenever $\{j\}\subset J\subset I$. 

% \begin{figure}[htbp]
%\centering
%%\includegraphics[width=12cm]{orb2i.jpgtex_t}
%\scalebox{.3}{\input{orb2i.jpgtex_t}}
%%\includegraphics[width=3in]{orb2fig.eps}
%\caption{The covering $\Vv$}
%\end{figure}

\begin{figure}[htbp] %  figure placement: here, top, bottom, or page
%xx   \centerline{\psfig{figure=orbfig3.jpg,width=3in}}
   \centering
   \includegraphics[width=3in]{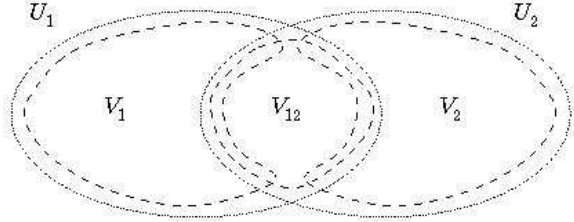} 
   \caption{An $\Aa$-compatible covering $\Vv$; the sets $V_I$ have dashed boundaries.}
   \label{fig3}
\end{figure}

\begin{defn}\labell{def:resol} Let $\Xx$ be an ep groupoid constructed from the good atlas  $\Aa$.  
A covering $\Vv=\{|V_I|\}_{I\subseteq A}$ of $|X|$ is called {\bf $\Aa$-compatible} if the following conditions hold:\SSS

\NI
{\rm (i)}  $|V_I|\sqsubset
|U_I|$  for all $I$,  \SSS

\NI
{\rm (ii)} $cl(|V_I|)\cap cl(|V_J|)\ne \emptyset $ iff $|V_I|\cap |V_J|\ne \emptyset$ iff  one of $I,J$ is contained in the other.\SSS

\NI
{\rm (iii)} for any $i\in I$, any two distinct components of $\pi_i^{-1}(|V_I|)\subset U_i$ have disjoint closures in $U_i$.\SSS

\NI
If $\Vv$ is $\Aa$-compatible, the {\bf $\Vv$-resolution} $\Xx_\Vv$ of $\Xx$ is a nonsingular  groupoid defined as follows.
Its set of objects is 
$$
(X_\Vv)_0: = \sqcup_I \HHat V_I,
$$
where  $ \HHat V_I = \pi^{-1}_I(|V_I|)\subset \HHat U_I$
and $I$ is any subset of $A$.  For $J\subseteq I$ the space of morphisms with source $\HHat V_I$ and target $\HHat V_J$ 
 is given by the restriction of $\HHat\pi_{JI}$ to $\pi_I^{-1}(|V_I|\cap|V_J|)$.  When $I=J$ these are identity maps. The category is completed by adding the inverses of these morphisms. The projection from the space of objects of $\Xx_\Vv$ to its Hausdorff  orbit space is denoted
 $\HHat\pi_{\Hsm}:(X_\Vv)_0\to |X_\Vv|_{\Hsm}$.
  \end{defn}
  
 \begin{figure}[htbp] %  figure placement: here, top, bottom, or page
 %xx   \centerline{\psfig{figure=orbfig4a.jpg,width=3.5in}}
    \centering
   \includegraphics[width=3.5in]{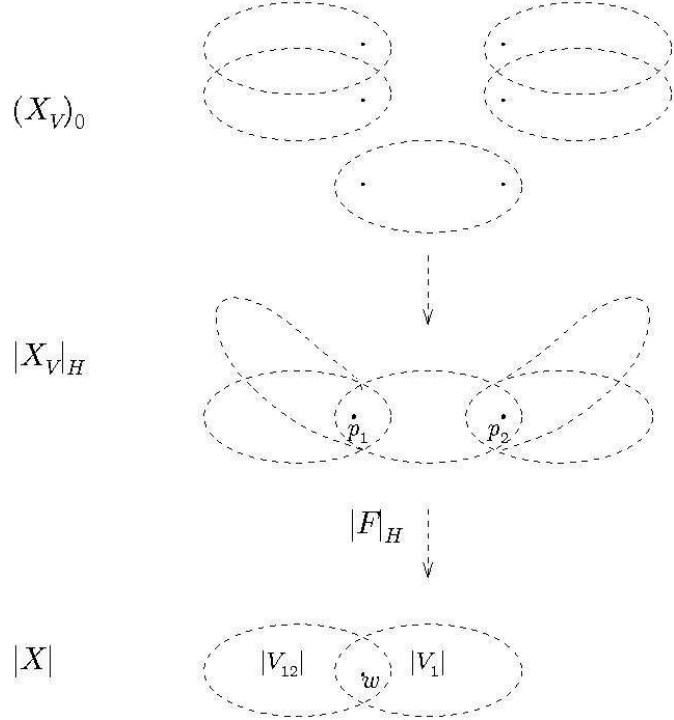} 
    \caption{$(X_\Vv)_0, |X_\Vv|_{H}$ and $|X|$ 
    with $|G_1|=|G_2|=2.$
 The points $p_1,p_2$ project down to $w.$}
    \label{fig4}
 \end{figure}
 
 Observe that we do not need to add the composites of the form $\HHat \pi_{JI}\circ (\HHat\pi_{KI})^{-1}$, for such morphisms would be defined over the intersection $|V_I|\cap|V_J|\cap|V_K|$ which is nonempty only if $K\subset J\subset I$ or $J\subset K\subset I$.
 In the former case this composite is  $(\HHat\pi_{KJ})^{-1}$
 while in the latter it is  $\HHat\pi_{JK}$.  Hence in either case it is already in the category.
 
  Since $\Xx_\Vv$ has at most one morphism between any two objects, it  is a nonsingular sse Lie groupoid.  It need not be proper.
  However, we now show that it does have a weighting.

  \begin{prop}\labell{prop:weigh} Let $\Xx, \Aa$ and $\Vv$ satisfy the conditions in Definition~\ref{def:resol}, and form  
 $\Xx_\Vv$ as above. 
%  Define $\la_\Vv:\pi_0((X_\Vv)_0) \to (0,\infty)$ by
% setting $\la_\Vv(\HHat V_I^\al): = k/|G_I|$ for each $\al$ and 
% $I$, where $k$ is the order of the stabilizer of a generic point in $X_0$.
 Then there is a weighting function $\La_\Vv: |X_\Vv|_{\Hsm}\to (0,\infty)$  such that %the triple
  $(\Xx_\Vv,  \La_\Vv)$ is a wnb groupoid.
  Moreover $(\Xx_\Vv,  \La_\Vv)$ is commensurate with $\Xx$.
  \end{prop}

The proof uses the following technical lemma.   For each component $|V_I|^\al$ of $|V_I|\subset |X|$ we set
$\HHat V^\al_I: = 
\pi^{-1}_I(|V_I|^\al)\cap \HHat U_I$, where $\pi_I:\HHat U_I \to |U_I|\subset |X|$ is the obvious projection.
% To this end, we 
% now investigate the morphisms in this category.  These are denoted 
% by the letter $\HHat\pi$ to distinguish them from the 
% projections with image $|X|$:
% $$
% \pi_I: \HHat U_I\to |X|,\qquad \pi_\ell:\HHat U_\ell \equiv U_\ell\to |X|.
% $$
%Further, for each component $|V_I|^\al$ of $|V_I|\subset |X|$ we set
%$\HHat V^\al_I: = 
%\pi^{-1}_I(|V_I|^\al)\cap \HHat V_I$. 
% 
 
\begin{lemma}\labell{le:wt} Assume $\Xx$ and 
$\Vv$ are as in Proposition~\ref{prop:weigh}. Then  the following statements hold.\SSS

\NI
 {\rm (i)} For each $\al\in \pi_0(|V_I|)$, 
 the group $G_I$ acts transitively on the  components of $\HHat V_I^\al$.
\SSS

\NI
{\rm (ii)}
For each $\ell\in I$ the group $G'_{\ell I}$ acts freely on the components of $\HHat U_I$.  

\NI
{\rm (iii)}   If $I\subset J$ and $q\in \HHat V_I$ projects down to $\pi_I(q)\in |V_I|\cap |V_J|$, there are 
exactly
$|G_J|/|G_I|$ components of $\HHat V_J$ whose image by $\HHat\pi_{IJ}$ contains $q$.
\end{lemma}
\begin{proof}  (i) holds because as we remarked above
 the projection $\pi_I:\HHat V_I^\al\to |V_I|^\al$ quotients out by
  the action of $G_I$.

We shall check (ii) for $\ell = i_1$. (The other cases are similar.) Then the element $ (\ga_k,\dots,\ga_2,1)$ of 
$G'_{i_1I}$ acts on $\HHat U_I$ via the maps
$$
(\de_{k-1},\dots,\de_1)\mapsto 
(\ga_k^{-1} \de_{k-1} \ga_{k-1},\dots, \ga_2^{-1} \de_1).
$$
If the source and target of this map lie in the same component of 
 $\HHat V_I$ then the morphisms $\de_i$ and $\ga_{i+1}^{-1}\de_i\ga_i$ lie in the same component of $X_1$ for each $i$.  
 In particular, $\de_1$ and $\ga_2^{-1} \de_1$ lie in the same component of $X_1$.  Hence composing with $\de_1^{-1}$ we see that $\ga_2\in G_2$ is isotopic to an identity map.  But because 
 $\id:X_0\to X_1$  is a section of the local covering map 
 $s:X_1\to X_0$ the set of identity morphisms form a connected component of $X_1$ in any groupoid $\Xx$.   Hence  we must have $\ga_2= 1$. Repeating this argument we see that $\ga_i=1$ for all $i$.  This proves (ii).
 
 Finally 
 (iii) holds because the set of components of $\HHat V_J$ whose image by $\HHat\pi_{IJ}$ contains $q$ form an orbit of an action by
the group $\prod_{i\in J\setminus I} G_i$ which is free by (ii).  \end{proof}

\NI
{\bf Proof of Proposition~\ref{prop:weigh}.}  Without loss of generality we assume that $\Xx$ is connected and  define 
%%D$k$ 
$\ka$ to be the order of the stabilizer of a generic point in $X_0$. Thus $\ka=1$ iff $\Xx$ is  effective.

We first investigate the relation of  $\Xx_\Vv$ to $\Xx$.  To this end,  
define the functor $F:\Xx_\Vv\to \Xx$ 
 on objects as the projection 
  $$
F:  \HHat V_I\to U_{i_1},\quad (\de_{k-1},\dots,\de_1)\mapsto s(\de_1).
  $$
The morphisms given by $\HHat\pi_{JI}$ are taken to identities if 
$i_1\in J$. If not, the projection $\HHat\pi_{IJ}$ with source
$(\de_{k-1},\dots,\de_1)$ is taken to the composite
$\de_{r-1}\circ\dots\circ\de_1$ where $i_r$ is 
the smallest element in $I\cap J$.  

\begin{figure}[htbp] %  figure placement: here, top, bottom, or page
%xx   \centerline{\psfig{figure=orbfig5.jpg,width=3in}}
   \centering
   \includegraphics[width=3in]{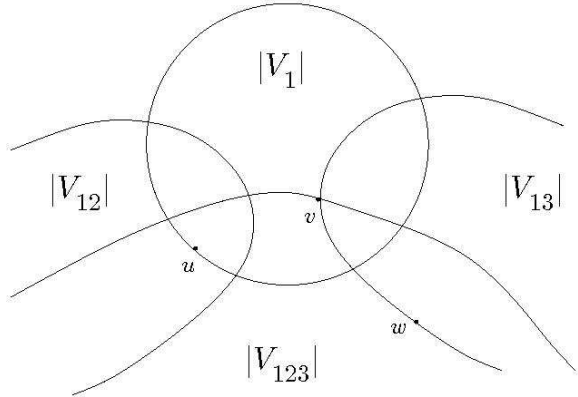} 
   \caption{This is a diagram of part of $|X|$. Here $I_u=\{1\},J_u=\{1,2\}$; $I_v=J_v=\{1\}$; and $I_w=\{1,3\}, J_w=\{1,2,3\}$.}
   \label{fig:5}
\end{figure}

Since $|X|$ is Hausdorff,  the induced map 
$|X_\Vv|\to |X|$ factors through $|F|_{\Hsm}:|X_\Vv|_{\Hsm}\to |X|$.
For each $w\in |X|$, let  $I= I_w$ (resp. $J = J_w$) be the subset of minimal length
  such that $w\in cl(|V_I|)$ (resp. $w\in |V_J|$). (See 
  Fig.~\ref{fig:5}.)
  %, and choose $\al=\al_w, \be=\be_w$ so that 
  %$w\in cl(|V_I|^\al)\cap |V_J|^\be$.
  Note that $I\subseteq J$ because $\Vv$ is $\Aa$-compatible.
  For each $\ell$ such that $w\in |U_\ell|$, 
  $w$ has precisely $|G_\ell|/|G_w|$ preimages $x$ in $U_\ell\subset X_0$, 
  where $G_w$ denotes the  isomorphism class of the stabilizer groups $G_x$.  If $w\in |V_L|$ and $\ell\in L$, then each such $x$ 
  lies in $U_\ell\cap \HHat \pi_{\ell L}(\HHat V_L)$. (Here we identify 
  $U_\ell$
  with  $\HHat U_\ell$.) By Lemma~\ref{le:wt}, each $x$ has
  $|G_{\ell L}'|$ preimages in $\HHat V_L$, each lying in a different component.  Thus
  $w$ has $|G_L|/|G_x|$ preimages in 
  $\HHat V_L$: cf. Fig.~\ref{fig4}.   Moreover the minimality of $J=J_w$ implies that
  $w$ has $|G_J|/|G_x|$ preimages 
  in $|X_\Vv|$.  On the other hand,  
 if  $I=I_w\ne J$ two such preimages map to distinct elements in  $|X_\Vv|_{\Hsm}$ 
   iff the corresponding elements of $\HHat V_J$
    map to distinct elements of $\HHat U_I$ under the projection $\HHat\pi_{IJ}$.  Hence $w$ has $|G_I|/|G_x|$ preimages $p$
  in $|X_\Vv|_{\Hsm}$.

  Now consider $p\in |X_\Vv|_{\Hsm}$ and let $|F|_{\Hsm}(p) = w$.  Define
  $\La_\Vv:|X_\Vv|_{\Hsm}\to \R$ by setting $\La_\Vv(p) = \ka/|G_{I}|$ where
  $I= I_w$ is as above.  Then, because $w$ has $|G_I|/|G_x|$ preimages $p$
  in $|X_\Vv|_{\Hsm}$,  $(|F|_{\Hsm})_*(\La_\Vv)(w) = 1$ provided that $\ka = |G_x|$, i.e. provided that $w$ is smooth.
  Hence $F$ is a layered equivalence, provided that $(\Xx_\Vv,
  \La_\Vv)$ is a branched manifold. 
 
 Thus it remains to check that the conditions of
  Definition~\ref{def:brorb} hold for $(\Xx_\Vv, \La_\Vv)$. 
In the following construction, 
 each local branch at  $p$
  lies in a component of 
  $\HHat V_J$, where $J = J_w$ as defined above. They will all be assigned the same weight $\ka/|G_J|$. 
  The remarks in the preceding paragraph give rise to the following characterization of the fiber in $|X_\Vv|$  over $p\in |X_\Vv|_{\Hsm}$.
\begin{quote}\it{
Define $I=I_w, J=J_w$ as above, where $w=|F|_{\Hsm}(p)$. 
  Since $p\in cl(|\HHat V_I|_{\Hsm})$ there is a convergent sequence $p_n\to p$
  of elements in $|\HHat V_I|_{\Hsm}$.  Lift this sequence to 
  $\HHat V_I$ and
  choose $q\in cl(\HHat V_I)\subset \HHat U_I$ to be one of its limit points.   Then
  the points in the fiber over $p$ in $|X_\Vv|$ are in bijective correspondence with
  the $|G_J|/|G_I| $ distinct elements in $\HHat V_J$ that are taken to $q$ by $\HHat\pi_{IJ}$.}
\end{quote}
\NI
In particular, the elements in the fiber over $p$ 
lie in different components of $|\HHat V_J|$.  
It follows that for any $K$ each component of $\HHat V_K$ maps bijectively  onto a closed subset of the inverse image 
$|\HHat V_K|_{\Hsm}: = (|F|_{\Hsm})^{-1}(|V_K|)$ of $|V_K|$ in $|X_\Vv|_{\Hsm}$. (Note that each component of $\HHat V_K$ injects into $|X_\Vv|$
 because $\Xx_\Vv$ is nonsingular.) 
% $\HHat V_J$,
% Since each of the elements in this fiber lies in a different component of
% $\HHat V_J$, it follows that each component of  $\HHat V_J$, which injects into $|X_\Vv|$ because $\Xx_\Vv$ is nonsingular, also injects into $|X_\Vv|_{\Hsm}$.
% 
% To see they inject into $|X_\Vv|_{\Hsm}$, note that 
%  
%   because for any $K$  each component of $\HHat V_K^\al$ 
%  injects  into $|X_\Vv|_{\Hsm}$ by Lemma~\ref{le:}(ii).
%  Further,

  To construct the local branches at $p$,
  first suppose that $I_w=J_w$.  
%  $p\in |\HHat V_{I}^\al|_{\Hsm}$ for some $\al$ where $I = I_w$ and $w:= |F|(p)\in |X|$. 
  This hypothesis implies that $p\in |\HHat V_{I}^\al|_{\Hsm}$
  for some component $|V_I|^\al$ of $|V_I|\subset |X|$,
 but that  $p\notin cl(|\HHat V_{K}|_{\Hsm})$ for any $K\subsetneq I$.  Take 
  $$
  N_p: =  |\HHat V_{I}^\al|_{\Hsm}\setminus \Bigl(\bigcup_{K\subsetneq I}
  cl(|\HHat V_{K}|_{\Hsm})\Bigr),
  $$
and choose a single local branch $U_1$ equal to the inverse image of
$N_p$ in any component of $\HHat V_{I}^\al$.   (Note that each
  component of $\HHat V_{I}^\al$ surjects onto the connected set
  $|\HHat V_{I}|_{\Hsm}^\al$ since the
  components of $\HHat U_I$ surject onto those of $|U_I|$.) 
  As mentioned above, we set
$m_1: = \ka/|G_I|$. Since
  $\La_\Vv$ equals $ \ka/|G_{I}|$ 
  at all points of $N_p$, the conditions are satisfied in this case.

% Hence 
% there is a unique $q\in C$ such that $\HHat \pi_{\Hsm}(q) = p$.   By
% Lemma~\ref{le:wt}(iv) there are precisely $|G_J|/|G_I|$ components of 
% $\HHat V_J$ whose image by $\Hat\pi_{IJ}$ contains $q$.  

If $p\in \p(|\HHat V_I|_{\Hsm}^\al)$  where 
$I={I_w}$, we choose $N$ to be a connected open
neighborhood  of $p$ in $|X_\Vv|_{\Hsm}$ that satisfies the following
conditions:\SSS

\NI
(i) $N\subset |\HHat V_L|_{\Hsm}$ for all $L$ such that $p\in |\HHat V_L|_{\Hsm}$,\SSS

\NI
(ii) $N$ is disjoint from all sets $cl(|\HHat V_K|_{\Hsm})$ such that $p\notin
cl(|\HHat V_K|_{\Hsm})$.\SSS

\NI
(This is possible because of the local finiteness of $\Vv$ and the
fact that two sets in $\Vv$ intersect iff their closures intersect.)
%Hence $p\in cl(|\HHat V_K|_{\Hsm} \Rightarrow |\HHat V_K|_{\Hsm}\cap |\HHat 
%V_I|_{\Hsm}\ne \emptyset \Rightarrow I\subset K$.)
Now define $L: = L_w\supset I$ to be the {\it maximal} set $L$ such that 
$p\in |\HHat V_L|_{\Hsm}$.  
Choose $p'\in (N\cap |\HHat V_I|_{\Hsm})$ and a
 lift $q'\in \HHat V_I^\al$ of $p'$, i.e. so that 
$\HHat\pi_{\Hsm}(q') = p'$.  Then by Lemma~\ref{le:wt}(iii) there are
precisely $|G_L|/|G_I|$ components of $\HHat V_L$ whose image by
$\HHat\pi_{IL}$ contains $q'$.  Choose their intersections with 
$(\pi_{\Hsm})^{-1}(N_p)$  to be the local branches
at $p$, where $\pi_{\Hsm}:(X_\Vv)_0 \to|X_\Vv|_{\Hsm}$ is the obvious projection. Then the covering and local regularity properties of Definition~\ref{def:brorb} hold by the above discussion of the fibers
of the map $|X_\Vv|\to |X_\Vv|_{\Hsm}$.  The weighting condition
also holds at $p$ since $p$
lies in the image  of each local branch.
 
 It remains to check the weighting condition at the other points
 $p''$  in $N$.  Define $I'': = I_{w''}$,
 the minimal index set $K$ such that $p''\in cl(|\HHat V_K|_{\Hsm})$.  Condition
 (ii) for $N$ implies that $p\in cl(|\HHat V_{I''}|_{\Hsm})$ and hence that
 $I\subseteq I''$.  On the other hand, because $N\subset |\HHat V_L|_{\Hsm} $ the
 minimality of $I''$ implies that
 $I''\subseteq L$.  
By Lemma~\ref{le:wt}(iii) there are precisely $|G_L|/|G_{I''}|$ components of
$\HHat V_L$ whose image in $|X_\Vv|_{\Hsm}$ meets $|\HHat V_{I''}|_{\Hsm}$ near
$p''$.  Since these correspond bijectively to the local branches 
that intersect $(\pi_{\Hsm})^{-1}(p'')$ the weighting
condition  holds for all $p''\in N$. This completes the proof.
\QED\MS
 
In order to show that every orbifold has a resolution we need to see that suitable coverings $\Vv$ do exist. This well known fact is established in the next lemma: cf. Fig.~\ref{fig3}.

\begin{lemma}\labell{le:cov}  Let $Y$ be a normal topological space
 with an
open cover $\{U_i\}$  such that each set 
$U_i$ meets only finitely many other $U_j$. Then
 there are open subsets
$U_i^0\,\sqsubset\, U_i$ and $V_I\subset U_I$ with the following
properties: \SSS

\NI
{\rm (i)}\,  $Y\,\subseteq \,\cup_i\, U_i^0$,\,\,\,\,\,$Y\,\subseteq\,
\cup_I V_I$;\SSS

\NI
{\rm (ii)}\,\,  $V_I\cap U_i^0 = \emptyset$ if $ i\notin I$;\SSS

\NI
{\rm (iii)}\,\,  if $cl(V_I)\cap cl(V_J)\ne \emptyset$ then one of the sets $I,J$
contains the other.
\end{lemma}

\begin{proof} First choose an open covering $\{U_i^0\}$ of $Y$ such that $U_i^0\,\sqsubset\,U_i$ for all $i$.  For each $i$, choose $k_i$
so that $U_i\cap U_J = \emptyset$ for all $J$ such that $|J|> k_i$. 
Then choose for $n = 1,\dots, k_i$ 
open subsets $U_i^n$, $W_i^n$ of $U_i$ such that
$$
U_i^0\,\sqsubset\, W_i^1 \,\sqsubset\, U_i^1
\,\sqsubset\, W_i^2 \,\sqsubset\,\dots \,\sqsubset\,
U_i^{k_i} = U_i,
$$
and set $U_i^j: = U_i^{k_i}$ for $j>k_i$.  Then, if $|I| = \ell,$ define
$$
V_I = W_I^\ell - \bigcup_{J:|J| > \ell} cl(U_J^{\ell+1}).
$$
We claim that this covering satisfies all the above conditions. 
For example, to
 see that $Y\subseteq \cup_IV_I$, observe that
the element $y\in Y$ lies in $V_K$, 
where $K$ is a set of maximal length such that $y\in W_K^{|K|}$.
To prove (iii) first consider the case 
when $y\in cl(V_I)\cap cl(V_J)$ where $|I|=|J|=\ell$, but $I\ne J$.  Then $y\in cl(W^\ell_{I\cup J})\subset U^{\ell+1}_{I\cup J}$; but this is impossible because $U^{\ell+1}_{I\cup J}$ is cut out of $V_I$, which implies that $cl(W^\ell_{I\cup J})$ does not intersect $cl(V_I)$.
The rest of the proof is similar and is left to the reader.
\end{proof}

%In this situation we will say that the cover $\{V_I\}$ is {\bf $\Uu$-subordinate}.\MS

\NI
{\bf Proof I of Proposition~\ref{prop:resol}.}\,
We must show that every orbifold $\un Y$ is commensurate with a tame branched manifold.
 Let $\Xx$ be a groupoid structure on $\un Y$ constructed 
 as before from an atlas.  By applying Lemma~\ref{le:cov} to the covering of $|X|$ by the charts $|U_i|$ of this atlas, we may find a covering $\{|V_I|\}$ of $|X|=Y$ that satisfies the conditions in Definition~\ref{def:resol}.  We may also choose the $|V_I|$ 
so that for one (and hence every) $i\in I$ the pullback
$\pi^{-1}(|V_I|)\cap U_i$ has piecewise smooth boundary in $U_i$.
Now consider the corresponding wnb groupoid $(\Xx_\Vv,\La_\Vv)$ 
constructed in Proposition~\ref{prop:weigh}.   This is  
commensurate to $\Xx$, and its branch locus is piecewise smooth.
Therefore as in Lemma~\ref{le:tame0} it has a tame refinement.
\QED
  
  \begin{rmk}\labell{rmk:resol}\rm (i)   If $\Xx$ itself is
   nonsingular, then the projection $\Xx_\Vv\to \Xx$ is an equivalence.  Indeed $\Xx_\Vv$ is just  the refinement of $\Xx$ corresponding to the covering $\Vv$ of $|X|$.\SSS

\NI
(ii)  
Consider the two noneffective groupoids $\Xx, \Zz$  of 
Remark 2.14. Cover  $|X| = |Z| = S^1$ by three 
open arcs $U_i, i = 1,2,3$, whose triple intersection is empty.  Then if we choose the same covering $\Vv$ in both cases, the groupoids $\Xx_\Vv$ and $\Zz_\Vv$ are isomorphic. See Fig. 7. Thus if one wants to preserve information about the topological structure of the trivially acting morphisms it might be better to define the resolution of a noneffective groupoid using the approach discussed in Remark 3.17.
\end{rmk}

\begin{figure}[htbp] %  figure placement: here, top, bottom, or page
   \centering
   \includegraphics[width=4in]{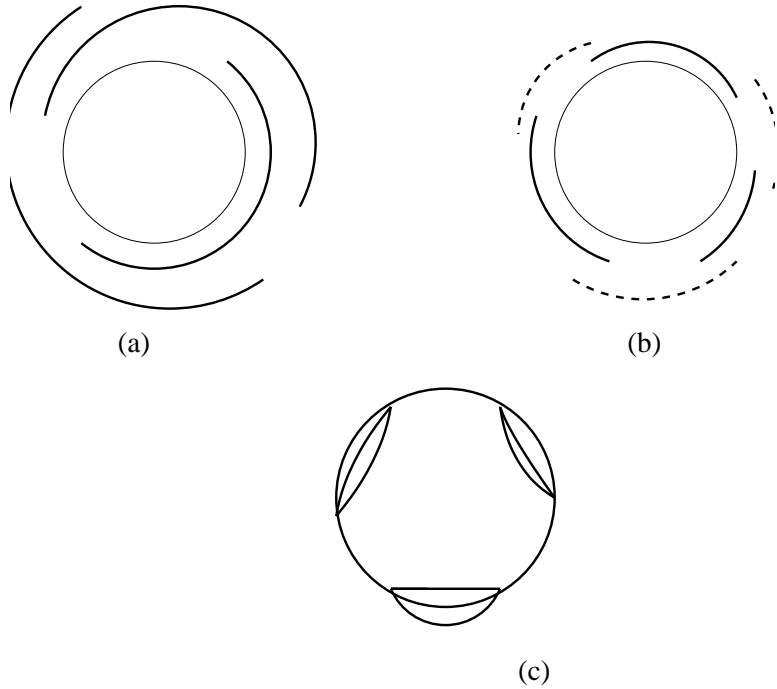} 
   \caption{(a) illustrates the covering $\Uu$ and (b) the covering $\Vv$, in which the $V_i$ are represented by solid arcs and the $V_{ij}$ by dashed arcs. (c) is the resolution $|X_\Vv|_H = |Z_\Vv|_H$.
   The arc is tripled over the parts of the circle not covered by the $V_i$ since each set $\widehat V_{ij}$ has three components.}
   \label{fig:example}
\end{figure}

%(ii)  
%Consider the two nonreduced groupoids $\Xx, \Zz$  of 
%Remark~\ref{rmk:nonred0}. Cover  $|X| = |Z| = S^1$ by three 
%open arcs $U_i, i = 1,2,3$, whose triple intersection is empty.
%Then if we take $|V_I| = |U_I|$  for all $I$, the groupoids $\Xx_\Vv$ and $\Zz_\Vv$ are isomorphic.  Thus if one wants to preserve information about the topological structure of the trivially acting morphisms it might be better to define the resolution of a nonreduced groupoid using the approach discussed in Remark~\ref{rmk:nonred}.
%\end{rmk}

\NI
{\bf Compatibility with operations in SFT.}\,\,
In Hofer--Wysocki--Zehnder~\cite{HWZ}, it is important that all  constructions are compatible with the natural operations of
 symplectic field theory. These arise from the special structure
 of the groupoids they consider, which live  
in a category of manifold-like objects called polyfolds
with boundary and corners.  Though polyfolds 
 are infinite dimensional, their corners have
 finite codimension.  Moreover, each groupoid  $\Xx$ has a finite
 formal dimension (or index), and its boundary
 is build up inductively from pieces of lower index.  For example, in the simplest situation where there is a boundary but no corners, the  boundary $\p\Xx$ of $\Xx$ decomposes as
a disjoint 
union of products $\Xx_1\times \Xx_2$ of ep groupoids $\Xx_\al$ of lower index (and without boundary).  Hence it is important that the operation of constructing a resolution is compatible with taking products and also extends from the boundary of $\Xx$ to $\Xx$ itself.

We now discuss these questions in the context of finite dimensional manifolds. First consider  products.
The construction of the resolution $\Xx_\Vv\to \Xx$ is determined by the choice of good atlas $\Aa$ (together with an ordering of the set $A$ of  charts)
 and the choice of a subcovering $\Vv$.  
 If $\Aa_\al = \{(U_{\al i}, G_{\al i}, \pi_{\al i}): i\in A_\al\}$ is an atlas for $\Xx_\al$, where $\al=1,2$, then there is a product  atlas 
  $\Aa_1\times \Aa_2$ for $\Xx_1\times \Xx_2$
with charts $(U_{1i}\times U_{2j}, G_{1i}\times G_{2j}, \pi_{1i}\times \pi_{2j})$ where $(i,j)\in A_1\times A_2)$. 
But the product covering by sets $|V_{1I}|\times |V_{2J}|$ is not
$\Aa_1\times \Aa_2$-compatible
since there will be nonempty intersections of the form
$\bigl(|V_{1I}|\times |V_{2J}|\bigr)\cap \bigl(|V_{1I'}|\times |V_{2J'}|\bigr)$ where $I\subsetneq I'$ and $J'\subsetneq J$.\footnote
{
The fact that there is  no  ordering of the new index set
$A_1\times A_2$ that is symmetric in $\al$ is less serious, since one could define the resolution without using such an order.  
For example, the objects of the new resolution would contain  
$|I|!$ copies of  $\HHat V_I$, one associated to each possible ordering of the index set $I$.}  
Hence our construction does not commute with taking the product.
On the other hand, the following lemma shows that we can take
the resolution of a product to be the product of resolutions of the factors.
 
 \begin{lemma}  For any resolutions $F_\al: (\Bb_\al, \La_\al)\to \Xx_\al$ for $\al=1,2$ 
 the product 
 $$
 F_1\times F_2: (\Bb_1\times \Bb_2, \La_1,\times \La_2) \to 
 \Xx_1\times \Xx_2
 $$
  is a 
 resolution.
 \end{lemma}
 \begin{proof}  This is immediate.
 \end{proof}
 
We now show how to extend a resolution from $\p\Xx$ to $\Xx$.
Note that $\p\Xx$ has  collar neighborhood in $\Xx$ that is
 diffeomorphic\footnote
 {Two groupoids $\Xx,\Yy$ are diffeomorphic if there is a functor $F:\Xx\to \Yy$ that is a diffeomorphism on objects and morphisms.}
  to the product $\p\Xx\times \Ii'$, 
   where $\Ii'$ 
 denotes the trivial groupoid with objects  $(-1,0]$
 as in Proposition~\ref{prop:commen}(iii).  Further $\Xx$ is diffeomorphic to the groupoid $\Xx'$ obtained by 
  extending the collar neighborhood of $\p\Xx$ by $\p\Xx\times \Ii$,
  where $\Ii$ is the trivial groupoid with objects  $[0,1]$.
  
 \begin{lemma}  Let $\Xx$ be an ep groupoid with boundary. Then
any resolution $F: \Bb\to \p \Xx$  extends to a resolution 
 $F': \Bb'\to \Xx$ of $\Xx$.
 \end{lemma}
 \begin{proof}   Construct
 a resolution $\Xx_\Vv\to \Xx$.  As in Proposition~\ref{prop:commen}(iii)
 there is a layered covering  $\Ww\to \p\Xx\times \Ii$ 
 that restricts to $\Bb\to \p \Xx$ over $\p\Xx\times \{1\}$ and to 
 $\p(\Xx_\Vv)\to \p\Xx$ over  $\p\Xx\times \{0\}$. Now take $\Bb'$ to be the union of $\Xx_\Vv$ with $\Ww$ and define $F'$ to be the induced
 functor $\Bb'= \Xx_\Vv\cup \Ww\to \Xx\cup(\p\Xx\times \Ii) = \Xx'\cong\Xx.$
\end{proof} 
   
%%%%%%%%%%%%%%%%%%%%%%%%%%%%%%%%%%%%%%%%%%%%%%%%%%%%%%%%%%
\subsection{Orbibundles and multisections}\labell{ss:multi}  
%%%%%%%%%%%%%%%%%%%%%%%%%%%%%%%%%%%%%%%%%%%%%%%%%%%%%%%%%%

We now show how branched 
%%Dorbifolds 
manifolds arise  as multivalued sections of (unbranched) orbibundles.  This is a groupoid version of the constructions in
Cieliebak--Mundet i Riera--Salamon~\cite{CRS} and 
Hofer--Wysock--Zehnder~\cite{HWZ}.

Recall from Moerdijk~\cite[\S5]{Moe} that an {\bf orbibundle} $\un E\to \un Y$
over an orbifold $\un Y$ is given by an equivalence class of functors
$\rho:\Ee\to \Xx$, where $\Xx$ is a orbifold structure on $Y$ and $\Ee$ is 
a groupoid constructed as follows.  The objects $E_0$ of $\Ee$ form a vector bundle $\rho_0:E_0\to X_0$ and its space $E_1$ of
morphisms arise from a right action of $\Xx$ on $E_0$.  In other words, $E_1$ is the fiber product 
$E\times_{X_0}X_1: = \{(e,\ga): \rho(e) = t(\ga)\}$, the target map is
$(e,\ga)\mapsto e$ and
the source map $$
\mu: E_0\times_{X_0}X_1\to E_0, \quad (e,\ga)\mapsto e\cdot\ga
$$
satisfies the obvious identities: namely
the diagram
$$
\begin{array}{ccc}
E_0\times _{X_0}X_1&\stackrel{\mu}\to& E_0\\
\rho\downarrow&& \rho\downarrow\\
X_0\times_{X_0}X_1&\stackrel{s'}\to & X_0
\end{array}
$$
commutes (here $s'(t(\ga),\ga))= s(\ga)$), the identities $\id_x\in X_1$ act
by the identity and composition is respected, i.e. $(e,\ga\de) = (e,\ga)\cdot\de$. 

This construction applies when given any  sse groupoid $\Xx$ and any vector bundle $E_0\to X_0$ 
that supports a right action of $\Xx$. We shall call $\rho:\Ee\to \Xx$ a vector bundle over $\Xx$.
Clearly, $\Ee$ is an ep groupoid if $\Xx$ is.  Also, 
if $(\Xx,\La)$ is a weighted branched manifold, we can pull $\La$ back to give a weighting on $\Ee$. Hence $\Ee$ is a wnb groupoid if $\Xx$ is. 

\begin{example}\labell{ex:bund}  \rm  (i)
Let $E_0= TX_0$, the tangent bundle of $X_0$.  Then $\Xx$ acts on $E_0$ in the obvious way.  The resulting bundle $T\Xx\to \Xx$ is
called the tangent bundle.  Note that $\Xx$ acts effectively on $TX_0$ iff $\Xx$ is effective. \SSS

\NI(ii) Every bundle $\Ee\to \Xx$ has an effective reduction $\Ee_{\it eff}\to\Xx'$ where $\Xx'$ is the quotient of $\Xx$ by the subgroup of $K\cong K_y$ that acts trivially on $E_0$. (For notation, see Lemma~\ref{le:red2}.)    Thus the functor $\Xx\to \Xx_{\it eff}$ factors through $\Xx'$.
\end{example}

A {\bf section} of the bundle $\rho: \Ee\to \Xx$ is a smooth functor
$\si:\Xx\to \Ee$ such that $\rho\circ\si = id$. 
In particular, the restriction of each such functor to the space of objects is a (smooth) section $\si_0$ of the vector bundle $\rho_0: E_0\to X_0$. 
The conditions imposed by requiring that $\si_0$ extend to a functor imply that $\si_0$ must descend to a section
of $|E|\to |X|$. 

    We denote by $\Ss(\Ee)$ the space of all such sections
    $\si$.  Since each section $\si$ is determined by the map
    $\si_0:X_0\to E_0$, we may identify $\Ss(\Ee)$ with a subset 
    $\Ss_0$ of the space $\Sect(E_0)$ of smooth sections of this      
    bundle with the usual Fr\'echet topology.\footnote
    {
Since the components of $X_0$ are noncompact, we take the topology given by uniform $C^k$-convergence on compact sets for all $k\ge 1$.
The induced topology on $\Ss(\Ee)$ is not very satisfactory, but it is good enough 
for the present purposes.}
Let us now suppose that $\Xx$ is compact.  By replacing $\Xx$ by an equivalent groupoid if necessary, we may
 suppose that $X_0$ and $X_1$ have finitely many components.  Hence the
    subset $\Ss_0$ is defined by a finite number
of smooth compatibility conditions and so is a submanifold of $\Sect(E_0)$.

\begin{defn}\labell{def:trans}  A section $\si:\Xx\to \Ee$ is said to be {\bf transverse to the zero section} if its image $\si(X_0)$ intersects  the zero set transversally and if for each intersection point
 $x\in X_0$ the induced map $d\si(x):T_\si\Ss(\Xx) \to E_x$ is surjective, 
where $E_x$ is the fiber of $E_0\to X_0$ at $x$ and $d\si(x)$ is the composite of the derivative of $\si$ with evaluation at $x$.
\end{defn}

We claim that if the base $\Bb$ is sufficiently nice there are enough sections.

\begin{lemma}\labell{le:sect0} Suppose that $\Ee\to \Bb$ is a bundle over a tame   wnb groupoid.  Then a generic element of
 $\Ss(\Ee)$ 
intersects the zero section transversally.
\end{lemma}
\begin{proof} We just sketch the proof since the techniques are standard.  The idea is the following. 
Consider an  open set $N_p\subset |B|_{\Hsm}$ with local branches $U_i^p, i\in A_p$ as in Definition~\ref{def:brorb}.  
As in Lemma~\ref{le:part},
the tameness condition
 allows one to construct enough sections to prove the result for
 the restriction of $\Ee$ to $\Bb_N$, the full subcategory of $\Bb$ with objects $U_i^p, i\in A_p$.  But any section with compact support in $\Bb_N$ extends uniquely to a section of $\Ee\to \Bb$.
The result follows.
\end{proof}

If $\Bb$ is a  wnb groupoid, $\Ee$ is oriented and the section $\si:\Bb\to \Ee$ is transverse to the zero section,
then the zero set $\Zz(\si)$ of $\si$ inherits a natural structure as 
a  wnb groupoid.  
   Standard arguments shows that its cobordism class is independent of the section chosen.  In particular, if 
the fiber dimension of $\Ee$ equals $\dim \Bb$ the  Euler number
$\chi(\Ee)$ can be calculated as the number of zeros of a generic section of $\Ee\to \Bb$, or equivalently as the number of points in the zero-dimensional wnb groupoid $\Zz(\si)$;  
cf Example~\ref{ex:brorb}(ii).  For example, the tangent bundle $T\Bb\to \Bb$ of the resolution of the teardrop orbifold described in Example~\ref{ex:tear} has a section with one positively oriented
 zero in each of $D^+, D^-$.   It follows that $\chi(Y) = 1+1/k$.
 
 \begin{defn}\labell{def:euler} The  {\bf (homology) Euler class} of a $d$-dimensional bundle $\Ee\to \Bb$ over an $n$-dimensional wnb groupoid is the 
singular homology class $e(\Ee)\in H_{n-d}(|B|_{\Hsm};\R)$ represented 
 by the image in $|B|_{\Hsm}$ of the fundamental class $[Z(\si)]
 \in H_{n-d}(|Z|_{\Hsm}; \R)$ 
 of the zero set $\Zz(\si)$ of a generic section. 
 \end{defn}
 
  This definition has the expected functorial properties  and is consistent with standard definitions; cf. 
  Proposition~\ref{prop:euler} and the discussion in 
  \S\ref{ss:euler} below.
  Moreover, if $\Bb$ has rational weights then 
$e(\Ee)\in H_{n-d}(|B|_{\Hsm};\Q).$

\begin{rmk}\labell{rmk:trans}\rm
 The condition on $d\si(x)$ in 
Definition~\ref{def:trans} ensures the existence of enough local deformations of $\si$ to have a good transversality theory, e.g.
one in which the zero sets of two generic sections intersect transversally.  
As an example,
consider the tangent bundle $T\Xx\to \Xx$ of the teardrop orbifold.
Every section of this bundle 
vanishes at the singular point $p$, and so
 no section is transverse in the sense of the above definition.
\end{rmk}

One way to deal with the lack of sections over singular $\Xx$
is to consider multivalued sections.
In Cieliebak--Mundet i Riera--Salamon~\cite{CRS} and Hofer--Wysock--Zehnder~\cite{HWZ} these are defined by means of the 
characteristic function of their graph. The following definition is a mild adaptation of that in~\cite{HWZ}. Although one could define multisections of bundles over wnb groupoids  we shall be content here with the case when the base is unbranched. Further, to be compatible 
with~\cite{HWZ} we assume 
from now on that all weighting functions take rational values.

\begin{defn}\labell{def:multi}  Let $\rho:\Ee\to \Xx$ be a bundle over an ep groupoid
and denote by
$\Qq^{\ge 0}$  the category with objects the nonnegative rational numbers and only identity morphisms.
A {\bf multisection} of $\Ee\to \Xx$
is a smooth functor $\Ll_S:\Ee\to \Qq^{\ge 0}$
that has the following local description: for
each point $x\in X_0$ there is an open neighborhood $U$ and a finite nonempty set  of smooth local sections $\si_{j}:U\to E_0$ with positive rational weights $m_{j}$ such that for all $x\in U$
$$
\Ll_S((x,e)) =\underset{j:\si_{j}(x) = e}\sum m_{j},\qquad
\underset{e\in E_x}\sum \Ll_S((x,e))=1,
$$
where by convention the sum over the empty set is $0$.
The triple $(U, \si_j, m_j)$ is called a {\bf local section structure for $\Ll_S$.}
We denote by $\Ss_m(\Ee)$  the set of all multisections of $\Ee\to \Xx$.
\end{defn}

The first condition implies that for each $x\in X_0$ there are a finite number of elements $e\in E_x$ such that $\Ll_S((x,e))\ne 0$.
The set of such elements $(x,e)\in E_0$ is called the {\bf support} of $\Ll_S$.
The second condition is equivalent to requiring that 
$\sum m_j = 1$ and implies
 that $\Ll_S$ has total weight $1$.  
The sum of two multisections $\Ll_S,\Ll_T$ is given  by the convolution product:
 $$
(\Ll_S + \Ll_T) (x,e): = \underset{e'+e'' = e}\sum \Ll_S(x,e')\Ll_T(x,e'').
$$
Similarly, $r\Ll_S(x,e): = \Ll_S(x,re)$ for $r\in \Q$.  Hence
multisections form a $\Q$-vector space.
Note also that because $\Ll_S$ is a functor it takes the same value on all 
 equivalent objects in $E_0$ and so descends to a 
 function on $|E|_{\Hsm}$. Hence when counting the number of zeros of a multisection, one should count the equivalence classes $c$ lying in the zero section, each weighted by
 the product $\oo_c\Ll_S(c)$, where $\oo$ is the orientation
  (cf. Example~\ref{ex:brorb} (ii)). 
 
 We say that $\Ll_S$ is {\bf single valued}  if there is just one local section $\si_j$ over each open set $U$.  In this case the support of $\Ll_S$ is the image of a single valued section $\si_0:X_0\to E_0$
 that extends to a functor $\si:\Xx\to \Ee$, i.e. a section as we defined it above. However, if $\Xx_{\it eff}$ is singular 
% then no orbibundle
%%D over $\Xx$  can support single valued sections.
the values of single valued sections at nonsmooth points  are restricted.
For example, they must vanish at $x\in X_0$ if every point 
of $E_x$ is moved by some element of $G_x$. 

 We now show that  the support of a multisection $\Ll_S$ is (under a mild hypothesis) the image of a wnb groupoid. 
We begin with a preliminary lemma that
 explains the relation between two different
  local section structures.

 \begin{lemma} \labell{le:locsec}
Consider two local  section structures $(U, \si_j, m_j), (U', \si_j', m_j')$  for $\Ll_S$, and let 
 $x\in U\cap U'.$  For each pair $i,j$ define $V_{ij}: = \{y\in U\cap U':
 \si_i(y) = \si_j'(y)\}$.  
Then for all $i$ there is $j$ such that 
% $(x,\si_i(x)) \in cl\bigl(\Int(\im \si_i\cap \im \si_j')\bigr)$.
$x\in cl\bigl(\Int V_{ij})\bigr)$.
 \end{lemma}
 \begin{proof} If not, there is an open neighborhood $V\subset U\cap U'$
  of $x$ in $B_0$ such that $V\cap \Int V_{ij} = \emptyset$ for 
  all $j$.
  %, where $A_j: = \si_i(U)\cap \im \si_j'(U')\subset E_0$.  
  But
$V\subset \cup_j V_{ij}$ and each set  $V\cap V_{ij}$ is relatively closed in $V$.  Since there are only finitely many $j$,
 at least one intersection $V\cap V_{ij}$ must have  nonempty interior.  But $\Int(V\cap V_{ij}) = V\cap \Int(V_{ij})$ since $V$ is open, so this  contradicts the hypothesis.
 \end{proof}
  
%  We now explain how this definition fits in with the concepts introduced earlier.
  
\begin{prop}\labell{prop:constr}  Let $\rho:\Ee\to \Xx$ be an oriented bundle over an ep groupoid $\Xx$.  Then:\SSS

\NI
{\rm (i)}
The support of any multisection $\Ll_S:\Ee\to \Qq^{\ge 0}$ 
 is the image of an
sse groupoid $\Yy_S$ by a functor $\Si_S:\Yy_S\to \Ee$
such that $\Ll_S$ pulls back to a (single valued) section $\si_S$ of the bundle
$F_S^*(\Ee)\to \Yy_S$, where $F_S: = \rho\circ \Si_S: \Yy_S\to \Xx$, 
$$
\xymatrix
{
F_S^*{\Ee} \ar[d] \ar[r]&
\Ee\ar[d]_{\rho}\\
\Yy_S \ar@/^1.1pc/[u]^{\si_S}\ar[ur]^{\Si_S\!\!}\ar[r]^{\,\,\,F_S}&{\Xx}.
}
$$

\NI
{\rm (ii)} The groupoid $\Yy_S$ is nonsingular iff no open subset of the support of $\Ll_S$
 is contained in the singular set $E_0^{sing}: = \{(x,e)\in E_0: |G_{(x,e)}|>1\}$. \SSS
 
 \NI
{\rm (iii)} If $\Yy_S$ is nonsingular, it may be given the structure of a wnb groupoid
so as to make  $F_S:\Yy_S\to  \Xx$ a layered covering.
\end{prop}
\begin{proof}
 To see this, choose a locally finite covering of $X_0$ by sets $U^\al,\al\in A,$ 
 with the properties of  Definition~\ref{def:multi}, 
 and then define $\Yy_S$ to be the category
 whose space of objects is the  disjoint union of copies of the sets $U^\al$, with one copy $U_{i}^\al$ of $U^\al$  for each section $\si_{i}$.  These 
 sections give a smooth immersion $\Si_S:Y_0\to E_0$ that is injective 
 on each component of the domain.  We define the morphisms in $\Yy_S$ to be the pullback by $\Si_S$ of the morphisms in the 
{\it interior} of the full subcategory of $\Ee$ with objects $\Si_S(Y_0)$.
(For example, if two local sections agree at an isolated point
 we ignore the corresponding morphism.  We define $Z_1$ to be the set of all such ignored morphisms.) Then $\Yy_S$  is an sse groupoid. The rest of (i) is clear, as is  (ii).

Because $|E|$ is Hausdorff, the functor  $\Si_S$ induces a continuous map 
 $$
 |\Si_S|_{\Hsm}: |Y_S|_{\Hsm}\to |\supp(\Ll_S)|\subset |E|.
 $$
It is injective over the open dense set $|V|_{\Hsm}: = |Y_S|_{\Hsm}\setminus |s(Z_1)|_{\Hsm}$ (where $s$ denotes the source map) because the restriction of $\Si_S$ to the full subcategory of 
$\Yy_S$ with objects $(Y_S)_0\setminus s(Z_1)$ is a bijection onto a
 full subcategory of $\Ee$.

If $\Yy_S$ is nonsingular, 
we define the local branches over the points in 
$$
N_U: = 
(|F_S|_{\Hsm})^{-1}(|U|)\subset |Y_S|_{\Hsm}
$$
 to be the
sets $U_i$ with weights $m_i$, and define 
$$
\La_Y(p): = |\Ll_S|\bigl(|\Si_S|_{\Hsm}(p)\bigr)\;\in\; (0,\infty)\cap {\mathbb Q}\subset |Q^{\ge 0}|.
$$
%
%and define $\La_Y$ on $N_U$
%by the weighting condition in Definition~\ref{def:brorb}.  To check
%that $\La_Y$ is well defined on the intersection of two different sets $N_1, N_2$ of the form $N_U$ note first that its pushforward by $|\Si_S|_{\Hsm}$  
%equals $|\Ll_S|$. Hence $\La_Y$ is well defined on the open dense set $|V|_{\Hsm}$.  
%But the construction of $\La_Y$ shows that it is determined by its values on $|V|_{\Hsm}$.
 It is now straightforward to check that $(\Yy_S,\La_Y)$ is a wnb groupoid.  In particular each $U_i^\al$ injects into $|Y_S|_{\Hsm}$ because the composite 
 $$
 U_i\stackrel{\pi_{\Hsm}}\longrightarrow |Y_S|_{\Hsm}\stackrel{|F_S|_{\Hsm}}\longrightarrow |U|,
 $$
 is injective.
 
  To see that $F_S:\Yy_S\to \Xx$ is a resolution, note first that 
  the covering property is immediate and that
  the weighting property holds because $\sum_{e\in E_x}\Ll_S((x,e))=1$. To see that $|F_S|_{\Hsm}$ is proper, it suffices to show that every sequence $\{p_k\}$ in $|Y_S|_{\Hsm}$ whose image
by $|F_S|_{\Hsm}$ converges has a convergent subsequence.  Since the covering $U_i^\al$ of $(Y_S)_0$ is locally finite, we may pass to a subsequence of  $\{p_k\}$ (also called  $\{p_k\}$) whose elements all lie in the same set $|U_i^\al|_{\Hsm}$.  Then, for each $k$, there is $y_k\in  U_i^\al$ such that $\pi_{\Hsm}(y_k) = p_k$.
 Now choose $x_k, x_\infty\in X_0$ so that $|x_k|: = \pi(x_k) = 
 |F_S|_{\Hsm}(p_k)$ and 
  $|x_\infty|: = \pi(x_\infty)\in |X|$ is the limit of $\{|x_k|\}$. If $|x_\infty|\in |U^\al|$, then, 
  because the map $U_i^\al\to |U^\al|\subset |X|$ is a diffeomorphism, the 
  sequence $\{y_k\}$ converges in $U^\al_i$ to the point $y_\infty$ corresponding to $x_\infty$.  Hence $\{p_k\}$ has the limit $\pi_{\Hsm}(y_\infty)\in |Y_S|_{\Hsm}$ as required. 
  
  Otherwise, choose
 a local section structure  $(U^\be,\si_j',m_j)$ for $\Ll_S$ near $x_\infty$.  Then $x_\infty\in cl(U^\al\cap U^\be)$, and we may suppose that  $x_k\in U^\be$ for all $k$.  By applying
  Lemma~\ref{le:locsec} to each point $y_k\in U_i^\al$ and passing to a further subsequence, we may suppose that there is $j$ such that
  $y_k\in cl\bigl(\Int V_{ij})\bigr)$ for all $k$, where 
  $$
  V_{ij}: = \{y\in U_i^\al\cap U_j^\be: \si_i(y) = \si_j'(y)\}.
  $$
 Therefore for each $k$ there is $z_k\in U_j^\be$ such that $y_k\approx z_k$.  It now follows from  Lemma~\ref{le:approx} that
 $\pi_{\Hsm}(z_k) = \pi_{\Hsm}(y_k) = p_k$.  But now $\{z_k\}$ has a limit in $U_j^\be$.  Hence $\{p_k\}$ does too. This completes the proof.
\end{proof}
 
Conversely,  suppose that $F:(\Bb,\La_B)\to \Xx$ is a resolution of an 
ep groupoid, and let $\Ee\to \Xx$ be a bundle.  Then any (single valued) section $\si_S:\Bb\to F^*(\Ee)$ of the pullback bundle 
pushes forward to 
a functor $\Si_S: \Bb\to \Ee$.  Note that
 its image is not in general the support of 
a multisection in the sense of Definition~\ref{def:multi} since it  
need not contain entire equivalence classes. 
%We say that two such functors
% $\phi_S: \Bb\to \Ee$ and $\phi'_T: \Bb'\to \Ee$ are equivalent if 
% there are layered coverings  
%$F:\Bb''\to \Bb, F':\Bb''\to \Bb'$ such that 
%$\psi_S\circ F = \psi_T\circ F'$.
However, because $\Ee$ is proper the induced map $|\Si_S|:|B|\to |E|$ factors through $|\Si_S|_{\Hsm}:|B|_{\Hsm}\to |E|$.  

\begin{lemma}\labell{le:multi1} Let $\Ee\to \Xx$ be a   bundle
over the ep groupoid $\Xx$ and $F: (\Bb,\La_B)\to \Xx$ be
a resolution. Each section $\si_S:\Bb\to F^*(\Ee)$
of the pullback bundle gives rise to a multisection $\Ll_S:\Ee\to \Qq^{\ge0}$ where
$$
\Ll_S(x,e): = \sum_{p\in |B|_{\Hsm}\,:\, |\Si_S|_{\Hsm}(p)=|(x,e)|}\;\La_B(p),
$$
and $\Si_S:\Bb\to \Ee$ is the composite of $\si_S$ with 
the push forward $F_*:F^*(\Ee)\to \Ee$. 
 \end{lemma}
%
%Suppose that $\Xx$ is locally finite.  
%Then $\Ll_S$ is a multisection of $\Ee\to \Xx$.
%%that depends only on the equivalence class of $\phi_S: \Bb\to \Ee$.
%\end{lemma}
\begin{proof}  The definition implies that 
$\Ll_S(x,e) = \Ll_S(x',e')$ whenever $(x,e)\sim (x',e')$ in $\Ee$.  Hence $\Ll_S$ is a functor.  
It has the required local structure at $x\in X_0$ by 
 Lemma~\ref{le:brstr1}. 
 \end{proof}

\begin{rmk}\labell{rmk:multi}\rm (i) The above construction of a 
multisection $\Ll_S$ of $\Ee\to \Xx$  from a section $\si_S$ of $F^*(\Ee)\to \Bb$
may be described 
 more formally as follows.
   Consider the (weak) fiber product $\Bb':= \Bb\times_\Xx\Xx$.
 This is a wnb groupoid by Lemma~\ref{le:comm} (ii), and there is 
  a layered equivalence $G': (\Bb', \La') \to (\Bb,\La)$.  
Consider the 
%homotopy commutative 
diagram
$$
\begin{array}{ccc}\Bb': =  \Bb\times _\Xx\Xx&\stackrel{F'}\longrightarrow& \Xx\vspace{.05in}\\
G'\downarrow\;\;\;\;&& =\downarrow\;\;\;\;\\
\Bb&\stackrel{F}\longrightarrow& \Xx.
\end{array}
$$
  The section $\si_S:\Bb\to F^*(\Ee)$ pulls back to a section 
$\si_S':\Bb'\to (F')^*(\Ee)$ that gives rise to a functor 
$\Si'_S:\Bb'\to \Ee$ whose image is precisely the support of $\Ll_S$.
(The effect of passing to $\Bb'$ is to saturate the image of $\Si_S$ under $\sim$.
Note that the two functors $\Si_S':\Bb'\to \Ee$ and
$\Si_S\circ G':\Bb'\to\Bb\to \Ee$ do not coincide, because the diagram only commutes up to homotopy.)  
Thus the two approaches give rise to essentially the same 
multisections. To distinguish them, we shall call  $\Si_S$ a
{\bf wp multisection}. 

One might think of $\Si_S$ as a stripped down version of $\Ll_S$, with inessential information removed.  For example, in the case 
of the teardrop with resolution $F:\Bb\to \Xx$, the pushforward of a section of $F^*(\Ee)\to \Bb$ is single valued over each component
of $X_0$, while the corresponding $\Ll_S$ is multivalued over $D_+$.
\SSS

\NI
(ii)   Suppose that $\Ll_S, \Ll_T$ give rise as above to 
the wnb groupoids $(\Yy_S,\La_S)$ and $(\Yy_T,\La_T)$.
Then it is not hard to check that their sum $\Ll_S + \Ll_T$ 
gives rise to the fiber product
$\Yy_S\times_{\Xx}\Yy_T$.  On the other hand, if $\Yy_S = \Yy_T$
there is a simpler summing operation given by adding the 
corresponding sections $\si_S, \si_T$ of $F^*(\Ee)\to \Yy_S$.
 \end{rmk}

 We say that a multisection $\Ll_S$ of $\Ee\to\Xx$ is {\bf tranverse to the zero section} if it is made from
 local sections $\si_i:U_i\to E_0$ that are transverse to the zero section.   It is easy to check that this is equivalent to saying that the corresponding single valued section $\si_S$ of $F^*(\Ee)\to \Yy_S$ is transverse to the zero section.  Hence the intersection 
 of $\si_S$ with the zero section
 has the structure of a wnb groupoid $(\Zz_S,\La_S)$ as in the discussion after Lemma~\ref{le:sect0}. 
 
 \begin{defn}\labell{def:euler1}  Let $\Ee\to\Xx$
 be a $d$-dimensional vector bundle $\Ee\to\Xx$ over an 
 $n$-dimensional ep groupoid.  If $\Ee$ is effective, we define its 
 {\bf (homology) Euler class} to be  the singular homology 
 class $e(\Ee)\in H_{n-d}(|X|;\Q)$ represented by the image under the composite map $\Zz_S\to \Yy_S\to \Xx$ of
 the fundamental class $[Z_S]\in H_{n-d}(|Z_S|_{\Hsm})$ of the
 zero set of a generic multisection $\Ll_S$. 
  In general,
  %%Df $\Ee$ is noneffective, 
  we define $e(\Ee)\in H_{n-d}(|X|;\Q)$ to be the Euler class of the corresponding effective bundle $\Ee_{\it eff}\to \Xx'$ 
(cf. Example~\ref{ex:bund} (ii)).
 \end{defn}

 It is not hard to prove
 that any two wnb groupoids 
 $(\Zz_S,\La_S), (\Zz_T,\La_T)$  constructed in this way from
 multisections  $\Ll_S, \Ll_T$ are cobordant
  since the pullbacks of $\si_S$ and $\si_T$ to 
  $\Yy_S\times_{\Xx}\Yy_T$ are homotopic. The next result states that Definitions~\ref{def:euler} and~\ref{def:euler1} are consistent.
  
  \begin{prop}\label{prop:euler}  Let $F:(\Bb,\La_B)\to \Xx$ be any resolution of the  $n$-dimensional
  ep groupoid $\Xx$ and let $\rho: \Ee\to\Xx$  be any 
   $d$-dimensional bundle.  Then 
   $$
   (|F|_{\Hsm})_*\bigl(e(F^*(\Ee)\bigr) = e(\Ee) \;\in\; H_{n-d}(|X|;\R).
   $$
    Moreover, if $\La_B$ takes rational values,
   this equality holds in $H_{n-d}(|X|;\Q)$.
   \end{prop}
  
\begin{proof}
This is an immediate consequence of Proposition~\ref{prop:fclass}; the details of its proof are left to the reader.
\end{proof}

 We conclude this section with some constructions. 
  First we explain how to 
 construct \lq\lq enough" wp multisections using the resolution $\Xx_\Vv\to \Xx$.   We shall explain this in the context of
 Fredholm theory and so shall think of the fibers of $E_0\to X_0$ as  infinite dimensional.  We start from a Fredholm section 
 $f:\Xx\to \Ee$ such that over each chart $U_i\subset X_0$ there is a (possibly finite dimensional) space $\Ss_i$ of local sections
 $U_i\to E_0$ that is large enough to achieve transversality over $U_i$.
 We then construct a
  vector space $\Ss$ of wp multisections with controlled branching that is large enough for global  transversality.  It is finite dimensional if the $\Ss_i$ are.
  This is an adaptation 
  of a  result 
 in Liu--Tian~\cite{LiuT} and was the motivation for their 
 construction of the resolution.\footnote
 {
 The topological aspects of their construction were explored earlier in McDuff~\cite{Mcv}.  However the current approach using groupoids, when combined with the Fredholm theory of polyfolds, allows
 for a much cleaner treatment.} 
  
 \begin{prop}\labell{prop:sect}  Let $\Xx$ be an ep groupoid with a finite good atlas
$\Aa = \{(U_i,G_i,\pi): i\in A\}$, and choose an $\Aa$-compatible cover
 $\Vv = \{|V_I|\}$  of $|X|$ as in Lemma~\ref{le:cov}.  Let $F: \Xx_\Vv\to \Xx$ be the (tame) resolution constructed in Proposition~\ref{prop:weigh}, and
let $\Ee\to \Xx$ be an orbibundle.  Then:\SSS

\NI
{\rm (i)}  Every  section $s$ of  the induced vector bundle
$\rho: E_i\to U_i$ whose support is contained in $U_i^0\sqsubset U_i$ extends to a global section $\si^s$ of the pullback bundle
$F^*(\Ee)\to \Xx_\Vv$.   \SSS

\NI
{\rm (ii)} Let $f:\Xx\to \Ee$ be a section. Suppose that for 
each $i\in A$ 
there is a  space $\Ss_i$ of sections of the vector bundle $E_0|_{U_i}\to  U_i$ such that $f+s: U_i\to E_0$ is transverse to the zero section over $U_i$ for sufficiently small
generic  $s\in \Ss_i$.  Then there is a corresponding
 space $\Ss$ of sections of the pullback bundle $\phi^*(\Ee)\to \Xx_\Vv$ such that
 $F^*(f)+\si:\Xx_\Vv\to F^*(\Ee)$ is transverse to the zero section for sufficiently small
 generic $\si\in \Ss$.
\end{prop}
\begin{proof} Suppose given a section $s:U_i\to E_i$.  If $i\in J$, we define $\si^s(x)$ on the object $x=(\de_{k-1},\dots,\de_1)\in \HHat V_J$ to be $s(\HHat\pi_{iJ}(x))$.  Otherwise, $|V_J|$ is disjoint from the support $|U_i^0|$ of $s$ by construction of $\Vv$, and we set $\si^s(x) = 0$.  It follows immediately from the definitions that $\si^s$ is compatible with the morphisms in $\Xx_\Vv$ and 
so extends to a functor $\Xx_\Vv\to \phi^*(\Ee)$. This proves (i).

To prove (ii), choose a smooth partition of unity $\be$ on $\Xx_\Vv$, which exists by Lemma~\ref{le:part}. 
Then define $\Ss$ to be the vector space generated by the sections $\be\si, s\in\Ss_i.$
 It is easy to check that it has the required properties.
\end{proof}

The previous result used the existence of a resolution to construct multisections. This procedure can be reversed.  The following argument applies to any  $\Xx$ that acts effectively on some vector bundle $E_0\to X_0$.\footnote
{
Whether one can always find suitable $\Ee$ is closely related to the
presentation problem discussed in~\cite{HM}.}  
\MS

\NI
{\bf Proof II of Proposition~\ref{prop:resol} for effective groupoids $\Xx$.} 
Let $\Ee\to \Xx$  be the tangent bundle of $\Xx$; cf.
 Example~\ref{ex:bund}(i). Then, because $\Xx$ is effective, the set
$E_0^{sing}$ is nowhere dense. (In fact, $E_0$ may be triangulated in such a way that $E_0^{sing}$ is a union of simplices 
of codimension $\ge 2$.)
Each point $x_0\in X_0$ has a neighborhood $U$ with a smooth compactly supported section
$s_U: U\to E_0|_U$ that is nonzero at $x_0$ and satisfies the condition 
in Proposition~\ref{prop:constr}(ii), i.e.  no open subset of $s_U(U)$ 
 is contained in $E_0^{sing}$. 
 If $U$ is part of a local chart $(U,G,\pi)$ for $\Xx$, then, as explained in~\cite{HWZ},
 we may extend $s_U$ 
 to a multisection  $\Ll_S$ as follows.  Over $X_0\setminus U$,
 $\Ll_S$ is the characteristic function of the zero section, i.e. $\Ll_S(x,e) = 0$ 
for $e\neq 0$ and  $\Ll_S(x,0) = 1$, while if $x\in U$,
 $$
% \begin{array}{ccll}
 \Ll_S(x,e) = \sum_{g\in G\,:\, s_U(x)\cdot g = (x,e)}\;1/|G|.
  %&\mbox { if } x\in U\\
% & = & 0&\mbox{ otherwise}.
% \end{array}
 $$
% and extend over the rest of $X_0$ by setting $\Ll_S$ equal to the characteristic function of the zero section, i.e. $\Ll_S(x,e) = 0$ 
%for $e\neq 0$ and  $\Ll_S(x,0) = 1$.
 Now choose a good atlas $\Aa = \{(U_i, G_i,\pi_i): i\in A\}$ for $\Xx$
 and a subordinate partition of unity $\be_i$.  Choose for each $i$
  and $g_i\in G_i$ a (non compactly supported) section $s_i:U_i\to E_0$ whose intersection with the nonsmooth set $E_0^{sing}$
  has no interior points. Then the sections
  $r\be_is_i$  
  for any $r\in \R$ have the same property over the support of $\be_i$, at least, since
 $E_0^{sing}$ consists of lines $(x,\la e), \la \in \R$.  Now construct $\Ll_{S_i}$ as above from the section $r_i\be_i s_i$, where $r_i\in \R$ and then define $\Ll_S: = \sum_i\Ll_{S_i}$.  Consider the corresponding groupoid $\Yy_S$ as defined in
Proposition~\ref{prop:constr}.  This will be nonsingular 
for generic choice of the constants $r_i$ and hence a wnb groupoid. 
We can tame it by 
Lemma~\ref{le:tame}.
\QED\MS

\begin{rmk}\labell{rmk:HWZ}\rm  We end with some remarks about
the {\bf infinite dimensional case.}
There are quite a few places in the above arguments where we used the local compactness of $|X|$.  For example we required that layered coverings
$F:\Bb'\to \Bb$ give rise to proper maps $|F|_{\Hsm}:|B'|_{\Hsm}\to |B|_{\Hsm}$, and used this to obtain the even covering property
of Lemmas~\ref{le:brstr} and \ref{le:brstr1}. 
In the infinite dimensional context 
the notion of layered covering must be formulated in such a way
 that these lemmas hold.

Both constructions for the resolution of an ep  groupoid work
in quite general contexts.  For example, 
the construction of $\Xx_\Vv$ works for groupoids in any category 
in which Lemma~\ref{le:cov}  holds.
  Similarly, the above construction of $\Yy_S$
works as long as there are enough local sections of $\Ee\to \Xx$.
In particular we need partitions of unity.

Hofer--Wysocki--Zehnder~\cite{HWZ} define a notion of properness 
that yields a concept of
 ep polyfold groupoid $\Xx$ that has all expected properties.
 In particular, these groupoids admit smooth partitions of unity
 since they are modelled on $M$-polyfolds built using Hilbert rather than Banach  spaces.
 The  Fredholm theory developed in \cite{HWZ}  implies that
a vector bundle $\Ee\to \Xx$ equipped with a Fredholm section $f$ has a good class of Fredholm multisections $\Ll_S: = f+s$ that perturb $f$ and meet the zero section transversally.  
 Since the kernel of the linearization of the operators $f+s$
 has a well defined orientation,\footnote
 {
 Note that in this paper we have assumed that all objects 
 (ep groupoids, branched manifolds, bundles) are oriented. We then  give the zero set of a 
 multisection the induced orientation.  However, 
 when constructing the Euler class
  it is not necessary to orient the ambient
 orbibundle $\Ee\to \Xx$ provided that one consider a class of multisections whose intersections with the zero set
 carry a natural orientation; cf. the definition of $G$-moduli problem in~\cite[Def.~2.1]{CRS}.}
  this intersection can be given the structure of a finite dimensional (oriented) wnb groupoid $\Zz$.  Since any two such multisections are cobordant in the sense that they extend to the pullback bundle over $\Xx\times \Ii$,
 the cobordism class 
 of $\Zz$ is independent of choices, as is the fundamental class $[Z]$
 defined in Proposition~\ref{prop:fclass}.  
\end{rmk}

%%%%%%%%%%%%%%%%%%%%%%%%%%%%%%%%%%%%%%%%%%%%%%%%%%%%%%%%%%%%%
\subsection{Branched manifolds and the Euler class.}\labell{ss:euler}
%%%%%%%%%%%%%%%%%%%%%%%%%%%%%%%%%%%%%%%%%%%%%%%%%%%%%%%%%%%%%

We now discuss the relation of our work on the Euler class to that of Satake and Haefliger and also to the paper~\cite{CRS} of
Cieliebak--Mundet i Riera--Salamon. 

First, we sketch a proof that the homological Euler 
class of Definition~\ref{def:euler1} is Poincar\'e dual to the usual Euler class for orbibundles.
As Haefliger 
points out in~\cite{Hae2}, one can define cohomology characteristic classes 
for orbifolds by adapting the usual constructions for manifolds.  For example, if $\un E\to \un X$ is an oriented $d$-dimensional orbibundle, choose the representing functor
$\Ee\to \Xx$ so that the vector bundle $\rho: E_0\to X_0$ is trivial over each component of $X_0$.  Each trivialization of this bundle 
defines a functor $F:\Xx\to \GGL$, where $\GGL$ is the 
topological category with one object and morphisms $GL(d,\R)^+$ corresponding to the group $GL(d,\R)^+$ of matrices of positive determinant.  Then  the classifying space $B\GGL$ is a model for the classifying space  $BGL(d,\R)^+$ of oriented $d$-dimensional bundles, and so carries a universal bundle.
The Euler class $\eps(\Ee)$ of $\Ee\to \Xx$ is defined to be the pullback 
by $BF: B\Xx\to B\GGL$ of the Euler class of this universal bundle.
Since the projection $B\Xx\to |X|$ induces an isomorphism on rational homology we may equally well think of $\eps(\Ee)$ as an element in
$H^d(|X|;\Q)$.  As such, it depends only on $\un E\to \un X$ and so may also be called $\eps(\un E)$.

We claim that the homology Euler class $e(\Ee)\in H_{n-d}(|X|;\Q)$ of
Definition~\ref{def:euler1} is Poincar\'e dual to $\eps(\Ee)$. One way to prove this is as follows. In what follows we assume for simplicity that $\Xx$ is effective.    Note that $\eps(\Ee)$ may be represented
on the orbifold $|X|$ 
in terms of de Rham theory  by $i_0^*(\tau)$,
 where $i:|X|\to |E|$ is the inclusion of the zero section and 
 $\tau\in H^d_c(|E|)$ is a compactly supported smooth form that represents the Thom class.\footnote{Thus, for each 
 $x\in X_0$, $\tau$ pulls back under the map $E_x\to |E|$ to a generator of the image of 
 $H^d(E_x,E_x\setminus\{0\};\Z)$ in  $H^d(E_x,E_x\setminus\{0\};\R)$.} 
 Therefore it suffices to show that for each smooth 
$(n-d)$-form $\be$ on $|X|$,
$$
\int_{|X|} i^*(\tau)\wedge \be = \int_Z \be,
$$
where the tame branched manifold $(\un Z,\La_Z)\to \un X$ represents
$e(\Ee)$ and the right hand integral is defined by equation~
(\ref{eq:int}).  This holds by adapting standard arguments from the smooth case.  For example, we may choose 
 a smooth triangulation of $|X|$ so that the singular points
 lie in the codimension $2$ skeleton, and then choose 
 the multisection $\Ll_S$ so that its zero set $Z_0\subset X_0$ 
is transverse to this triangulation.  
Then $|X|^{sm}$ and $|Z|^{sm}: = |Z|\cap |X|^{sm}$ are
both pseudocycles in $|X|$.  
Then the left hand integral above  equals the integral of $\tau\wedge\rho^*(\be)$ over 
$i(|X|^{sm})$.  
But the pseudocycle $i(|X|^{sm})$ is cobordant to 
the weighted pseudocycle given by the image of 
the multisection $\Ll_S$. (This is just the image of the fundamental class of $\Yy_S$ under $|\Si|_S$.) We now perturb the latter cycle straightening it out near the zero set $|Z|^{sm}$ so that its intersection with a neighborhood 
$U\subset |E|$ of  $i(|X|^{sm})$ is precisely 
 $|\rho|^{-1}(|Z|^{sm})\to |Z|^{sm}$ (with the obvious rational weights.)  We may then choose $\tau$ so that it vanishes outside $U$.  
Then
$$
\int_{|X|} i^*(\tau)\wedge \be = 
\int_{i(|X|^{sm})} \tau\wedge \rho^*\be = \int_C\tau\wedge\rho^*\be =
\int_Z \be.
$$
  
Finally, we compare our approach to that of
Cieliebak--Mundet i Riera--Salamon.
Although they work in 
the category of  Hilbert manifolds, we shall restrict attention here to the finite dimensional case.
They consider an orbifold to be a quotient $M/G$, where $G$ is a compact Lie group acting on a smooth manifold $M$ with finite stabilizers.
Hence, for them an orbibundle is a $G$-equivariant bundle $\rho: E\to M$, where again $G$ acts on $E$ with finite stabilizers.  Their Proposition 2.7
shows that there is a homology  Euler class $\chi(E)$ for such bundles that lies  
 in the equivariant homology group $H_{n-d}^G(M;\Q)$ and has
all the standard properties
 of such a characteristic class, such as naturality, the expected relation to the Thom isomorphism,
and so on.

We claim that the Euler class described in
 Definition~\ref{def:euler1} is the same as theirs.  
 To see this, first note the following  well known 
 lemma; cf. \cite{HM}.
 
 \begin{lemma} Every effective orbifold  may be identified with 
 a quotient $M/G$, where $G$ is a compact Lie group acting on 
 a smooth finite dimensional manifold $M$ with finite stabilizers.
 \end{lemma}
\begin{proof}[Sketch of Proof.]
Given such a quotient 
one can define a corresponding ep groupoid $\Xx$ by taking a complete set of local slices for the $G$ action;
  cf. the discussion at the beginning of
\S3 in Moerdijk~\cite{Moe}. In the other direction, given an ep groupoid $\Xx$ one takes $M$ to be the orthonormal frame bundle of $|X|$ with respect to some Riemannian metric on its tangent bundle.
The group $G$ is either ${\rm O}(d)$ or ${\rm SO}(d)$, where
$d=\dim X$.  
\end{proof}

Similarly, every orbibundle $\Ee\to \Xx$ can be identified with
a $G$-equivariant bundle $\rho: E\to M$.  
The equivariant homology $H_*^G(M;\Q)$ is, by definition, the homology of the
homotopy quotient $EG\times_G M$ (where $EG\to GB$ is the universal $G$-bundle, i.e. $G$ acts freely on $EG$).  Because $G$ acts with finite stabilizers, the natural quotient map\footnote
{
If the ep groupoid $\Xx$ is an orbifold structure on $M/G$, then
 $EG\times_G M$ can be identified with the  
 classifying space $B\Xx$ (denoted $|\Xx_{\bullet}|$ in~\cite{Moe}),
 and the projection $EG\times_G M\to M/G$ can  be identified with
 the natural map $B\Xx\to |X|$.}
$EG\times_G M\to M/G$ induces an isomorphism on {\it rational} homology.  Hence, in this case, their Euler class may be considered to lie in $H_*(M/G;\Q) = H_*(|X|;\Q)$.

In~\cite[\S10]{CRS} the Euler class is constructed as the zero set of a multivalued section of the orbibundle $\rho:E\to M$, much as
above.  However,  the definitions in \cite{CRS} are all somewhat different from ours.
For example the authors do not give an abstract definition of a weighted branched manifold but rather think of it as a subobject of some high dimensional manifold $M$ on which the compact Lie group  $G$ acts.
They also treat orientations a little differently, in that they do not 
assume the local branches are consistently oriented but rather incorporate the orientation into a signed weighting function on the associated oriented Grassmanian bundle; cf. their Definitions~9.1 and~9.11.  As pointed out in their Remark~9.17, their Euler class can be defined in the slightly more restrictive setting in which the orientation is given by
a consistent orientation of the branches.  Hence below we shall 
assume this, since that is the approach taken here. 
A third difference is that
their Definition~9.1 describes the analog of a weighted branched groupoid, i.e. they do not restrict to the nonsingular case as we did above.\footnote
{
We made this restriction in order to simplify the exposition; 
certain technical details become harder to describe if one allows the local branches over $N_p$ to have completely arbitrary orbifold structures.
This problem is avoided in~\cite{CRS} because in this case the local orbifold structures are determined by the global $G$ action.} 

However to check that the two Euler classes are the same, it suffices to check that 
%their branched manifolds do have wnb groupoid structures and that 
the multivalued sections used in their definition
can be described by functors $\Ll_S$ as 
in Definition~\ref{def:multi}.   
Here is the definition of a multivalued section given in \cite[Def.~10.1]{CRS} in the oriented case.

\begin{defn}\label{def:CRS}  Consider an oriented finite dimensional locally trivial bundle $\rho: E\to M$ over an oriented smooth finite dimensional manifold $M$.  Assume that a compact oriented Lie group $G$ acts 
smoothly, preserving orientation and with finite stabilizers on $E$ and $M$, and that $\rho$ is $G$-equivariant.    Then
a {\bf multivalued section} of $E$ is a function
$$
\si:E\to \Q\cap [0,\infty)
$$
such that \SSS

\NI
{\bf (Equivariance)} $\si(g^*x,g^*e) = \si(x,e)$ for all $
x\in M, e\in E_x, g\in G$,\SSS

\NI
{\bf (Local Structure)} for each $x_0\in M$ there is an open neighborhood $U$ of $x_0$ and finitely many smooth sections $s_1,\dots, s_m:U\to E$ with weights $m_i\in \Q\cap(0,\infty)$ such that $$
\sum m_i = 1,\qquad
\si(x,e) = \sum_{s_i(x) = e} m_i\;\; \mbox{for all }x\in U,
$$
where by convention the sum over the empty set is $0$.
\end{defn}

To see that every such multisection can be described in terms of a functor $\Ll_S$ we argue as follows.  Choose a locally finite covering of $M$ by sets $U^\al$ that have local section structures in the
 sense of Definition~\ref{def:multi}.
By~\cite{CRS} Proposition 9.8(i), the sets $U^\al$ and the sections $s_i$ are locally $G$-invariant.  In other words, we may assume that each $U^\al$ is the image of a finite to one map  $V^\al\times \Nn^\al_G\to U^\al$, where $V^\al\subset M$ is a local slice for the $G$ action and  $\Nn^\al_G$ is a neighborhood of the identity in $G$.  We may use the
local slices $V^\al$, with their induced orientations,  to build an ep groupoid $\Xx$  representing $M/G$
with objects $X_0: = \sqcup V^\al$ and morphisms induced by the $G$ action.  There is a similar groupoid $\Ee$ representing $E/G$ with objects $\rho^{-1}(V^\al) = \cup_{x\in V^\al} E_x$.  Then, because $\si$ is $G$-invariant,
we may define a functor 
$\Ll_S: \Ee\to \Qq^{\ge0}$ by setting
$$
\Ll_S(x,e) = \si(x,e),\quad x\in V^\al, e\in E_x.
$$
This clearly satisfies the conditions of Definition~\ref{def:multi}.
It is now straightforward to check that the definition of Euler 
class in \cite{CRS} is consistent with the one given here.


\begin{thebibliography}{CCCCC}



\bibitem{CHu}  B. Chen and S. Hu, A de Rham model for Chen--Ruan   
       cohomology ring of abelian orbifolds,
       SG/0408265
       
\bibitem{CR}  W. Chen and Y. Ruan, A new cohomology theory of 
    orbifold, AG/0004129, {\it Comm. Math. Phys.} {\bf 248} (2004),    
    1--31.

\bibitem{Ch}  W. Chen, Pseudoholomophic curves in 
$4$-orbifolds and some applications, SG/0410608, to appear in Fields Institute Communications.

\bibitem{CRS} K. Cieliebak, I. Mundet i Riera and D. Salamon,
Equivariant moduli problems, branched manifolds and the Euler class, 
{\it Topology} {\bf 42} (2003) 641--700.

\bibitem{Hae1} A. Haefliger, Homotopy and Integrability, {\it Manifolds, Amsterdam, 1970}, Springer Lecture Notes in Math, {\bf 197}, (1971), 133-163.

\bibitem{Hae2} A. Haefliger, Holonomie et Classifiants, {\it Ast\'erisque} {\bf 116} (1984), 70-97.

\bibitem{Hae3} A. Haefliger, Groupoids and Foliations, {\it Contemporary Math} {\bf 282} (2001), 83--100.

\bibitem{HM} A. Henriques and D. Metzler, Presentations of 
        Noneffective Orbifolds, AT/0302182, {\it Trans. Amer. Math. Soc.} {\bf 356} (2004), 2481-99.

\bibitem{H} H. Hofer, A general Fredholm theory and applications,
SG/0509366. 

\bibitem{HWZ} H. Hofer, C. Wysocki, and E. Zehnder, 
       {\it Polyfolds and Fredholm Theory}, Parts I and II, 
       preprint (2005), FA/0612604.

\bibitem{Ler} E. Lerman, Orbifolds as a localization of the $2$-    
        category of groupoids, DG/0608396.

\bibitem{LiuT}
Gang Liu and Gang Tian, Floer homology and Arnold
conjecture, {\it Journ. Diff. Geom.}, {\bf 49} (1998), 1--74.

\bibitem{LuT} Guangcun Lu and Gang Tian, Constructing virtual 
      Euler cycles and classes, preprint (Dec. 2005). 

\bibitem{Mcv} D. McDuff, The virtual moduli cycle, 
{\it Amer. Math. Soc. Transl.} (2) {\bf 196} (1999), 73 -- 102.


\bibitem{Moe} I. Moerdijk, Orbifolds as groupoids, an introduction.    
     DG/0203100, in {\it Orbifolds in Mathematics and Physics}, ed   
     Adem, {\it Contemp Math} {\bf 310}, AMS (2002), 205--222.

\bibitem{MM} I. Moerdijk and Mr\v cun, {\it Introduction to Foliations and Lie Groupoids}, Cambridge Studies {\bf 91} (2003), CUP.

\bibitem{MPr} I. Moerdijk and D.A. Pronk, Simplicial cohomology of     
     orbifolds, {\it Indag. Math.} {\bf 10}  (1999), 269--293.

\bibitem{RobS}  J. Robbin and D. Salamon, A construction of the          
      Deligne--Mumford orbifold, SG/0407090

\bibitem{SalFl}   D. Salamon, Lectures on Floer theory, 
       {\it Proceedings of the IAS/Park
       City Summer Institute, 1997}, ed Eliashberg and Traynor, 
       Amer. Math. Soc., Providence, RI. (1999), 143--229.

\bibitem{Sat} I. Satake, On a generalization of the notion of     
       manifold, {\it Proc. Nat. Acad. Sci.} {\bf 42} (1956), 
       359--363.

\bibitem{Sat2} I. Satake, The Gauss--Bonnet Theorem for $V$-manifolds, {\it J. Math. Soc. Japan} {\bf 9} (1957), 464--492.

\bibitem{Sch} M. Schwarz, Equivalences for Morse homology, in {\it 
       Geometry and Topology in Dynamics} ed M. Barge, K. Kuperberg, 
       Contemporary Mathematics {\bf  246}, Amer. Math. Soc. (1999),   
       197--216.  

\bibitem{Z}  A. Zinger, Pseudocycles and Integral Homology, AT/0605535.
\end{thebibliography}
 \end{document}